\documentclass{elsart}
\usepackage{graphicx}
\usepackage{amssymb}
\usepackage{lineno}
\usepackage[latin1]{inputenc}
\usepackage[francais,english]{babel}
\usepackage{amsmath}
\usepackage{array}
\usepackage{amsfonts}
\usepackage{amssymb}
\usepackage{latexsym}
\usepackage{mathrsfs}
\usepackage{pstricks}
\usepackage{algorithm}
\usepackage{algorithmic}  
\usepackage{psfrag}
\usepackage[center]{caption}
\newcommand{\disp}{\displaystyle}
\newcommand{\eqdef}{\ensuremath{\stackrel{\mathrm{def}}{=}}}
\newcommand{\ud}{\frac{1}{2}}
\newcommand{\h}{\mathfrak{h}}
\newcommand{\Th}{\mathcal{T}_\mathfrak{h}}
\newcommand{\A}{\mathcal{A}}
\newcommand{\B}{\mathcal{B}}
\newcommand{\I}{\mathcal{I}}
\newcommand{\W}{\mathcal{W}}
\newcommand{\WB}{\W\cup\B}
\newcommand{\Fi}{\mathcal{F}^{\rm{i}}_{\h}}
\newcommand{\F}{\mathcal{F}_{\h}}
\newcommand{\FWB}{\mathcal{F}^{\mathcal{WB}}_{\h}}
\newcommand{\FI}[1]{\mathcal{F}_{\h}^{{\rm{w}},#1}}
\newcommand{\FII}[1]{\mathcal{F}_{\h}^{{\rm{d}},#1}}
\def \[{ \mbox{[ \hspace{-.5em}[}}
\def \]{ \mbox{] \hspace{-.5em}]}}
\newcommand{\Id}[1]{\I^{{\rm{d}},#1}}
\newcommand{\Iw}[1]{\I^{{\rm{w}},#1}}
\newcommand{\Ihd}[1]{\I_{\h}^{{\rm{d}},#1}}
\newcommand{\Ihw}[1]{\I_{\h}^{{\rm{w}},#1}}
\newcommand{\M}{{\mathbb{M}}}

\newcommand{\DM}{\Delta{\M}}
\newcommand{\Ih}{\I_{\h}}
\newcommand{\vr}{v_{\rm{r}}}

\newtheorem{theo}{Property}
\begin{document}

\begin{frontmatter}
\title{Mass conservative BDF-discontinuous Galerkin/explicit finite volume schemes for coupling subsurface and overland flows}
%\author{P. Sochala, A. Ern, S. Piperno\corauthref{cor1}.}
\author{P. Sochala, A. Ern, S. Piperno}
\address{Universit\'e Paris-Est, CERMICS, Ecole des Ponts, Champs sur Marne, 77455 Marne la Vall\'ee Cedex 2, France}
%\corauth[cor1]{corresponding author: serge.piperno@cermics.enpc.fr}
\begin{abstract}
Robust and accurate schemes are designed to simulate the
coupling between subsurface and overland flows. The coupling
conditions at the interface enforce the continuity of both
the normal flux and the pressure. Richards' equation
governing the subsurface flow is discretized using a
Backward Differentiation Formula and a symmetric interior
penalty Discontinuous Galerkin method. The kinematic wave
equation governing the overland flow is discretized using a
Godunov scheme. Both schemes individually are mass
conservative and can be used within single-step or multi-step
coupling algorithms that ensure overall mass conservation
owing to a specific design of the interface fluxes in the
multi-step case. Numerical results are presented to
illustrate the performances of the proposed algorithms.  
\end{abstract}
\end{frontmatter}

\section{Introduction}
The interactions of subsurface and overland flows are an
important ingredient for a comprehensive understanding of
hydrology processes. While there is an extensive bulk of
literature devoted to the numerical study of water flows in
single-phase and variably saturated porous media, the issue
of coupling such flows with surface flows generated by
rivers, tides or floods has received less attention. One of
the most popular ways to couple Darcy and Stokes flows is
through the well-known Beavers--Joseph--Saffman condition
\cite{Beavers67,Saffman71,Mikelic00}. This condition was
used for instance in \cite{Discacciati02,Miglio02} in the
mathematical and numerical study of the coupling of Darcy
flow with a three-dimensional non-hydrostatic shallow-water
model. Another approach used in numerical hydrology (see
among others \cite{VanderKwaak01}) considers discontinuous
pressures at the interface and evaluates an interface flux
as the pressure difference, modulated by a multiplicative
exchange coefficient depending on the soil relative permeability. A
third approach consists of assuming both normal flux and
pressure continuity. This means that the hydraulic head of
the subsurface flow matches the depth of the overland flow
at the interface, while the normal ground flow velocity is
used as a source term in the governing equation of the
overland flow. Examples of studies based on this approach
include coupling one-dimensional surface flow with vertical
soil columns \cite{Singh98}, coupling the two-dimensional
Richards' equation with a one-dimensional kinematic or
diffusive wave approximation for the overland flow
\cite{Kollet05,Beaugendre06}, and coupling the
two-dimensional Darcy's equation with one-dimensional
shallow-water equations \cite{Dawson06} .

In the present work, we assume that the subsurface flow
occurs in a variably saturated porous medium and that this
flow can be described by Richards' equation, entailing in
particular that there are no trapped air pockets in the
soil; otherwise more general multi-phase models should be
used \cite{Bastian99}. Furthermore the kinematic wave
approximation is used to describe the overland flow. This
choice is solely made for ease of exposition and more
general shallow water models can also be used. Concerning
the coupling conditions, we adopt the third approach
described above, namely enforcing the continuity of both
normal flux and pressure at the interface. These coupling
conditions are generally valid when the overland flow is
mainly produced by exfiltration from the soil, so that
normal flux and pressure equilibrium can be expected to hold
at all times. A different situation, which falls beyond the
present scope, is that of a runon surface wave rapidly
propagating over a dry soil. 

Many methods can be employed to discretize in space
Richards' equation, namely finite differences
\cite{Woodward,Celia90}, finite volumes (FV)
\cite{Manzini04}, finite elements (FE)
\cite{Celia90,kavanagh-nonsmooth} or mixed finite elements
(MFE) \cite{Knabner02,Bause04}. These methods are generally
combined with an implicit Euler time scheme. An
alternative approach for space discretization is to use a
discontinuous Galerkin (DG) method. Advantages offered by DG
methods include local (elementwise) conservation (as FV and
MFE), high-order accuracy (as FE and MFE) and flexibility in
the use of non-matching meshes (as FV), in particular within
multi-physics and multi-domain approaches. Various forms of
DG methods can be used for Richards' equation and more
generally for two-phase flows in porous media. Examples
include the so-called Local Discontinuous Galerkin method
\cite{Fagherazzi04,Bastian07}  and  the non-symmetric or the
symmetric interior penalty DG method
\cite{Klieber06,Bastian02,Bastian04}. In the present work,
we choose the symmetric interior penalty DG method (in short
SIPG), because it preserves the natural symmetry in the
discrete diffusion operator. Regarding time discretization,
the common approach when working with DG methods is to
employ Runge-Kutta (RK) explicit schemes
\cite{Cockburn_Shu98} or diagonally implicit ones
\cite{Bastian02}. Here, we propose instead to use a backward
differentiation formula (BDF). We think that this approach
offers several advantages, such as high-order accuracy in
the time discretization, circumventing the CFL condition
which can be very restrictive for explicit schemes when
diffusion processes are dominant, and in general higher
computational efficiency than implicit RK schemes for
problems where the nonlinear solver is expensive. Typically,
if piecewise polynomials of degree $p$ are used in the DG
method, a BDF of order $(p+1)$ can be employed. 

The main objective of this work is to design robust and
accurate schemes for coupling subsurface and overland
flows. While Richards' equation is discretized by a BDF-SIPG
method, the overland flow governing equation is discretized
by a Godunov scheme and advanced in time with a different
time step if the overland flow time scale is quite different
from the subsurface flow time scale. Two important issues
are addressed in the design of our coupling
algorithms: 1) satisfy as accurately as possible
the coupling conditions which impose certain specific
inequality and equality constraints on the pressures and
normal fluxes, similarly to the boundary conditions
encountered in Signorini problems, and 2) ensure
overall mass conservation for the whole system consisting of
subsurface and overland flow. This point deserves some
particular attention. Indeed, although mass conservative
schemes are used for both subsurface and overland flows, the
interface flux must be chosen appropriately when working
with multi-step methods such as BDFs. For simplicity, we
will discuss in detail only the design of the interface flux
for the second-order BDF. Finally, although the material
will be presented in a 2D/1D setting (that is, a
two-dimensional subsurface flow coupled to a one-dimensional
overland flow), the results extend naturally to the 3D/2D
setting. In particular, the wet part of the interface is not
tracked directly, but is determined at each time step by a
cell-oriented procedure within an iterative loop that solves
consecutively the overland and subsurface flow governing
equations. 

This paper is organized as follows. In
Section~\ref{sec:model}, we present the physical problem. In 
Section~\ref{sec:discretization}, we describe the time and
space discretization of the model problem and design the
coupling algorithms for both first-order and second-order
BDFs. Finally, in Section~\ref{sec:results}, we present
numerical results assessing the performance of the proposed
algorithms on three test cases. 

\section{Model problem}
\label{sec:model}
\subsection{The setting}
\label{sec:setting}
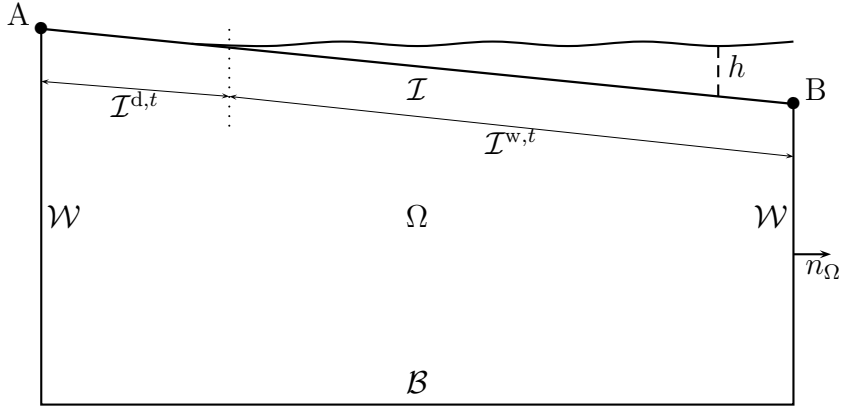
\begin{figure}[htb]
\centering
\begin{pspicture}(-1,0)(11,5.5)
%--- Boite ---------------------------
\psline(0,5)(0,0)(10,0)(10,4)
\psline[linewidth=0.3mm](0,5)(10,4)
\pscurve(2,4.80)(3,4.77)(4,4.83)(5,4.77)(6,4.83)(7,4.77)(8,4.83)(9,4.77)(10,4.83)
\psline[linestyle=dashed](9,4.1)(9,4.77)
\psline{->}(10,2)(10.5,2)
\psline[linewidth=0.1mm]{<->}(0,4.3)(2.5,4.1)
\psline[linewidth=0.1mm]{<->}(2.5,4.1)(10,3.3)
\psline[linestyle=dotted](2.5,5)(2.5,3.7)
%--- Les indic des points -------------
\rput(5,4.2){$\I$}
\rput(1.25,4){$\Id{t}$}
\rput(6.25,3.5){$\Iw{t}$}
\rput(0,5){$\bullet$}
\rput(10,4){$\bullet$}
\rput(-0.3,5.2){$\rm{A}$}
\rput(10.3,4.2){$\rm{B}$}
\rput(5,2.5){$\Omega$}
\rput(0.3,2.5){$\W$}
\rput(9.7,2.5){$\W$}
\rput(5,0.3){$\B$}
\rput(9.25,4.5){$h$}
\rput(10.4,1.8){$n_{\Omega}$}
\end{pspicture}
\vspace{0.3cm}
\caption{Schematic of the computational domain with basic notation.\label{Fig:Couplage1}} 
\end{figure}

Let $\Omega \subset \mathbb{R}^2$ denote the bounded
subsurface flow domain with outward normal unit vector
$n_{\Omega}$. The boundary of $\Omega$ is divided into three
parts (see Figure~\ref{Fig:Couplage1}): $\I$ is the upper
part of the boundary where overland flow can occur, $\W$ are
lateral walls and $\B$ represents the lower part of the
boundary. At any time $t$, the set $\I$ is divided into
``wet'' and ``dry'' parts $\Id{t}\cup \Iw{t}$, with 
\begin{equation}
\Iw{t}\eqdef \{x\in \I;\ h(x,t)>0\},\quad \Id{t} \eqdef
\{x\in \I;\ h(x,t)=0\}, 
\end{equation}
where $h$ is the depth of the overland flow. Observe that
the above partition of $\I$ is time-dependent. 

\subsection{Subsurface flow}
\label{sec:subsurfaceflow}
The soil is modeled as a non-deformable porous medium in
which the pores can contain both water and air (unsaturated
zone) or only water (saturated zone). We assume that water
is incompressible and that air pressure does not affect
the flow. The water conservation equation takes the form 
\begin{equation}\label{Eq:Richards_1}
\partial_t[\theta(\psi)] + \nabla \cdot v(\psi) = f,
\end{equation}
where $\partial_t$ denotes partial time-derivative, $\psi$
is the hydraulic head ($m$), $\theta(\psi)$ the 
volumetric water content (dimensionless), $v(\psi)$ the flow
velocity ($ms^{-1}$), and $f$ a volumetric water source or
sink ($s^{-1}$). In the sequel, we assume that there are no
volumetric sources or sinks, so that $f=0$. The flow
velocity depends on the hydraulic head through the
generalized Darcy law 
\begin{equation}
\label{Eq:Richards_2}
v(\psi) = -K(\psi) \nabla (\psi + z ),
\end{equation}
where $K(\psi)$ is the (possibly tensor-valued) hydraulic
conductivity ($ms^{-1}$) and $z$ the vertical coordinate
($m$). Substituting (\ref{Eq:Richards_2}) into
(\ref{Eq:Richards_1}), Richards' equation is obtained in the
form \cite{Ri31}
\begin{equation}
\label{Eq:Richards_3}
\partial_t[\theta(\psi)] - \nabla \cdot ( K(\psi) \nabla (\psi + z )) = 0.\
\end{equation}
Given at each time $t\in [0,T]$, where $T$ is the total
simulation time, the partition $\{\Iw{t},\Id{t}\}$ of $\I$,
a Dirichlet datum $\omega_{\psi}$ defined on $\Iw{t}$ and a
Neumann datum $\omega_{v}$ defined on $\Id{t}$, the
subsurface flow is governed by  
\begin{equation}
\label{Eq:Couplage1}
\left\{\begin{array}{ll}
\partial_t[\theta(\psi)] + \nabla \cdot v(\psi) = 0
&\quad\mbox{in}\ \Omega \times [0,T],\\ 
v(\psi)= -K(\psi)\nabla (\psi+z)&\quad\mbox{in}\ \Omega \times [0,T],\\
\psi(\cdot,0)=\psi^0 &\quad\mbox{in}\ \Omega,\\
v(\psi)\cdot n_{\Omega} = v_N &\quad\mbox{on}\ (\W \cup \B) \times [0,T],\\
v(\psi)\cdot n_{\Omega} = \omega_{v} &\quad\mbox{on}\ \{(x,t),x \in \Id{t}\},\\
\psi = \omega_{\psi} &\quad\mbox{on}\ \{(x,t),x \in \Iw{t}\},
\end{array}\right.
\end{equation}
where $v_N$ is the possibly time-dependent normal flow
velocity prescribed on $\W\cup\B$ and $\psi^0$ the initial
condition. Thus, Richards' equation is a nonlinear
parabolic equation which degenerates into a nonlinear
diffusion equation in the saturated zone where $\theta$ and
$K$ are constant. Examples for the two constitutive laws
$\psi \mapsto \theta(\psi)$ and $\psi \mapsto K(\psi)$,
which are necessary to close the subsurface flow model, are
specified in \S\ref{sec:results}.

\subsection{Overland flow}
\label{sec:overlandflow}
Water surface flows are often modeled by a simplified form
of the free boundary Navier-Stokes equations. Assuming
hydrostatic pressure, negligible vertical velocity gradients
and mild variations of the free surface leads to the
well-known shallow-water equations; see,
e.g. \cite{Gerbeau01} for a derivation of these
equations. Neglecting turbulence effects, the equations
expressing the conservation of mass and momentum reduce to 
\begin{equation}\label{Eq:SV1}
\partial_t h  + \partial_x q = (v(\psi)-\vr)\cdot n_{\Omega},
\end{equation}
\begin{equation}\label{Eq:SV2}
\partial_t q + \partial_x \bigg[ \frac{q^2}{h} +
\frac{gh^2}{2} \bigg] = gh(S-J), 
\end{equation}
where $q$ is the discharge ($m^2s^{-1}$), $v(\psi)\cdot
n_{\Omega}$ the source or sink term ($ms^{-1}$) resulting
from mass transfer between subsurface and overland flows,
$\vr$ the possibly time-dependent prescribed rainfall
intensity ($ms^{-1}$), $g$ the gravity acceleration
($ms^{-2}$), $S$ the possibly space-dependent bottom slope
(dimensionless) and finally $J$ (dimensionless) results from
friction effects. Note that the mass transfer term
$v(\psi)\cdot n_{\Omega}$ in the mass conservation equation
(\ref{Eq:SV1}) involves the subsurface flow velocity
resulting from (\ref{Eq:Couplage1}); infiltration occurs if
$v(\psi)\cdot n_{\Omega} < 0$ whereas exfiltration occurs if
$v(\psi)\cdot n_{\Omega} > 0$. The Manning-Strickler uniform
flow formula is chosen to link $J$ and $q$ and assuming the
flux to be uni-directional from left to right so that
$q\ge0$, this yields 
\begin{equation}
\label{Eq:SV3}
q = \mathcal{K} h^{5/3} J^{1/2},
\end{equation}
where $\mathcal{K}$ is the Strickler coefficient of roughness
($m^{1/3}s^{-1}$).

A common assumption is to neglect inertia and potential
energy effects in (\ref{Eq:SV2}), so that momentum balance
is governed by the equilibrium between slope and friction,
that is 
\begin{equation}
\label{Eq:SV4}
S=J.
\end{equation}
Substituting ($\ref{Eq:SV4}$) into ($\ref{Eq:SV3}$) yields
\begin{equation}
\label{Eq:SV5}
q = \varphi(h,S) \eqdef \mathcal{K} h^{5/3} S^{1/2}.
\end{equation}
Finally, using ($\ref{Eq:SV5}$) in ($\ref{Eq:SV1}$) leads to
the so-called kinematic wave approximation \cite{Ponce77} 
\begin{equation}
\partial_t h + \partial_x \varphi(h,S) = (v(\psi)-\vr)\cdot n_{\Omega}.
\end{equation}
This scalar conservation law is strictly hyperbolic wherever
$h>0$. In the present case, waves travel rightwards and an
upstream boundary condition in $\rm{A}$ (see
Figure~\ref{Fig:Couplage1}) must be set. Let $h^0$ be the
initial condition and let $h_{\rm{A}}$ be the upstream
boundary condition on the surface water depth prescribed at
point $\rm{A}$. Then, the overland flow is governed by 
\begin{equation}
\label{Eq:Couplage2}
\left\{\begin{array}{ll}
\partial_t h +\partial_x \varphi(h,S) = (v(\psi)- \vr)\cdot
n_{\Omega} &\quad\mbox{on}\ \I \times [0,T],\\
h(\cdot,0)=h^0 &\quad\mbox{on}\ \I, \\
h(\rm{A},\cdot)=h_{\rm{A}} &\quad\mbox{at}\ \rm{A}\times [0,T]. \\
\end{array}\right.
\end{equation}

\subsection{Admissible set}
\label{sec:admissibleset}
\begin{figure}[htb]
\centering
\begin{pspicture}(0,0)(10,4)
%--- Boite ---------------------------
\psline[linewidth=0.1mm]{->}(1,0)(9,0)
\psline[linewidth=0.1mm]{->}(5,0)(5,3)
\psline[linewidth=0.3mm](1,0)(5,0)(8,3)
%--- Les indic des points -------------
\rput(5,3.5){$h$}
\rput(9.5,0){$\psi$}
\rput(3,0.5){dry}
\rput(7,2.5){wet}
\end{pspicture}
\vspace{0.3cm}
\caption{The admissible set $\A$.\label{Fig:Couplage2}}
\end{figure}
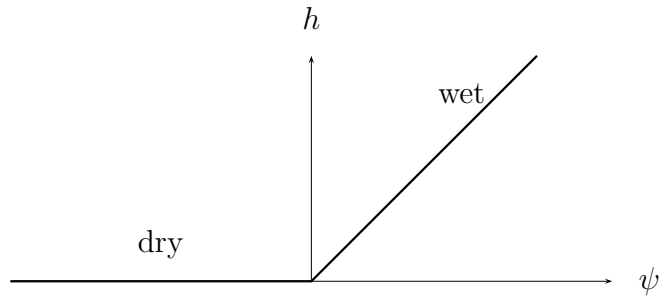

We refer to the quadruplet
$\{\Iw{t},\Id{t},\omega_{\psi},\omega_{v}\}$ as the coupling
variables. The model problem considered hereafter for
coupling subsurface and overland flows consists of finding
functions $(\psi,h)$ and the above coupling variables such
that 
\begin{equation}
\label{Eq:Couplage3}
\left\{\begin{array}{ll}
\psi\ \text{solves}\ (\ref{Eq:Couplage1}) & \mbox{in}\ \Omega \times [0,T],\\
h\ \text{solves}\ (\ref{Eq:Couplage2}) & \mbox{on}\ \I \times [0,T],\\
(\psi,h) \in \A &\mbox{on}\ \I \times [0,T],
\end{array}\right.
\end{equation}
where $\A$ denotes the set of physically admissible states
$\{\psi,h\}$. The admissible set $\A$
(see Figure~\ref{Fig:Couplage2}) has two branches, the branch
$\{h=0\}$ is associated with the dry surface where the soil
hydraulic head is less than or equal to zero corresponding
to unsaturated conditions, while the branch $\{h=\psi\}$ is
associated with the wet surface where the soil is saturated
and the hydraulic head is in hydrostatic equilibrium with
the overland flow pressure. Thus, the admissible set $\A$ is
defined as 
\begin{equation}
\A\eqdef\{ (\psi,h)\in \mathbb{R}^2,\ h=\psi^+\},
\end{equation}
where $\psi^+=\frac{1}{2}(\psi+|\psi|)$ is the positive part
of $\psi$.

We mainly focus here on situations where the overland flow
is produced by exfiltration. Indeed in this situation, a
smooth behavior on the admissible set can be expected. More
drastic situations like runon surface waves on unsaturated
soils can in many cases lead to a departure from the
admissible set especially if the soil is too dry. In these
limit situations, other models can be more suitable: for
instance when infiltration processes are very slow, a model
where surface flow coexists with an unsaturated soil can be
envisaged. 

\section{Discretization}
\label{sec:discretization}
\subsection{Discretization of Richards' equation}
\label{sec:Richardsdiscretization}
Let $\{\Th\}_{\h>0}$ be a shape-regular family of
unstructured meshes of $\Omega$ consisting for simplicity of
affine triangles. The meshes can possess hanging nodes. For
an element $\tau \in \mathcal{T}_{\h}$, let $\partial \tau$
denote its boundary and $n_{\tau}$ its outward unit normal. 
The discontinuous finite element space $V_{\h}$ is defined as
\begin{equation}
\label{FEspace}
V_{\h} \eqdef \{ \phi \in L^2(\Omega),\ \forall \tau
\in \mathcal{T}_{\h},\ \phi{|_{\tau}} \in \mathbb{P}_p(\tau) \},
\end{equation}
where $\mathbb{P}_p(\tau)$ is the set of polynomials of
degree less than or equal to $p$ on an element $\tau$. We
observe that the functions in $V_{\h}$ need not be
continuous. This fact is exploited by selecting basis
functions which are locally supported in a single mesh
element. The set $\F$ of mesh faces is partitioned into $\Fi
\cup \FWB \cup \mathcal{F}^{\I}_{\mathfrak{h}}$ where $\Fi$
is the set of internal faces, $\FWB$ the set of faces
located on $\WB$ and $\mathcal{F}^{\I}_{\mathfrak{h}}$ the
set of faces located on $\I$. For a face $F\in\Fi$, there
are $\tau^+$ and $\tau^-$ in $\mathcal{T}_{\h}$ such that
$F=\partial \tau^+ \cap \partial \tau^-$ and we define the
average operator $\{\}_F$ and the jump operator $\[\]_F$ as
follows: for a function $\xi$ which is possibly two-valued
on $F$, 
$$\{ \xi \}_F \eqdef \frac{1}{2}(\xi^+ + \xi^-)\quad
\text{and}\quad \[ \xi \]_F \eqdef \xi^- - \xi^+,$$
where $\xi^{\pm}=\xi|_{\tau^{\pm}}$. For vector-valued
functions, average and jump operators are defined
componentwise. We define $n_F$ to be the unit normal vector
to $F$ pointing from $\tau^-$ to $\tau^+$. The
arbitrariness in the sign of the jump is irrelevant in the
sequel. 

In the present work, faces on $\I$ can exclusively be
flagged either as dry or as wet, that is, we do not track
the wet/dry interface inside such faces. As a result, the
set $\mathcal{F}^{\I}_{\mathfrak{h}}$ can be further divided
into $\FII{t}$ and $\FI{t}$, where $\FII{t}$ collects the
faces flagged as dry and $\FI{t}$ those flagged as
wet. These two sets of faces induce a partition of $\I$ as
$\Id{t}_{\h}\cup\Iw{t}_{\h}$, where 
$$\Id{t}_{\h} \eqdef \{x\in \I,\ \exists F \in  \FII{t},\ x \in F\}
\text{ and }
  \Iw{t}_{\h} \eqdef \{x\in \I,\ \exists F \in  \FI{t},\ x \in F\}. $$

\subsubsection*{Space discretization}
Let $\psi_{\h}$ be the discrete approximation of $\psi$. The
symmetric interior penalty discontinuous Galerkin method for
Richards' equation can be concisely written as 
\begin{equation}
\label{Richards_un_14}
\forall \tau \in \mathcal{T}_{\h},\quad \forall \phi \in
\mathbb{P}_p(\tau),\quad \int_{\tau} \partial_t[\theta(\psi_{\h})] \phi +
a_{\tau}(\psi_{\h},\psi_{\h},\phi) =
b_{\tau}(\psi_{\h},\phi),
\end{equation}
where for $(\zeta,\psi,\phi)\in V_{\h} \times V_{\h} \times \mathbb{P}_p(\tau),$
\begin{align}
a_{\tau}(\zeta,\psi,\phi)\eqdef & \int_{\tau} K(\zeta)\nabla
\psi\cdot \nabla \phi + \int_{\partial \tau}K(\zeta)\ \nabla
\phi \cdot n_{\tau}\ ( \widehat{\psi}(\psi) - \psi) \nonumber\\ 
& + \int_{\partial \tau} \widehat{u}(\zeta,\psi)\cdot n_{\tau}\ \phi,\\ 
b_{\tau}(\zeta,\phi)\eqdef & \int_{\tau} \nabla 
\cdot(K(\zeta) \nabla z)\ \phi + \tilde{b}_{\tau}(\zeta,\phi).
\end{align}
Here, $\widehat{\psi}(\psi)$ is the numerical flux
associated with the hydraulic head 
$$\forall F \in \F,\ \widehat{\psi}(\psi)|_F \eqdef \left\{\begin{array}{ll}
\{\psi\}_{F} &\mbox{if}\ F \in \Fi, \\[3pt]
\ 0 & \mbox{if}\ F \in \FI{t}, \\[3pt]
\ \psi &\mbox{if}\ F \in \FII{t}\cup\FWB,
\end{array}\right.$$
and $\widehat{u}(\zeta,\psi)$ the numerical flux associated
with $u\eqdef -K(\psi)\nabla\psi$, 
$$\forall F \in \F,\ \widehat{u}(\zeta,\psi)|_F \eqdef \left\{\begin{array}{ll}
-\{K(\zeta)\nabla \psi\}_F + \eta K_s d_{F}^{-1} \[\psi\]_F n_F &\mbox{if}\ F \in \Fi, \\[3pt]
-\ K(\zeta)\nabla \psi + \eta K_s d_{F}^{-1}\psi n_{\Omega} &\mbox{if}\ F \in \FI{t}, \\[3pt]
\hspace{2cm} 0  &\mbox{if}\ F \in \FII{t}\cup\FWB,
\end{array}\right.$$
where $\eta$ is a positive parameter (to be taken larger
than a minimal threshold depending on the shape-regularity
of $\Th$), $K_s$ the hydraulic conductivity at saturation
and $d_{F}$ the diameter of the face $F$ which is defined as
the largest diameter of the triangle(s) of which $F$ is a
face. Observe that for a flow in a porous medium with variable
conductivity (as in variably saturated flows because of the
dependence of the conductivity on the hydraulic head), the
penalty coefficient at a given interface should scale as the
harmonic means of the normal hydraulic conductivity on both
parts of the interface, see \cite{DiPietro08,Ern08}. Here,
the variations of $K$ are sufficiently mild to use simply
the hydraulic conductivity at saturation. Furthermore,
$\tilde{b}_{\tau}(\zeta,\phi)$ collects the parts of the
numerical fluxes on boundary faces which are independent of
$\psi$, namely 
\begin{align*}
\tilde{b}_{\tau}(\zeta,\phi)\eqdef&  
\int_{\partial\tau \cap \FI{t}} \big(-K(\zeta)\ \nabla \phi \cdot n_{\Omega}\ + \eta K_s d_{F}^{-1} \phi \big)\omega_{\psi} \\
&-\int_{\partial\tau \cap \FII{t}} (\omega_v + K(\zeta)\nabla z \cdot n_{\Omega})\phi - \int_{\partial\tau \cap \FWB} (v_N + K(\zeta)\nabla z \cdot n_{\Omega})\phi.
\end{align*}
Summing $a_{\tau}(\zeta,\psi,\phi)$ over all mesh elements
yields the global form 
\begin{align}
&a_{\h}(\zeta,\psi,\phi) = \sum_{\tau \in \Th} \int_{\tau} K(\zeta)\nabla \psi\cdot \nabla \phi \nonumber \\
&\quad - \sum_{F \in \Fi } \int_{F} \Big( \{K(\zeta)\nabla \phi\}\[ \psi \]\cdot n_F + \{K(\zeta)\nabla \psi\}\[ \phi \]\cdot n_F - \eta K_s d_{F}^{-1} \[\psi \] \[\phi \] \Big) \nonumber \\
&\quad - \sum_{F \in \FI{t} } \int_{F} \Big( K(\zeta)\nabla \phi\ \psi \cdot n_F + K(\zeta)\nabla \psi\ \phi \cdot n_F - \eta K_s d_{F}^{-1}\psi\phi \Big).
\end{align}
\noindent The parameter $\eta$ must be chosen large enough
to ensure that the form $a_{\h}$ is coercive, in the sense
that there is $\alpha > 0$ such that for all $\phi \in
V_{\h}$,  
$$a_{\h}(\phi,\phi,\phi) \ge \alpha\Big( \sum_{\tau \in
\Th}\int_{\tau}K(\phi)|\nabla \phi|^2 + \sum_{F \in \Fi}K_s
d_F^{-1}\int_{F} \[\phi\]^2 +  \sum_{F \in \FI{t}}K_s
d_F^{-1}\int_{F} \phi^2 \Big).$$

\subsubsection*{Time discretization}
Let $N_T$ be the total number of time steps and let $\delta
t$ be the time step taken to be constant for the sake of
simplicity and such that $N_T\eqdef T/ \delta t$ is an integer. For any
function of time $\chi$ and for any integer $n \ge 0$,
$\chi^n$ denotes the value taken by $\chi$ at time $n\delta
t$. Furthermore, the time derivative of $\chi$ can be
approximated by a backward differentiation formula
\cite{Quarteroni00} in the form 
\begin{equation}
\label{schematemp_2}
(\partial_t \chi)^n\simeq \sum_{r=0}^{q}
\frac{\alpha_r^q}{\delta t}\chi^{n-r},
\end{equation}
where $q$ is the order of the formula and $\{ \alpha_r^q
\}_{0\le r\le q}$ are suitable coefficients. Using the
approximation (\ref{schematemp_2}) in (\ref{Richards_un_14})
for each $n \in \{1\cdots N_T\}$ leads to 
\begin{align}
&\forall \tau \in \mathcal{T}_{\h},\quad \forall \phi \in
\mathbb{P}_p(\tau),\nonumber\\
&\quad\frac{\alpha_0^q}{\delta t}
\int_{\tau} \theta(\psi_{\h}^n) \phi + 
a_{\tau}(\psi_{\h}^n,\psi_{\h}^n,\phi) =
b_{\tau}(\psi_{\h}^n,\phi) - \sum_{r=1}^{q}
\frac{\alpha_r^q}{\delta t} \int_{\tau} \theta(\psi_{\h}^{n-r}) \phi.
\label{schematemp_3}
\end{align} 
For the first few time steps, a BDF of lower order or a
one-step implicit scheme can be used, for example a
diagonally implicit Runge-Kutta scheme or the Crank--Nicolson
scheme. The former can present the drawback that the last
stage can lead to difficulties where the soil is being
saturated. 

\subsubsection*{Nonlinear iterative solver}
The nonlinear equation (\ref{schematemp_3}) is solved by the
iterative algorithm outlined in
Algorithm~\ref{alg1}. The discrete functions
$\big\{\psi_{\h}^{n-r}\big\}_{1\le r\le q}$
being known, successive approximations $\psi_{\h}^{n,m}$ of
$\psi_{\h}^{n}$ are computed using a quasi-Newton procedure
of the form  
\begin{equation} 
\label{schematemp_4}
\begin{array}{l}
\forall \tau \in \mathcal{T}_{\h},\ \forall \phi \in \mathbb{P}_p(\tau),\\
\quad \quad
\displaystyle{\frac{\alpha_0^q}{\delta t}\int_{\tau} \big( \theta(\psi_{\h}^{n,m}) + 
\partial_{\psi} [\theta(\psi_{\h}^{n,m})](\psi_{\h}^{n,m+1}-\psi_{\h}^{n,m})
\big) \phi }\\ 
\quad \displaystyle{+\ a_{\tau}(\psi_{\h}^{n,m},\psi_{\h}^{n,m+1},\phi) = b_{\tau}(\psi_{\h}^{n,m},\phi) -
\sum_{r=1}^{q}\frac{\alpha_r^q}{\delta t} \int_{\tau}
\theta(\psi_{\h}^{n-r}) \phi.}
\end{array} 
\end{equation} 
Let $\delta \psi_{\h}^{n,m} =
\psi_{\h}^{n,m+1}-\psi_{\h}^{n,m}$ and let $d_{\tau}$ be
defined as $\displaystyle{
d_{\tau}(\zeta,\psi,\phi) \eqdef \int_{\tau} \partial_{\psi}[\theta(\zeta)]\psi \phi,}$
so that equation ($\ref{schematemp_4}$) can be written as

\begin{equation} 
\label{schematemp_5}
\begin{array}{l}
\forall \tau \in \mathcal{T}_{\h},\quad \forall \phi \in \mathbb{P}_p(\tau),\\
\quad \displaystyle {\frac{\alpha_0^q}{\delta t}
d_{\tau}(\psi_{\h}^{n,m},\delta \psi_{\h}^{n,m},\phi) +
a_{\tau}(\psi_{\h}^{n,m},\delta\psi_{\h}^{n,m},\phi)}=b_{\tau}(\psi_{\h}^{n,m},\phi)\\
\displaystyle{\quad \quad -
\sum_{r=1}^{q}\frac{\alpha_r^q}{\delta t} \int_{\tau}
\theta(\psi_{\h}^{n-r}) \phi- \frac{\alpha_0^q}{\delta t}
\int_{\tau} \theta(\psi_{\h}^{n,m}) \phi
-a_{\tau}(\psi_{\h}^{n,m},\psi_{\h}^{n,m},\phi).}
\end{array}
\end{equation} 

\begin{algorithm}%[H]
 \caption{\hspace{5pt} Iterative algorithm at each time step for solving Richards' equation}
 \vspace{3pt}
\label{alg1}
 \begin{algorithmic}%[1]
   \REQUIRE $\psi_{\h}^{n-1},\psi_{\h}^{n-2},\cdots,\psi_{\h}^{n-q},\psi_{\h}^{n,0},\epsilon_{alg1}$ \\[3pt]
   \STATE set $m=0$ \\[2pt]
   \REPEAT
   \STATE \vspace{3pt} Find $\delta\psi_{\h}^{n,m}\in V_{\h}$ solving ($\ref{schematemp_5}$)\\[2pt]
   \STATE set $\psi_{\h}^{n,m+1} = \psi_{\h}^{n,m} + \delta\psi_{\h}^{n,m}$  \\[2pt]
   \STATE $m \gets m+1$ \\[2pt]
   \UNTIL{$E \le \epsilon_{alg1}$}  \\[4pt]
%   \STATE set $\psi_{\h}^{n} = \psi_{\h}^{n,m}$  \\[2pt]
   \ENSURE $\psi_{\h}^{n}=\psi_{\h}^{n,m}$
\end{algorithmic}
\end{algorithm}

\vspace{3mm}
\noindent The simplest initialization of Algorithm~\ref{alg1} consisting of choosing the approximation of the
solution at the previous time step ($\psi_{\h}^{n,0} =
\psi_{\h}^{n-1}$), but a higher order initialization can
also be used (see~\S\ref{sec:results}). The error measure
$E$ is the relative Euclidean norm of the component vector
associated with $\delta \psi_{\h}^{n,m}$, and
$\epsilon_{alg1}$ is a user-defined convergence criterion.  

\subsection{Discretization of the kinematic wave equation}
\label{sec:kinematicwave-discretization} 
\begin{figure}[htb]
\centering
\psset{unit=0.8}
\begin{pspicture}(-2,-1)(12,2.3)
%--- Les indic des points -------------
\rput(2,-0.7){$t$}
\rput(2,1.7){$t$}
\rput(8.1,-0.7){$t+\delta t$}
\rput(8.2,1.7){$t+n'\delta t'$}
\rput(5,0.4){$\delta t$}
\rput(2.5,1.45){$\delta t'$}
\rput(12,0.5){time}
\rput(-1.5,0.3){subsurface flow}
\rput(-1.5,1.3){overland flow}
\rput(3,1){$|$}
\rput(4,1){$|$}
\rput(5,1){$|$}
\rput(6,1){$|$}
\rput(7,1){$|$}
%--- Boite ---------------------------
\psline[linestyle=dashed](-2,0)(1,0)
\psline[linestyle=dashed](-2,1)(1,1)
\psline[linestyle=dashed]{->}(9,0)(12,0)
\psline[linestyle=dashed]{->}(9,1)(12,1)
\psline(1,0)(9,0)
\psline(1,1)(9,1)
\psline[linewidth=0.1mm](2,-0.4)(2,1.4)
\psline[linewidth=0.1mm](8,-0.4)(8,1.4)
\pscurve[linestyle=dotted]{->}(2,1)(2.5,1.2)(3,1)
\pscurve[linestyle=dotted]{->}(3,1)(3.5,1.2)(4,1)
\pscurve[linestyle=dotted]{->}(4,1)(4.5,1.2)(5,1)
\pscurve[linestyle=dotted]{->}(5,1)(5.5,1.2)(6,1)
\pscurve[linestyle=dotted]{->}(6,1)(6.5,1.2)(7,1)
\pscurve[linestyle=dotted]{->}(7,1)(7.5,1.2)(8,1)
\pscurve[linestyle=dotted]{->}(2,0)(5,0.2)(8,0)
\end{pspicture}
\caption{Multiple time-stepping for subsurface and overland flows.\label{Fig:sV2}}
\end{figure}
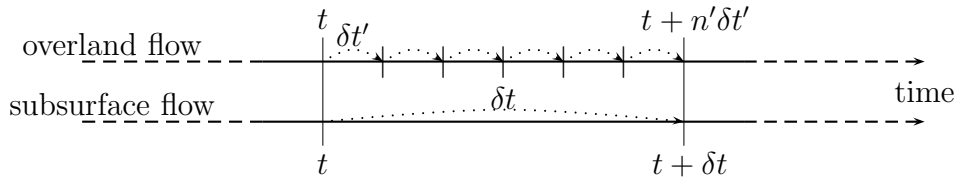

\begin{figure}[htb]
\centering
\psset{unit=0.8}
\begin{pspicture}(-1,-0.5)(14,3.5)
%--- Boite ---------------------------
\psline(3.5,0.9)(12.5,0)
\psline[linestyle=dashed](-1,1.35)(1,1.15)
\psline[linestyle=dotted](1,1.15)(3.5,0.9)
\psline[linestyle=dotted](12.5,0)(14,-0.15)
\psline[linewidth=0.1mm]{<->}(6.5,0.2)(9.7,-0.12)
\psline[linewidth=0.1mm](-1,2.7)(1,2.7)
\psline[linewidth=0.1mm](3.5,2)(6.5,2)
\psline[linewidth=0.1mm](6.5,1.7)(9.5,1.7)
\psline[linewidth=0.1mm](9.5,1)(12.5,1)
%--- Les indic des points -------------
\rput(1,1.15){$\bullet$}
\rput(3.5,0.9){$\bullet$}
\rput(6.5,0.6){$\bullet$}
\rput(9.7,0.28){$\bullet$}
\rput(12.5,0){$\bullet$}
\rput(8.1,-0.3){$l_i$}
\rput(1,1.47){$x_{-\ud}$}
\rput(1,0.83){$\rm{A}$}
\rput(3.7,1.22){$x_{i-\frac{3}{2}}$}
\rput(6.5,0.92){$x_{i-\ud}$}
\rput(9.7,0.6){$x_{i+\ud}$}
\rput(12.5,0.32){$x_{i+\frac{3}{2}}$}
\rput(0,3){$h_{-1}$}
\rput(5,2.3){$h_{i-1}$}
\rput(8.1,2){$h_{i}$}
\rput(11,1.3){$h_{i+1}$}
\rput(8.1,0.44){$\times$}
\rput(8.1,0.76){$x_i$}
\end{pspicture}
\caption{Space discretization at the ground surface.\label{Fig:sV1}}
\end{figure}
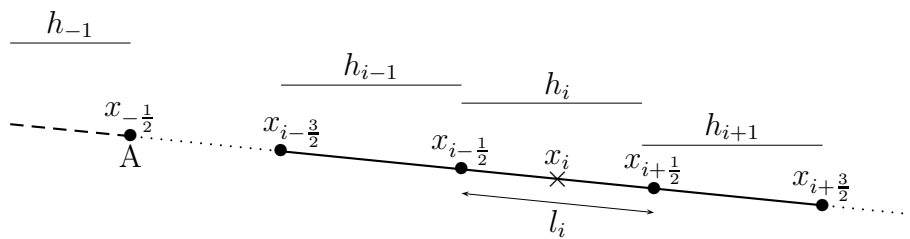

The kinematic wave equation is discretized on a surface mesh
on $\I$ which is simply the trace of the mesh $\Th$ on
$\I$. Let $N_{\I}$ be the number of mesh faces covering
$\I$. A finite volume scheme with Godunov flux and time step
$\delta t'$ is used. The time step is taken less than or
equal to the time step for Richards' equation, that is,
$\delta t' = \delta t / n'$ with $ n'\ge1$
(see Figure~\ref{Fig:sV2}). This choice is made because the
explicit FV scheme is, as usual, restricted by a CFL
condition to ensure its stability. This is not the case for
the discrete Richards' equation where, owing to the use of a
BDF, a larger time step can be employed. This leads to the
following notation: $h_{\h}^{n,k}$ for $n\in\{1\cdots N_T\}$
and $k\in\{0\cdots n'\}$ denotes the discrete approximation
of $h$ at time $n\delta t + k\delta t'$ and for brevity we
write $h_{\h}^{n}\eqdef h_{\h}^{n,0}=h_{\h}^{n-1,n'}$.
Let $x_i$, $l_i$, $x_{i-\ud}$ and $x_{i+\ud}$ be defined on
a generic mesh face $e_i$ on $\I$ respectively as the center, the length,
and the left and right vertices of $e_i$
(see Figure~\ref{Fig:sV1}). $S_i$ denotes the slope of the face 
$e_i$. Since the flux function $\varphi$ is convex and the
surface water depth is nonnegative, the Godunov flux
coincides with the upwind flux, yielding 
\begin{align}
\label{Eq:SV6}
&\forall k\in\{1\cdots n'\},\ \forall i\in\{1\cdots N_{\I}\},\nonumber \\
& h_{i}^{n-1,k} =  h_{i}^{n-1,k-1} + \frac{\delta t'}{l_i}\Big(\varphi(h_{i-1}^{n-1,k-1},S_{i-1}) - \varphi(h_{i}^{n-1,k-1},S_{i}) \\
& \hphantom{h_{i}^{n-1,k} =  h_{i}^{n-1,k-1} + \frac{\delta t'}{l_i}\Big(} -l_i \vr^{n-1,k-1}\cdot n_{\Omega} + \int_{e_i} v^{\star,n}_{\h}\Big), \nonumber 
\end{align}
where for all $i\in\{1\cdots N_{\I}\},\
h_i^{n,k}\eqdef h_{\h}^{n,k}|_{e_i}$ and $v^{\star,n}_{\h}$ is a
discrete interface flux yet to be defined (see $\S
3.3$). Observe that a fixed interface flux is used for the
multiple time steps comprised in a single time step of
Richards' equation.  
Equation (\ref{Eq:SV6}) requires the knowledge of the
surface water depth at $t=0$ (initial condition) and to the
left of the first face on a fictitious cell (boundary
condition) at all discrete times, $\forall i \in \{1\cdots
N_{\I}\}$, $h_i^0 = h^0(x_i)$ and $\forall n\in\{1\cdots
N_T\},\forall k\in\{0\cdots n'-1\}$, $h_{-1}^{n,k} = h_{\rm{A}}^{n,k}$.

The CFL condition for the explicit scheme (\ref{Eq:SV6}) can
be expressed as 
\begin{equation}
\label{Eq:SV7}
\delta t'\le \min_{1\le i\le N_{\I}}\big(\frac{l_{i}}{\partial_h\varphi(h_{\max},S_i)}\big),
\end{equation}
where $h_{\max}$ is an \textit{a priori} bound for the
surface surface water depth $h$ on $\I\times[0,T]$. By
definition of the flux function $\varphi$, this yields 
$$\delta t'\le
\frac{3}{5\mathcal{K}h_{\max}^{2/3}}\cdot\min_{1\le i\le
N_{\I}}\Big(l_{i}S_i^{-\ud}\Big). $$
In the absence of rainfall and coupling terms, the
satisfaction of the CFL condition implies a discrete maximum
principle and a decrease in the total variation for the
discrete surface water depth.

\subsection{Single-step coupling algorithm}
\label{single-step-algorithm}
We consider in this section the case where Richards'
equation is discretized in time using a first-order BDF
(that is, the Euler implicit scheme). Together with the
finite volume scheme described in
\S\ref{sec:kinematicwave-discretization} for the kinematic
wave equation, this yields a scheme to approximate the
coupled system (\ref{Eq:Couplage3}) provided we specify the
time evolution of the coupling variables
$\{\Ihd{n},\Ihw{n},\omega_{v}^n,\omega_{\psi}^n\}$ for
$n\in\{1\cdots N_T\}$ (see~\S\ref{sec:admissibleset}). Here,
as before, the superscript $n$ stands for the value at
$n\delta t$, so that $\Ihd{n}=\Ihd{n\delta t}$ and so
on. This time evolution is designed with the twofold
objective to ensure that a suitable approximation of
($\psi,h$) lies in the admissible set $\A$ at all discrete
times and to ensure overall mass conservation for the whole
system (subsurface and overland flows). The resulting
algorithm is outlined in Algorithm~\ref{alg2}. It is termed
single-step coupling algorithm in reference to the use of
the first-order BDF which spans a single time step
interval. For simplicity in the presentation of
Algorithm~\ref{alg2}, we define
\begin{itemize}
\item $\psi_{\h}^n\gets$
\texttt{Richards$\_$BDF1}$(\Ihd{n},\Ihw{n},\omega_{v}^n,\omega_{\psi}^n,\psi_{\h}^{n-1})$
as the resolution by Algorithm~\ref{alg1} of Richards'
equation on a time step by the SIPG method, the first-order
BDF and boundary data on $\I$ determined from
$\{\Ihd{n},\Ihw{n},\omega_{v}^n,\omega_{\psi}^n\}$,\\[-8pt]
\item $h_{\h}^n\gets$
\texttt{Kinematic$\_$wave}$(h_{\h}^{n-1},n',v_{\rm{r}},v_{\h}^{\star,n})$
as the resolution of the kinematic wave equation by using
($\ref{Eq:SV6}$) $n'$ times,\\[-8pt] 
\item $v^{\star,n}_{\h} \gets$
\texttt{Normal$\_$Velocity}$(\Ihd{n},\Ihw{n},\omega_{v}^n,\omega_{\psi}^n,\psi_{\h}^n)$
as the evaluation of the interface normal velocity
$v^{\star,n}_{\h}$ on $\I$ defined as 
\begin{equation}
v^{\star,n}_{\h}|_{F} \eqdef
\left\{\begin{array}{ll}
\omega_v^n|_{F} &\mbox{if}\ F \in \Ihd{n},\\[3pt]
v(\psi_{\h}^n|_{F})\cdot n_{\Omega} + \eta K_s d_F^{-1}(\psi_{\h}^n|_{F} - \omega_{\psi}^n|_{F}) & \mbox{if}\ F \in \Ihw{n}.
\end{array}\right.
\end{equation}
Note that the expression for $v^{\star,n}_{\h}$ on
$\Iw{n}_{\h}$ corresponds to the normal component of the
$H(\rm{div},\Omega)$-conforming velocity reconstruction
derived in \cite{Ern07} for DG methods. 
\end{itemize}
\begin{algorithm}%[H]
\caption{\hspace{5pt} Single-step coupling algorithm}
\vspace{3pt}
\label{alg2}
 \begin{algorithmic}%[1]
   \REQUIRE $\psi_{\h}^{n-1}$ and $h_{\h}^{n-1}$\\[4pt]
   \STATE $\tilde{h}_{\h}^{n} \gets$\texttt{Kinematic$\_$wave}$(h_{\h}^{n-1},n',v_{\rm{r}}^{n-1},0)$\\[4pt]
   \STATE Set $p=0$ and $h_{\h}^{n,1}=\tilde{h}_{\h}^{n}$, \\[4pt]
   \REPEAT
   \STATE $p \gets p+1$
   \STATE \vspace{5pt} Partition of $\I$: $\Ihd{n,p} = \{ e_i \in \mathcal{F}^{\I}_{\mathfrak{h}}, \forall k\in\{1\cdots p\}, h^{n,k}_i < 0 \}$ and $\Ihw{n,p} =    \I\backslash \Ihd{n,p}$\\[4pt]
   \STATE $\omega_{v}^{n,p}  \gets -\tilde{h}_{\h}^{n}/\delta t\ \mbox{on}\  \Ihd{n,p}$\\[4pt]
   \STATE $\omega_{\psi}^{n,p}  \gets \tilde{h}_{\h}^{n}\ \mbox{on}\  \Ihw{n,p}$\\[4pt]
   \STATE $\psi_{\h}^{n,p} \gets$ \texttt{Richards$\_$BDF1}$(\Ihd{n,p},\Ihw{n,p},\omega_{v}^{n,p},\omega_{\psi}^{n,p},\psi_{\h}^{n-1})$\\[4pt]
   \STATE $v^{\star,n,p}_{\h} \gets$ \texttt{Normal$\_$Velocity}$(\Ihd{n,p},\Ihw{n,p},\omega_{v}^{n,p},\omega_{\psi}^{n,p},\psi_{\h}^{n,p})$\\[4pt] 
   \STATE $\forall i \in \{1\cdots N_{\I}\}, h_{i}^{n,p} = \tilde{h}_{i}^{n}+\delta t/l_i \int_{e_i} v^{\star,n,p}_{\h}$ \\[5pt]
   \UNTIL{$\forall\ i \in \{1\cdots N_{\I}\}, h^{n,p}_i \ge 0 $}  \\[5pt]
   \ENSURE $\psi_{\h}^{n} = \psi_{\h}^{n,p}$ and $h_{\h}^{n} = h_{\h}^{n,p}$
\end{algorithmic}
\end{algorithm}

The principle of Algorithm~\ref{alg2} is the
following. Firstly, the surface water depth is predicted
without subsurface coupling term
($v_{\h}^{\star,n}=0$). This predicted surface water depth 
$\tilde{h}_{\h}^n$ then serves as a Dirichlet boundary
condition for Richards' equation. Because the Godunov scheme
satisfies a discrete maximum principle, $\tilde{h}_{i}^n\ge
0$ for all $i\in\{1\cdots N_{\I}\}$, so that $\Ihd{n,1} =
\emptyset$ and $\Ihw{n,1} = \I$. That is, we begin the
iterations by assuming that $\I$ is totally wet. Thus, for
$p=1$, the determination of $\omega_v^{n,p}$ is
irrelevant. Then, Richards' equation is solved and a first
estimate of the normal velocity $v_{\h}^{\star,n,p}$ is used to evaluate
the surface water depth $h_{\h}^{n,p}$. The sign of
$h_{\h}^{n,p}$ is subsequently checked on the faces of
$\I$. If $h_{i}^{n,p}$ is nonnegative on all faces, the
evaluation of the hydraulic head and of the surface water
depth can be accepted as the solution to the coupled system
at time $n\delta t$. Otherwise, a new partition of $\I$ is
determined and a Neumann condition is enforced on those
faces where the surface water depth is negative. This
Neumann condition is evaluated in such a way that at the
corresponding interface cells, the surface water is
completely infiltrated into the soil since $\omega_{v}^{n,p}= -\tilde{h}_{\h}^{n}/\delta t$. A new hydraulic head
and a new surface water depth are then calculated and the loop is
repeated until convergence. Note that convergence occurs
since the set $\Ihd{n,p}$ increases with $p$ while the set
$\Ihw{n,p}$ decreases.

\subsubsection*{Admissibility of $(\psi,h)$}
An important point is that Algorithm~\ref{alg2} delivers
nonnegative surface water depths. Moreover, on the wet part
of the interface, there holds 
$$\forall n\in\{1\cdots N_T\},\quad \forall F\in
\Ihw{n},\quad \psi_{\h}^n|_F=\tilde{h}_{\h}^{n}|_F,$$ 
since the value of the Dirichlet data $\omega_{\psi}^{n,p}$
on $\Ihw{n,p}$ is fixed during the loop. This is not the
condition $\psi=h$ enforced by the admissible set $\A$ but
an $\mathcal{O}(\delta t)$ approximation of it. Furthermore,
on the dry part of the interface, the surface water depth is
equal to zero and there holds 
$$\forall n\in\{1\cdots N_T\},\quad \forall F\in
\Ihd{n},\quad \psi_{\h}^n|_F\le\tilde{h}_{\h}^{n}|_F.$$
Again, this is an $\mathcal{O}(\delta t)$ approximation of
the condition $\psi\le 0$ enforced by the admissible
set. Furthermore, we observe that, if on a given face $e_i$,
the surface water depth $h_{i}^{n-1}$ is zero as well as the
upwind fluxes over the time step $[(n-1)\delta t,n\delta t]$,
the Neumann condition on Richards' equation is equal to the
rainfall intensity. Moreover, we observe that in contrast to
front tracking schemes, Algorithm~\ref{alg2} does not use
any information from the previous time step to determine the
wet portion of the interface. This offers the advantage of
robustness and ease of extension to 3D/2D settings, but can
entail higher computational costs than those incurred by
front tracking schemes in the absence of exfiltration (see
for instance \cite{Beaugendre06}).  

\subsubsection*{Overall mass conservation}
The total volume of water in the domain $\Omega$ at time
$n\delta t$ is obtained by integrating the volumetric water
content in $\Omega$
$$V_{\mathrm{grnd}}^n \eqdef \int_{\Omega}\theta(\psi_{\h}^n).$$
Taking the test function $\phi$ equal to $1$ in the SIPG
scheme (\ref{schematemp_5}), summing over the mesh elements
and using the first-order BDF to approximate the non
stationary term yields 
\begin{equation}
\label{Eq:SC1}
V_{\mathrm{grnd}}^{n}-V_{\mathrm{grnd}}^{n-1}=
\big(F_{\I}^{n} + F_{\mathcal{WB}}^{n}\big) \delta
t+\epsilon^{n}, 
\end{equation}
where $F_{\I}^{n}$ (resp. $F_{\mathcal{WB}}^{n}$) is the
flux over the time step $[(n-1)\delta t,n\delta t]$ across
the interface $\I$ (resp. across the bottom and lateral
walls), 
\begin{equation}
\label{Eq:SC1bis}
F_{\I}^{n} \eqdef -\int_{\I}
v^{\star,n}_{\h}\quad\text{and}\quad F_{\mathcal{WB}}^{n}
\eqdef -\int_{\WB} v_N^n, 
\end{equation}
and $\epsilon^{n}$ represents the numerical error in the
resolution of the nonlinear system. Recall that
$|\epsilon^n|\le C\epsilon_{alg1},$ where $\epsilon_{alg1}$
is the convergence tolerance of Algorithm~\ref{alg1} and $C$ a
constant due to the fact that the convergence criterion in
Algorithm~\ref{alg1} limits the norm of the variation of the
hydraulic head $\psi$ rather than the one of the volumetric
water content $\theta(\psi)$. The total volume of water in the overland 
flow at time $n\delta t$ is obtained by integrating the
surface water depth over $\I$ 
$$V_{\mathrm{over}}^n \eqdef \int_{\I} h_{\h}^{n}.$$
The total variation of water volume in the kinematic wave
equation over the time step $[(n-1)\delta t,n\delta t]$ is
obtained by summing the elementary contributions in equation
(\ref{Eq:SV6}), yielding 
\begin{equation}
\label{Eq:SC2}
V_{\mathrm{over}}^{n} - V_{\mathrm{over}}^{n-1} =
\big(-F_{\I}^n + F_{\rm{ABr}}^n \big)\delta t, 
\end{equation}
where $F_{\I}^{n}$ is already defined above and where
$F_{\rm{ABr}}^{n}$ represents the water flux over the
time step $[(n-1)\delta t,n\delta t]$ due to the rain and the
discharge at points $\rm{A}$ and $\rm{B}$, with
$F_{\rm{ABr}}^{n} \eqdef F_{\rm{A}}^n + F_{\rm{B}}^n +
F_{\rm{r}}^n$, and
$$F_{\rm{A}}^n \eqdef \frac{\delta t'}{\delta
t}\sum_{k=1}^{n'}\varphi(h_{\rm{A}}^{n-1,k}),\ 
F_{\rm{B}}^n \eqdef -\frac{\delta t'}{\delta
t}\sum_{k=1}^{n'}\varphi(h_{N_{\I}}^{n-1,k}),\ 
F_{\rm{r}}^n \eqdef -\frac{\delta t'}{\delta t}\sum_{k=1}^{n'}\int_{\I} \vr^{n-1,k} \cdot n_{\Omega}.$$
The total volume of water contained in the coupled system is
the sum of the volume of each system, $V^{n} \eqdef
V_{\mathrm{grnd}}^{n} + V_{\mathrm{over}}^{n}$. When
($\ref{Eq:SC1}$) and ($\ref{Eq:SC2}$) are summed, the
interface flux cancels, yielding 
\begin{equation}
\label{Eq:SC3}
V^{n} - V^{n-1} = \big(F_{\mathcal{WB}}^n + F_{\rm{ABr}}^n \big)\delta t + \epsilon^n.
\end{equation}
This relation readily implies the following overall water
volume conservation result for the single-step algorithm.
\begin{theo}
\label{Theo:MC1}
Let $\delta V^n$ be the overall water volume defect over the
time step $[(n-1)\delta t,n\delta t]$ defined as 
$\delta V^n \eqdef V^n - V^{n-1} -
(F_{\mathcal{WB}}^n+F_{\rm{ABr}}^{n})\delta t$. 
Let $\triangle V^n$ be the overall water volume defect over
the time interval $[0,n\delta t]$ defined as $\triangle V^n
\eqdef \sum_{i=1}^{n} \delta V^i$. Then, 
\begin{equation}
|\triangle V^n| \le nC\epsilon_{alg1}, 
\end{equation}
where $C$ is a constant and $\epsilon_{alg1}$ is the
tolerance in Algorithm~\ref{alg1}.
\end{theo}

\subsection{Two-step coupling algorithm}
\label{two-step-algorithm}
We consider in this section the case where Richards'
equation is discretized in time using a second-order BDF for
which 
\begin{equation}
\label{BDF2}
\Big(\frac{\partial \chi}{\partial t}\Big)^n\simeq \frac{1}{\delta t}
\Big(\frac{3}{2}\chi^{n} -2 \chi^{n-1} +\frac{1}{2} \chi^{n-2} \Big).
\end{equation}
The single-step coupling algorithm is not conservative when
the non-stationary term of Richards' equation is
approximated by the second-order BDF owing to the fact that
the Euler explicit scheme used to solve the kinematic wave
equation spans only a single time step. Consequently, to
obtain a mass conservative scheme, the interface flux
$F_{\I}^n$ used in the kinematic wave equation needs to be
transformed into a new interface flux $\Phi_{\I}^n$, so that
(\ref{Eq:SC2}) becomes 
\begin{equation}
\label{Eq:SC5}
V_{\mathrm{over}}^{n} - V_{\mathrm{over}}^{n-1} =
\big(-\Phi_{\I}^n + F_{\rm{ABr}}^{n}\big)\delta t. 
\end{equation}
To identify the expression for $\Phi_{\I}^n$, observe that
using a second-order BDF modifies ($\ref{Eq:SC1}$) into 
$$\frac{3}{2}V_{\mathrm{grnd}}^{n}-2
V_{\mathrm{grnd}}^{n-1}+\frac{1}{2}V_{\mathrm{grnd}}^{n-2} =
\big(F_{\I}^{n} + F_{\mathcal{WB}}^{n}\big) \delta
t+\epsilon^{n},$$ 
which can be rewritten as 
\begin{equation}
\label{Eq:SC6}
\frac{3}{2}\big( V_{\mathrm{grnd}}^{n}-V_{\mathrm{grnd}}^{n-1}\big)-\frac{1}{2}\big(
V_{\mathrm{grnd}}^{n-1} - V_{\mathrm{grnd}}^{n-2} \big) -
F_{\mathcal{WB}}^{n}\delta t = F_{\I}^{n} \delta t +
\epsilon^{n}, 
\end{equation}
where the fluxes $F_{\I}^{n}$ and $F_{\mathcal{WB}}^{n}$ are
still defined by (\ref{Eq:SC1bis}). Moreover it results from
(\ref{Eq:SC5}) that 
\begin{align}
&\frac{3}{2}\big(V_{\mathrm{over}}^{n} -
V_{\mathrm{over}}^{n-1}\big)-\frac{1}{2}\big(V_{\mathrm{over}}^{n-1}
- V_{\mathrm{over}}^{n-2}\big)
+\Big(-\frac{3}{2}F_{\rm{ABr}}^{n}+
\frac{1}{2}F_{\rm{ABr}}^{n-1}\Big) \delta t\nonumber\\
&\qquad\qquad\qquad\qquad\qquad\qquad\qquad\qquad\qquad
=\Big(-\frac{3}{2} \Phi_{\I}^{n} +\frac{1}{2}
\Phi_{\I}^{n-1}\Big)\delta t. \label{Eq:SC7}
\end{align} 
The new interface flux $\Phi_{\I}^n$ is determined so that
the mass flux $F_{\I}^n$ is exactly counter-balanced by the
interface flux in (\ref{Eq:SC7}), whence 
$$F_{\I}^n = \frac{3}{2} \Phi_{\I}^{n} -\frac{1}{2}
\Phi_{\I}^{n-1}\quad \Longrightarrow \quad \Phi_{\I}^{n} =
\frac{2}{3}F_{\I}^n + \frac{1}{3} \Phi_{\I}^{n-1}. $$
At the first time step where a one-step implicit scheme is
used, water volume conservation is directly enforced by
setting $\Phi_{\I}^{1} \eqdef F_{\I}^{1}$. 

\begin{algorithm}%[H]
 \caption{\hspace{5pt} Two-step coupling algorithm}
 \vspace{3pt}
\label{alg3}
 \begin{algorithmic}%[1]
   \REQUIRE $\psi_{\h}^{n-1}, \psi_{\h}^{n-2}$ and $h_{\h}^{n-1}$\\[4pt]
   \STATE $\vdots$\\[4pt]
   \REPEAT
   \STATE $\vdots$\\[4pt]
   \STATE $\omega_{v}^{n,p} \gets -(3\tilde{h}_{\h}^{n}/\delta t+v^{\star,n-1}_{\h})/2\ \mbox{on}\  \Ihd{n,p}$\\[4pt]
   \STATE $\psi_{\h}^{n,p} \gets$ \texttt{Richards$\_$BDF2}$(\Ihd{n,p},\Ihw{n,p},\omega_{v}^{n,p},\omega_{\psi}^{n,p},\psi_{\h}^{n-1},\psi_{\h}^{n-2})$\\[4pt] 
   \STATE $\vdots$\\[4pt]
   \STATE $\forall i \in \{1\cdots N_{\I}\},\ h_{i}^{n,p} = \tilde{h}_{i}^{n}+\delta t/l_i \int_{e_i} (2v^{\star,n,p}_{\h}+v^{\star,n-1}_{\h})/3 $ \\[5pt]
   \STATE $\vdots$\\[4pt]
   \UNTIL{$\forall\ i \in \{1\cdots N_{\I}\}, h^{n,p}_i \ge 0 $}\\[5pt]
   \STATE $\vdots$\\[4pt]
   \ENSURE $\psi_{\h}^{n}, h_{\h}^{n}$ and $v^{\star,n}_{\h} = v^{\star,n,p}_{\h}$
\end{algorithmic}
\end{algorithm}

The resulting algorithm, referred to as two-step coupling
algorithm, is outlined in Algorithm~\ref{alg3}. Only the
differences with Algorithm~\ref{alg2} are indicated. The key
modification concerns the evaluation of the interface normal
velocity in the calculation of $h_{i}^{n,p}$. The Neumann
data $\omega_{v}^{n,p}$ is also modified to ensure that the
Neumann condition indeed leads to a dry state in the
corresponding cell. Also, the discrete approximation
$\psi_{\h}^{n-2}$ at time $(n-2)\delta t$ is added to the
input and the interface normal velocity $v_{\h}^{\star,n}$
is added to the output at each time step since it is used in
the subsequent time step. 

The main result concerning the overall water volume
conservation for Algorithm~\ref{alg3} is the following. 
\begin{theo}
Let $\delta V^n$ be the overall water volume defect over the
time step $[(n-1)\delta t,n\delta t]$ defined as  
\begin{equation}
\label{Theo:MC2}
\delta V^n \eqdef V^n - V^{n-1} - (\tilde{F}_{\mathcal{WB}}^n+F_{\rm{ABr}}^{n})\delta t,
\end{equation}
where $\tilde{F}_{\mathcal{WB}}^n \eqdef
\frac{2}{3}F_{\mathcal{WB}}^n+\frac{1}{3}\tilde{F}_{\mathcal{WB}}^{n-1}$.
Let $\triangle V^n$ be the overall water volume defect over
the time interval $[0,n\delta t]$ defined as before. Then 
$$|\triangle V^n| \le \frac{1}{2}|\delta V^1|+nC\epsilon_{alg1},$$ 
where $C$ is a constant and $\epsilon_{alg1}$ is the
tolerance in Algorithm~\ref{alg1}. 
\end{theo}
\textit{Proof}: Owing to (\ref{Eq:SC7}), the coupling terms
are eliminated when (\ref{Eq:SC5}) et (\ref{Eq:SC6}) are
summed, leading to 
$$\frac{3}{2}\big(V^{n}
-V^{n-1}\big)-\frac{1}{2}\big(V^{n-1} - V^{n-2}\big) - \Big(
F_{\mathcal{WB}}^{n}\delta t
+\frac{3}{2}F_{\rm{ABr}}^{n}\delta t -\frac{1}{2}F_{\rm{ABr}}^{n-1} \Big) \delta t = \epsilon^{n}.$$
Using the definition of $\delta V^{n}$ yields the recurrence
relation $\delta V^{n} = \frac{1}{3}\delta V^{n-1} +
\frac{2}{3}\ \epsilon^{n}$, so that 
$$\delta V^{n} = \frac{1}{3^n}\delta V^{1}+\frac{2}{3}\sum_{i=1}^{n}3^{i-n}\times \epsilon^{i}.$$
Owing to the triangle inequality, it is inferred that
$$|\triangle V^{n}| \le
\frac{3}{2}\Big(\frac{1}{3}-\frac{1}{3^{n+1}}\big)|\delta
V^{1}|+\sum_{i=1}^{n}\Big(1-\frac{1}{3^{i+1}}\Big)|\epsilon^{i}|,$$ 
whence
$$|\triangle V^{n}| \le
\frac{1}{2}|\delta V^{1}|+\sum_{i=1}^{n}|\epsilon^{i}| \le
\frac{1}{2}|\delta V^{1}| + nC\epsilon_{alg1}. \qquad\qquad\qquad\Box$$
Finally, we observe that similar developments can be
considered for arbitrary order BDFs with additional
technicalities and longer recursion formulas. 

\section{Results}
\label{sec:results}
\begin{figure}[htb]
\hspace{0.2cm}
\psfrag{Water content }{\hspace{-0.9cm}\raisebox{0.2cm}{\footnotesize{Water content $\theta$}}}
\psfrag{Hydraulic head }{\hspace{-0.8cm}\raisebox{-0.2cm}{\footnotesize{Hydraulic head $\psi$}}}
\includegraphics[height=4.4cm]{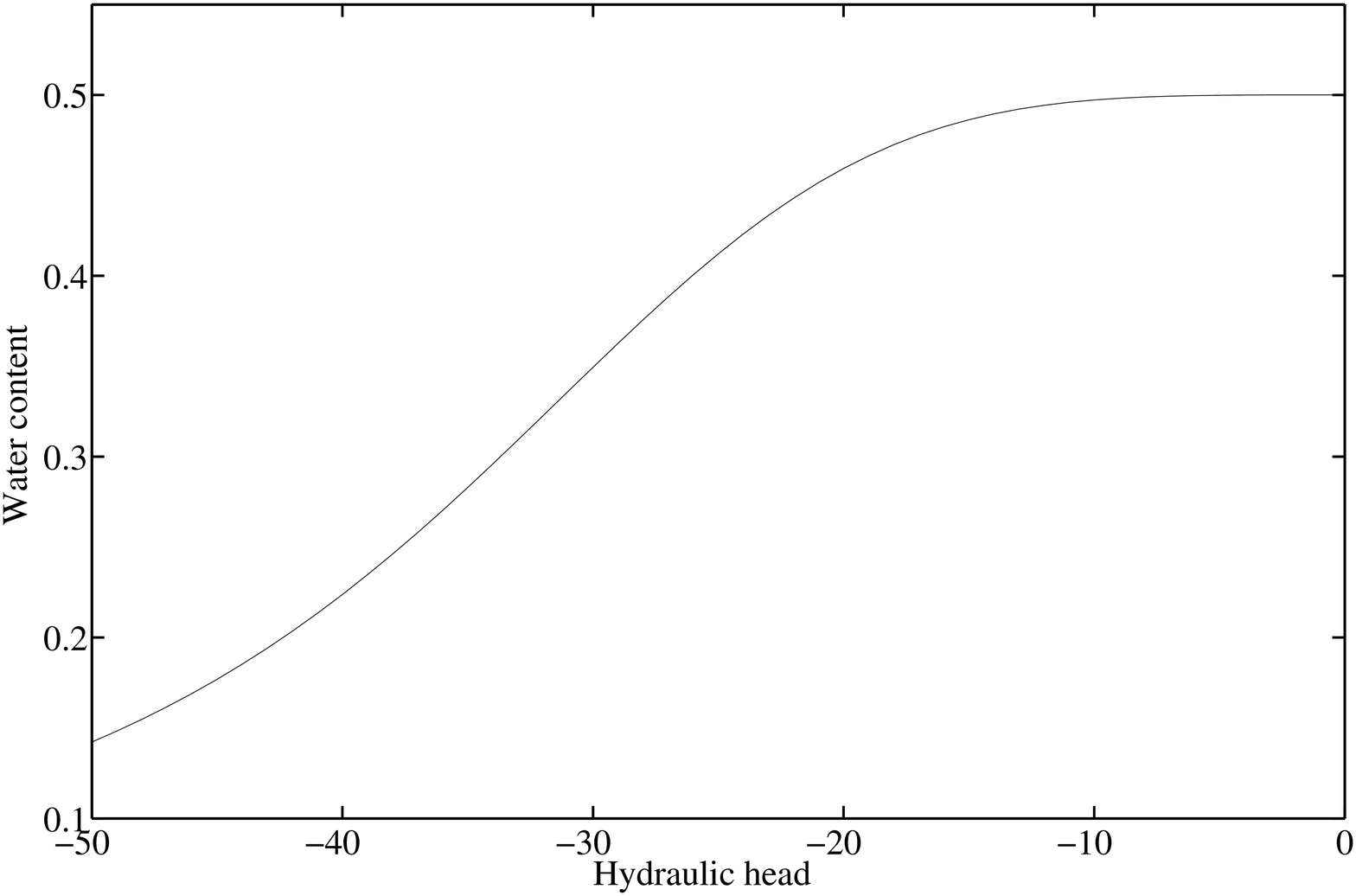}
\hspace{0.5cm}
\psfrag{Conductivite hydraulique}{\hspace{-1.4cm}\raisebox{0.2cm}{\footnotesize{Hydraulic conductivity $K$}}}
\psfrag{Charge hydraulique}{\hspace{-0.8cm}\raisebox{-0.2cm}{\footnotesize{Hydraulic head $\psi$}}}
\includegraphics[height=4.4cm]{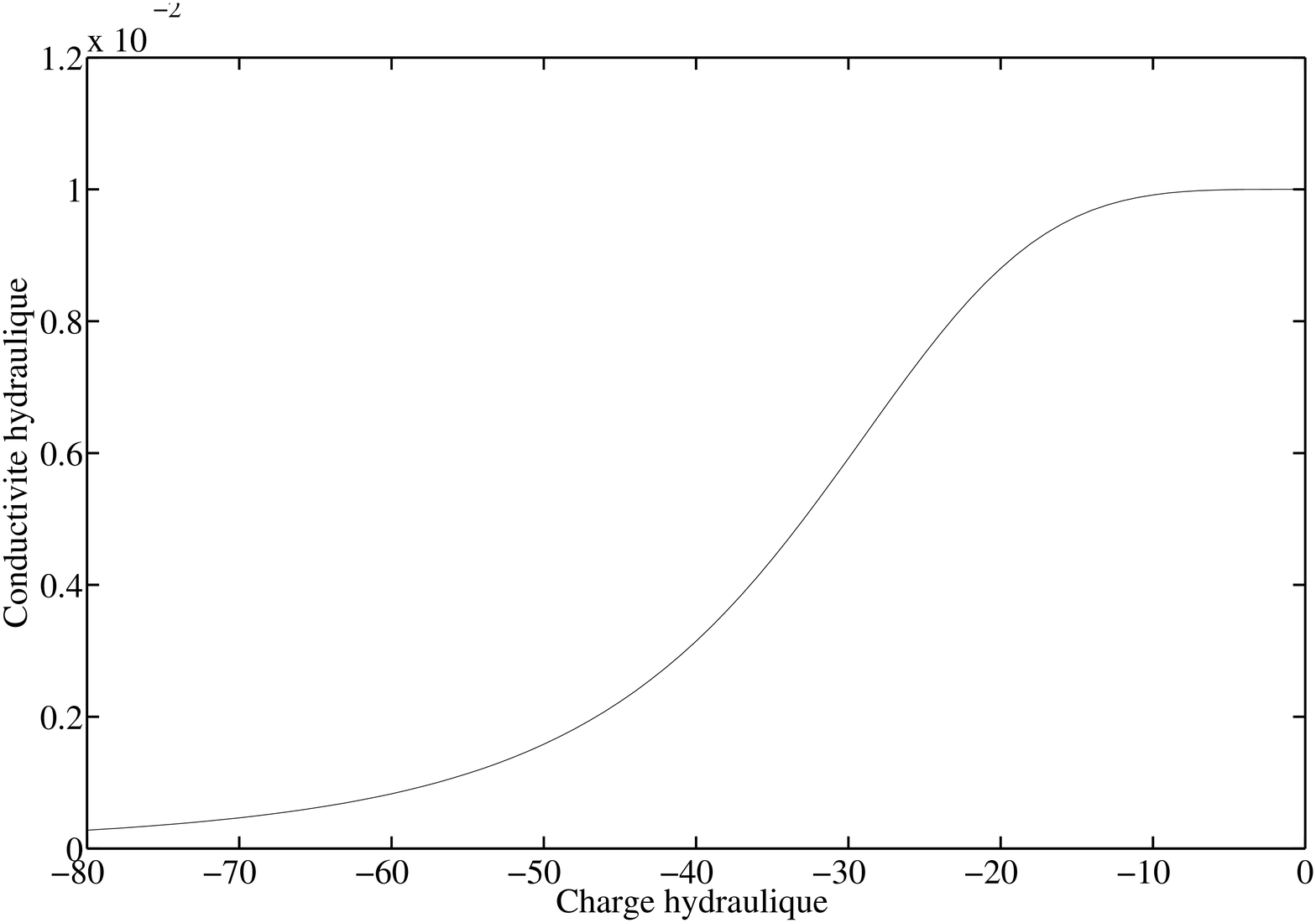}
\caption{Hydraulic properties of the soil used in the test cases.}\label{Fig:Couplage3}\vspace{0.1cm}
\end{figure}

Algorithm~\ref{alg3} is assessed on three test cases: the
first one concerns overland flow over a variable topography,
the second one infiltration due to rainfall and the third
one exfiltration resulting from injected water at the bottom
of the aquifer. The soil consists of sand and is
parameterized by the Haverkamp's constitutive relations
\cite{Haverkamp77} 
\begin{equation*}
\begin{array}{lll}
\disp{\theta(\psi)}=\disp{\frac{\theta_s-\theta_r}{1
+|\tilde{\alpha}\psi|^{\beta}}+\theta_r}\quad \text{and}&\
\disp{K(\psi)}=\disp{\frac{K_s}{1
+|\tilde{A}\psi|^{\gamma}},} 
\end{array}
\end{equation*}
with parameters
\begin{equation*}
\begin{array}{llll}
\theta_s = 0.5,   & \tilde{\alpha} = 0.028 cm^{-1}, & K_s = 10^{-2}cm.s^{-1}, &\gamma = 4,\\
\theta_r = 0.05,  & \beta = 4,  & \tilde{A} = 0.030 cm^{-1}. &\\
\end{array}
\end{equation*}
Figure~\ref{Fig:Couplage3} presents the volumetric water
content and the hydraulic conductivity as a function of the
hydraulic head. Furthermore, the Strickler coefficient
$\mathcal{K}$ is set to $60m^{1/3}s^{-1}$. For all test
cases, the bottom boundary $\B$ is located at $z=0$.
 
Piecewise affine finite elements are used ($p=1$ in
(\ref{FEspace})) along with the usual local Lagrangian basis
functions. For the first time step, the Crank--Nicolson
scheme is used. A direct solver based on the LU
decomposition is employed to solve the linear systems. The
convergence tolerance $\epsilon_{alg1}$ in Algorithm~\ref{alg1} is
set to $10^{-6}$ and the parameter $\eta$ is set to
$10$. Moreover, we focus on the use of the second-order
coupling algorithm. This choice is motivated by the fact
that it yields second-order discretization errors in time
along with second-order discretization errors in space in
the $L^2$-norm since $p=1$. In addition, a second-order
initialization of Algorithm~\ref{alg1} is chosen in the form 
\begin{equation}
\label{schematemp_6}
\forall n\ge3,\quad \psi_{\h}^{n,0} = 3\psi_{\h}^{n-1}
-3\psi_{\h}^{n-2} +\psi_{\h}^{n-3}, 
\end{equation}
except for the second time step where the first-order
initialization $\psi_{\h}^{2,0} =
2\psi_{\h}^{1}-\psi_{\h}^{0}$ is used. The second-order
initialization (\ref{schematemp_6}) can decrease
significantly the CPU time in comparison with the
initialization $\psi_{\h}^{n,0} = \psi_{\h}^{n-1}$.  

\subsection{Test case 1 (TC1)}
\label{sec:TC1}
\begin{figure}[htb]
\centering
\psset{unit=1}
\begin{pspicture}(-2,-0.5)(10,3)
%--- Boite ---------------------------
\psline[linestyle=dashed,linewidth=0.1mm](-0.6,2.22)(8,2.22)
\psline[doubleline=true](0,2.3)(0,0)(8,0)(8,2)
\psline[linestyle=dashed](8,0)(10,0)
\psline(8,2)(4.3,2.05)(3.7,2.15)(0,2.3)
\psline[linewidth=0.1mm]{<-}(0,-0.5)(3.5,-0.5)
\psline[linewidth=0.1mm]{->}(4.5,-0.5)(8,-0.5)
\psline[linewidth=0.1mm]{<-}(8.5,0.0)(8.5,0.7)
\psline[linewidth=0.1mm]{->}(8.5,1.4)(8.5,2)
\psline[linewidth=0.1mm]{<-}(-0.6,0.0)(-0.6,0.7)
\psline[linewidth=0.1mm]{->}(-0.6,1.4)(-0.6,2.22)
%--- Les indic des points -------------
\rput(4,-0.5){$3m$}
\rput(8.5,1){$0.3m$}
\rput(-0.85,1.05){$0.3025m$}
\rput(4,0.3){$v_N=0$}
\rput(7,2.42){$\bigtriangledown$}
\rput(1.5,2.8){$J_1=0.1\%$}
\rput(4,2.8){$J_2=0.3\%$}
\rput(6.5,2.8){$J_3=0.1\%$}
\rput(9.5,0.3){$z=0$}
\end{pspicture}
\vspace{0.3cm}
\caption{TC1 - Geometry, initial water table position and $h_{\rm{A}}$.}\label{Fig:Couplage4}
\end{figure}
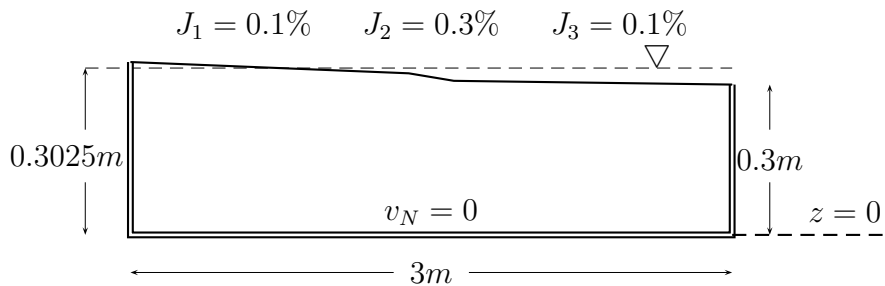

In this first test case, the runoff flow and the drainage of
the subsurface domain is induced by the presence of the
outlet, located below the initial height of the water
table. The geometry is presented in Figure~\ref{Fig:Couplage4}. 
The interface $\I$ is divided into three parts, 
$\I_1=\{ (x,z)\in\I,x\in[0,1.4]\}$ (slope $J_1 = 0.1\%$), 
$\I_2=\{ (x,z)\in\I,x\in[1.4,1.6]\}$ (slope $J_2 = 0.3\%$), and
$\I_3=\{ (x,z)\in\I,x\in[1.6,3]\}$ (slope $J_3 = 0.1\%$). 
The final simulation time is $T = 300s$. The initial
condition is an horizontal water table located at $0.3025m$
with an hydrostatic pressure profile and the boundary
condition on walls and bottom is a zero flux, 
\begin{align*}
\psi^0 = - z + 0.3025m &\quad \text{in}\ \Omega, \\
v_N = 0&\quad \text{on}\ (\WB) \times[0,T].
\end{align*}
For the overland flow, the initial condition is a horizontal
free surface and the boundary condition is a zero water depth 
\begin{align*}
h^0\ = - z + 0.3025m &\quad \text{on}\ I, \\
h_{\rm{A}}\ = 0 &\quad \text{at}\ \rm{A}\times[0,T].
\end{align*}

A mesh with 2063 triangles (corresponding to a typical
mesh-size of $3.5cm$) along with time steps $\delta t = 2.5s$ and
$\delta t' = 0.25s$ have been used. We have verified that
the interface normal velocity obtained with $\delta t =
\delta t' = 0.25s$ can be superimposed to that reported
below. In this case, the use of $\delta t = 2.5s$ instead of 
$\delta t = \delta t' = 0.25s$ leads to a gain of 89\% in
the computational time (i.e. the computation effort required
for performing one time-step in the solution of Richards'
equation is eighty times more expensive than the one of the
kinematic wave).

\begin{figure}[htb]
\begin{minipage}[c]{6.5cm}
\begin{center}
\includegraphics[width=6.5cm]{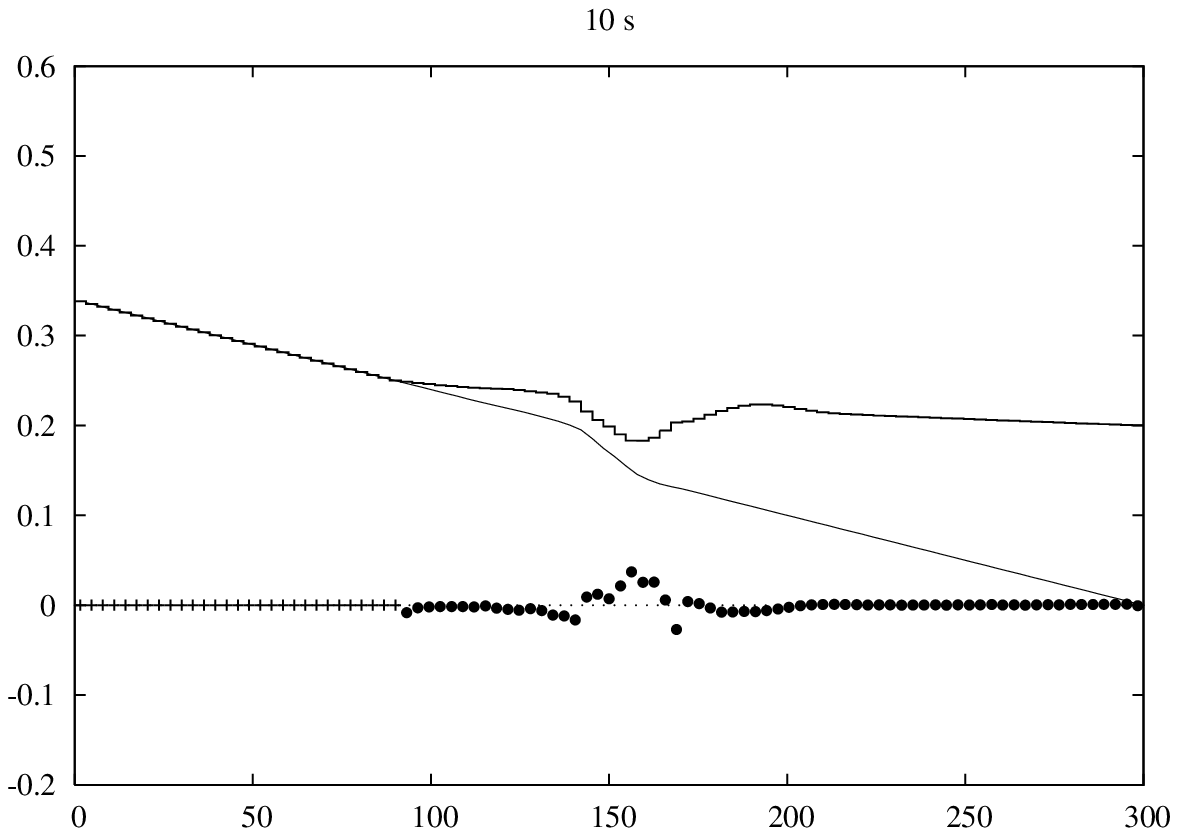}\\[6pt]
\includegraphics[width=6.5cm]{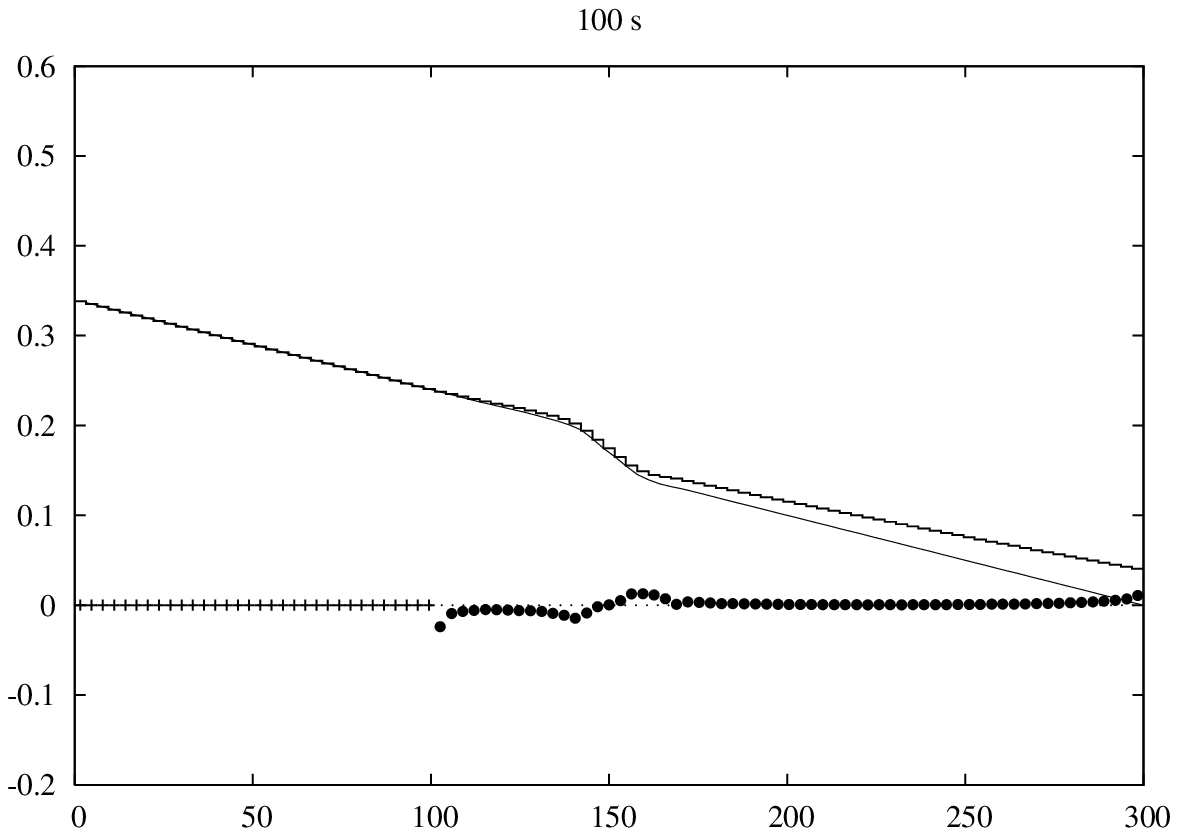}\\[6pt]
\includegraphics[width=6.5cm]{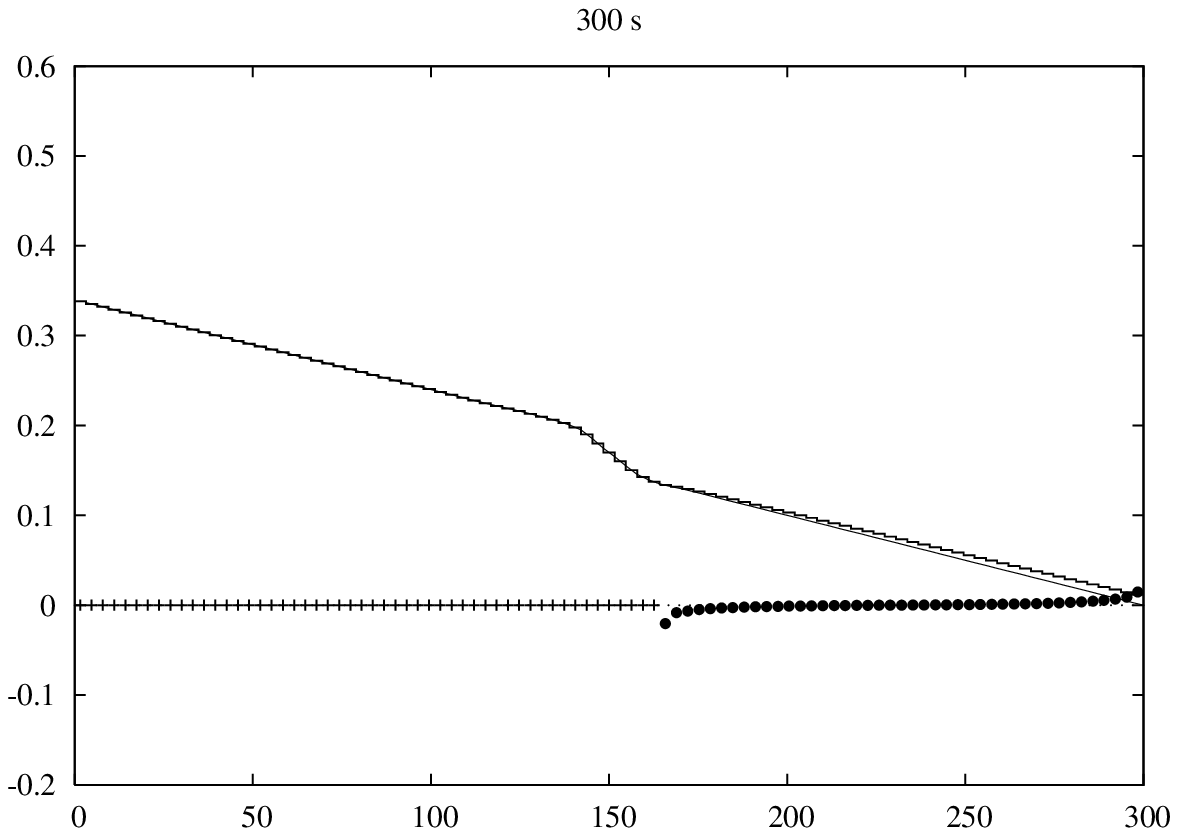}
\caption{TC1 - Free surface (solid line) and interface
normal velocity ($cm/12min$) plotted with circles if
interface is wet and with crosses if interface is
dry.\hfill\label{Fig:Couplage5}} 
\end{center}
\end{minipage}
\hspace{0.5cm}
\begin{minipage}[c]{6.5cm}
\begin{center}
\includegraphics[width=6.5cm]{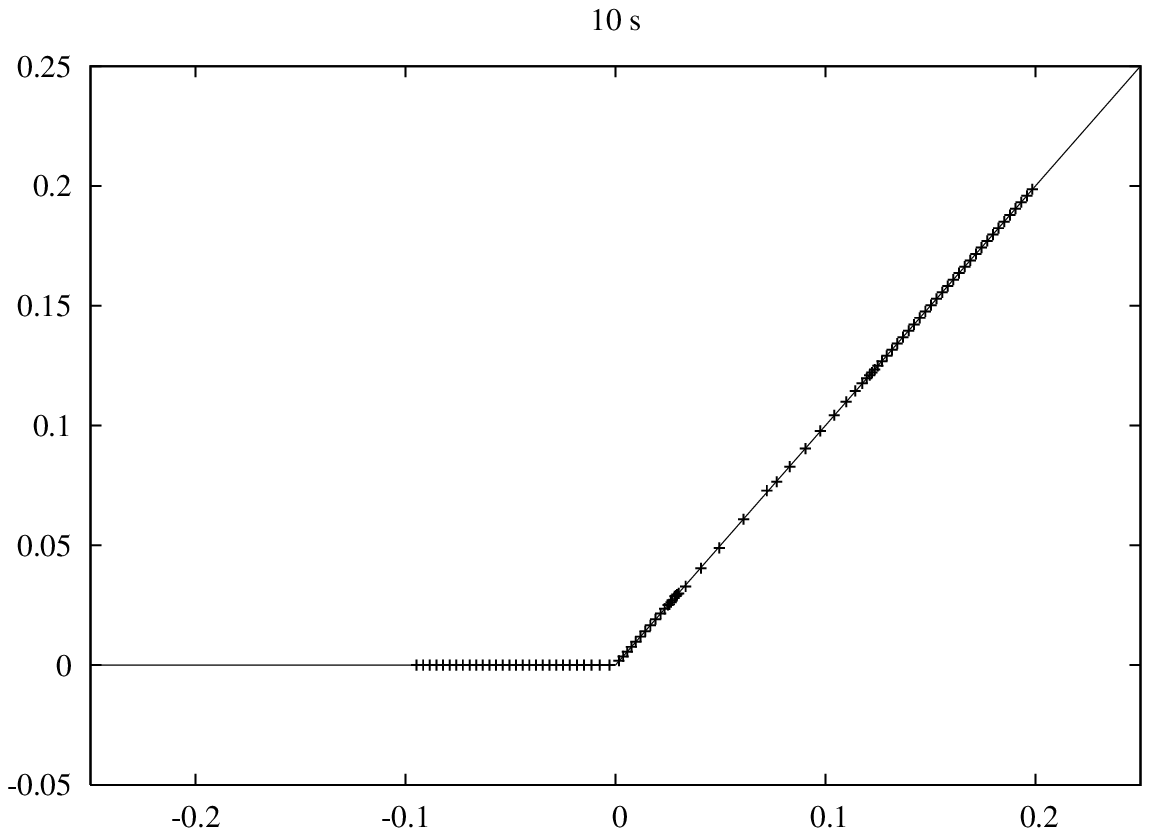}\\[6pt]
\includegraphics[width=6.5cm]{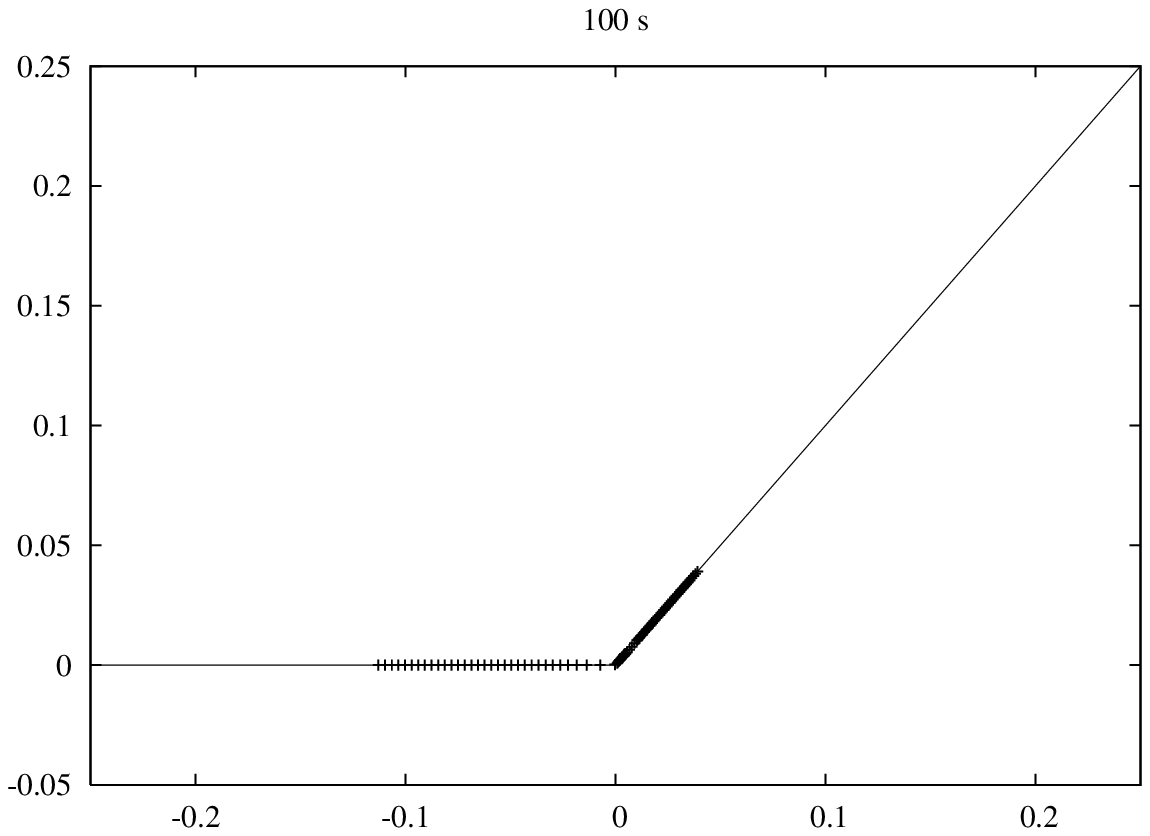}\\[6pt]
\includegraphics[width=6.5cm]{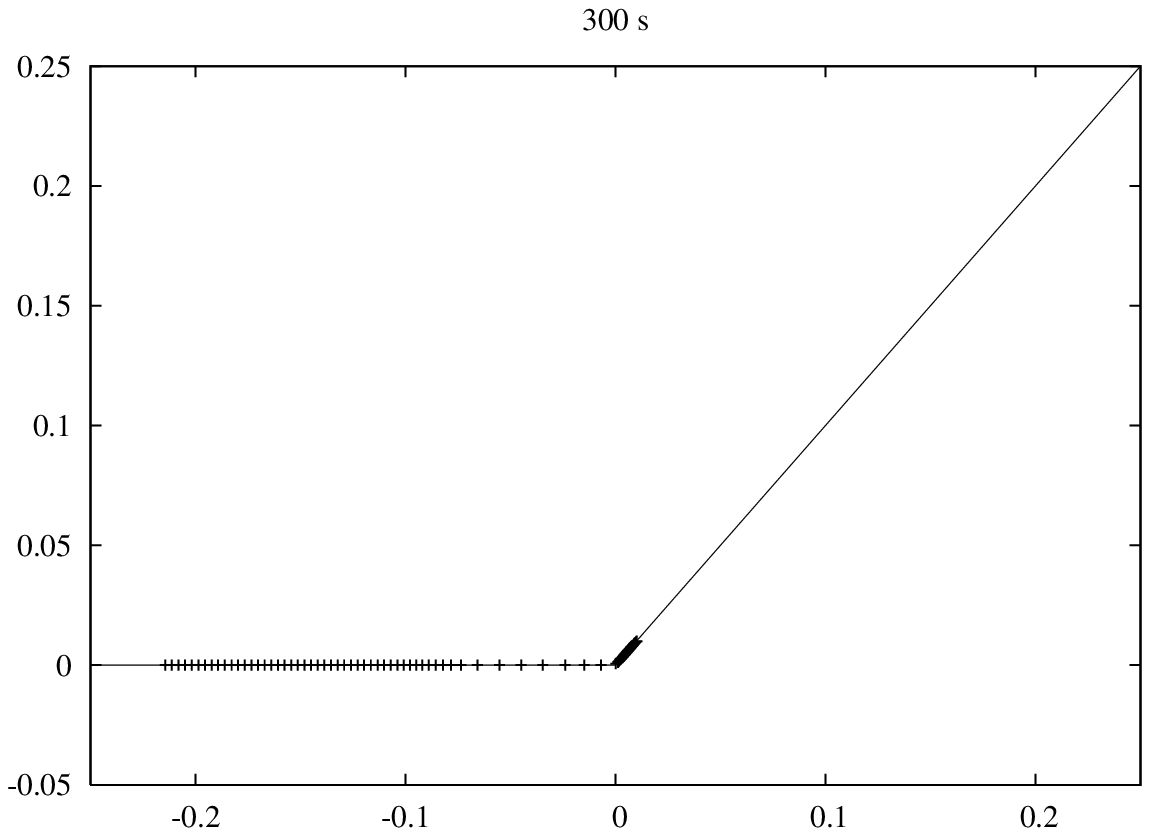}\vspace{0.3cm}
\caption{TC1 - Cloud of points ($\psi_{\h}^n,h_{\h}^n$) on
the admissible set $\A$ at different
times.\hfill\label{Fig:Couplage6}}
\vspace{0.8cm} 
\end{center}
\end{minipage}
\end{figure}

Figure~\ref{Fig:Couplage5} presents the free surface of the
overland flow ($h_{\h}^n$ + topography) and the interface normal
velocity $v^{\star, n}_{\h}$ along the interface at three
characteristic times of the simulation ($10s, 100s$ and
$300s$). The free surface being piecewise constant, it is
depicted on each interface cell by a solid line. The
interface normal velocity $v^{\star,n}_{\h}$ is plotted with circles if the
interface is wet (that is, on $\Ihw{n}$) and with crosses if
the interface is dry (that is, on $\Ihd{n}$).  

Figure~\ref{Fig:Couplage6} provides a closer insight at the
issue of staying on the admissible set $\A$. For the same
times as in Figure~\ref{Fig:Couplage5} and for each face of
$\mathcal{F}_{\h}^{\I}$, each couple
($\psi_{\h}^n,h_{\h}^n$) is represented by a cross (the
mean-value of $\psi_{\h}^n$ is considered on each face). The
admissible set $\A$ is also plotted with a solid line. 

The hydraulic jump in the overland flow is visible at the
beginning of the simulation on Figure~\ref{Fig:Couplage5} at
10s. Moreover, exfiltration appears on some faces located on
$\I_2$ and $\I_3$. During the simulation, a Neumann boundary
condition is imposed on the faces where the water becomes
equal to zero. It is confirmed by Figure~\ref{Fig:Couplage6}
where the number of couple ($\psi_{\h}^n,h_{\h}^n$) situated
on the branch $\{h=0\}$ increases. 

\subsection{Test case 2 (TC2)}
\label{sec:TC2}
\begin{figure}[htb]
\centering
\psset{unit=1}
\begin{pspicture}(-2,-0.5)(10,5)
%--- Boite ---------------------------
\psline[linestyle=dashed,linewidth=0.1mm](0,2.55)(8,2.55)
\psline[doubleline=true](0,3.1)(0,0)(8,0)(8,3)
\psline[linestyle=dashed](8,0)(10,0)
\psline(8,3)(0,3.1)
\psline[linewidth=0.1mm](0,4.5)(8,4.5)
\psline[linewidth=0.1mm]{->}(0,4.5)(0,4)
\psline[linewidth=0.1mm]{->}(0.5,4.5)(0.5,4)
\psline[linewidth=0.1mm]{->}(1,4.5)(1,4)
\psline[linewidth=0.1mm]{->}(1.5,4.5)(1.5,4)
\psline[linewidth=0.1mm]{->}(2,4.5)(2,4)
\psline[linewidth=0.1mm]{->}(2.5,4.5)(2.5,4)
\psline[linewidth=0.1mm]{->}(3,4.5)(3,4)
\psline[linewidth=0.1mm]{->}(3.5,4.5)(3.5,4)
\psline[linewidth=0.1mm]{->}(4,4.5)(4,4)
\psline[linewidth=0.1mm]{->}(4.5,4.5)(4.5,4)
\psline[linewidth=0.1mm]{->}(5,4.5)(5,4)
\psline[linewidth=0.1mm]{->}(5.5,4.5)(5.5,4)
\psline[linewidth=0.1mm]{->}(6,4.5)(6,4)
\psline[linewidth=0.1mm]{->}(6.5,4.5)(6.5,4)
\psline[linewidth=0.1mm]{->}(7,4.5)(7,4)
\psline[linewidth=0.1mm]{->}(7.5,4.5)(7.5,4)
\psline[linewidth=0.1mm]{->}(8,4.5)(8,4)
\psline[linewidth=0.1mm]{<-}(0,-0.5)(3.5,-0.5)
\psline[linewidth=0.1mm]{->}(4.5,-0.5)(8,-0.5)
\psline[linewidth=0.1mm]{<-}(8.5,0.0)(8.5,1.1)
\psline[linewidth=0.1mm]{->}(8.5,1.9)(8.5,3)
\psline[linewidth=0.1mm]{<-}(-0.5,0.0)(-0.5,0.9)
\psline[linewidth=0.1mm]{->}(-0.5,1.7)(-0.5,2.6)
%--- Les indic des points -------------
\rput(-0.6,1.3){$0.85m$}
\rput(4,-0.5){$6m$}
\rput(8.5,1.5){$1m$}
\rput(4,0.3){$v_N=0$}
\rput(4,3.5){$J=0.5\%$}
\rput(4,4.9){$|\vr \cdot n_\Omega|= 0.1K_s=3.6mmh^{-1}$}
\rput(7,2.75){$\bigtriangledown$}
\rput(9.5,0.3){$z=0$}
\end{pspicture}
\vspace{0.3cm}
\caption{TC2 - Geometry, initial water table position and constant rainfall intensity.\label{Fig:Couplage9}}
\end{figure}
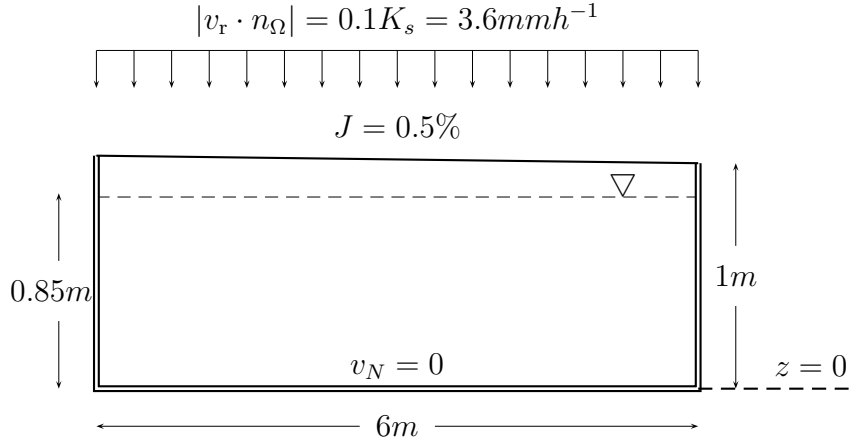

The principle of this test case is inspired by the work of
Abdul and Gilham \cite{Abdul84}: a constant rainfall
intensity is imposed on the upper part of the domain for a
fixed period of time, whereas the lateral and lower
boundaries are impermeable. In our case, the geometry is
shown in Figure~\ref{Fig:Couplage9} and the final simulation
time is $T = 360s$. The initial condition is an horizontal
water table located at $0.85m$ with an hydrostatic pressure
profile and the boundary condition on walls and bottom is a
zero flux, 
\begin{align*}
\psi^0 = - z + 0.85m &\quad \text{in}\ \Omega, \\
v_N = 0&\quad \text{on}\ (\W \cup \B) \times[0,T].
\end{align*}
For the overland flow, the initial condition and the
boundary condition are 
\begin{align*}
h^0\ = 0 &\quad \text{on}\ \I, \\
h_{\rm{A}}\ = 0 &\quad \text{at}\ \rm{A}\times[0,T].
\end{align*}
A constant rainfall intensity equal to 10$\%$ of the
hydraulic conductivity at saturation is imposed during $180s$
and is stopped afterwards, 
\begin{align*}
\vr\cdot n_\Omega = -0.1K_s&\quad \text{on}\ \I\times[0,180], \\
\vr\cdot n_\Omega = 0 &\quad      \text{on}\ \I\times[180,T].
\end{align*}

A mesh with 2049 triangles (corresponding to a typical
mesh-size of $10cm$) and time step $\delta t = \delta t' =
1s$ have been used. We have verified that the interface
normal velocity obtained with a finer mesh (8763 elements)
and a smaller time step ($\delta t = \delta t' = 0.5s$) can
be superimposed to that reported below. Also, observe that
$\delta t'=1s$ roughly corresponds to the CFL condition, so
that, for the present test case, the accuracy limit on the
time step for Richards' equation is comparable to the CFL
restriction.

\begin{figure}[htb]
\begin{center}
\includegraphics[width=6.5cm]{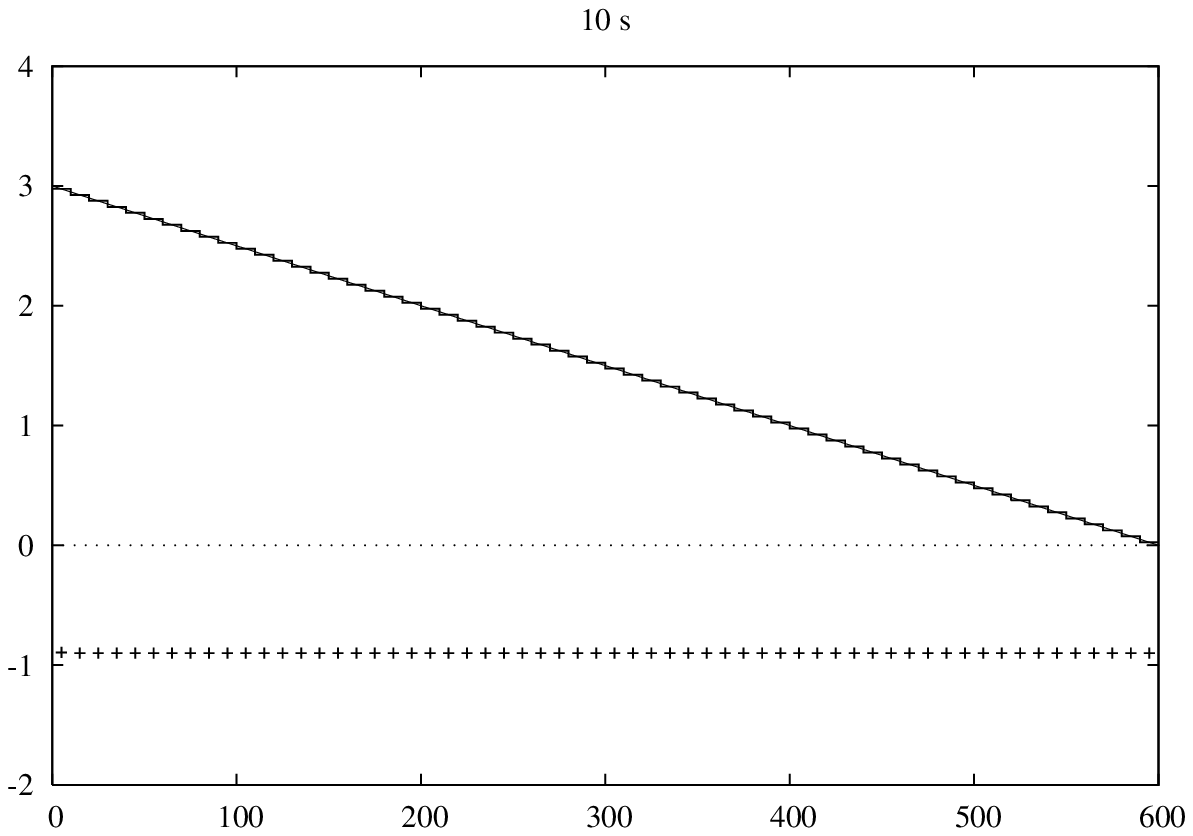}\hspace{3pt}
\includegraphics[width=6.5cm]{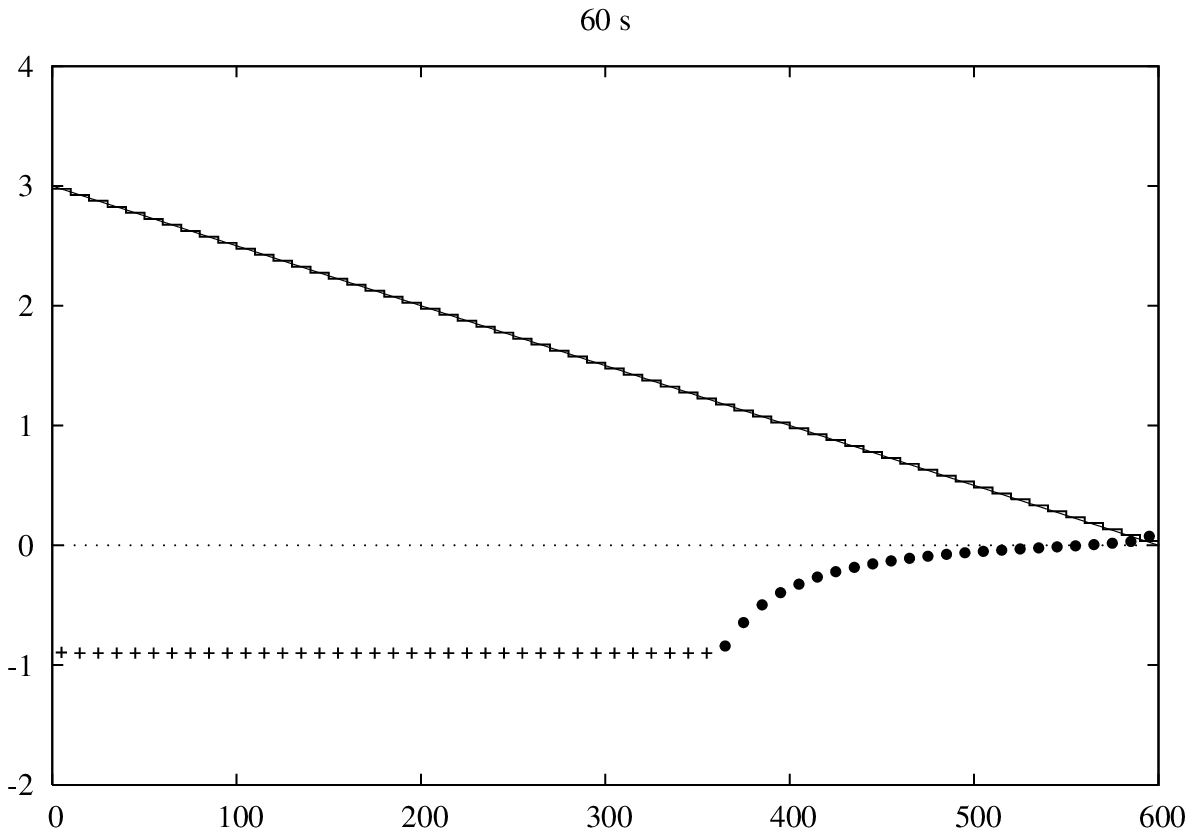}\\[6pt]
\includegraphics[width=6.5cm]{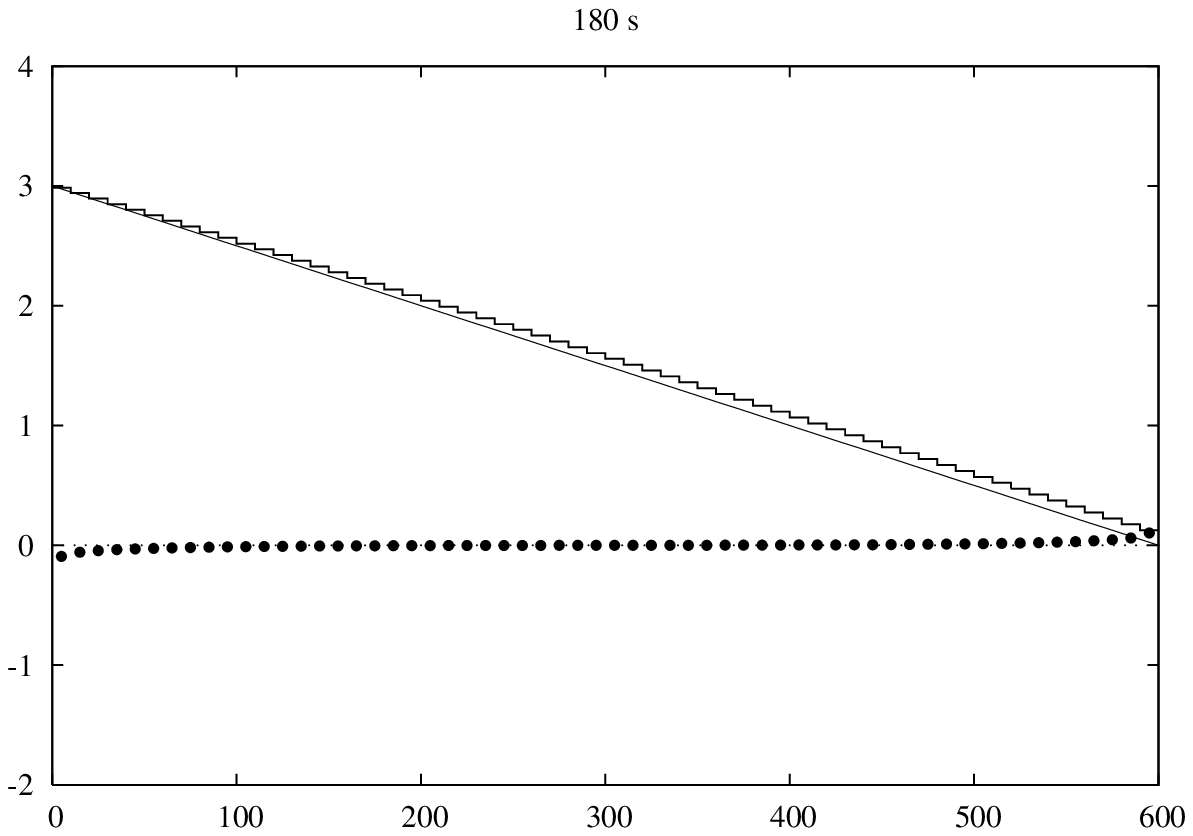}\hspace{3pt}
\includegraphics[width=6.5cm]{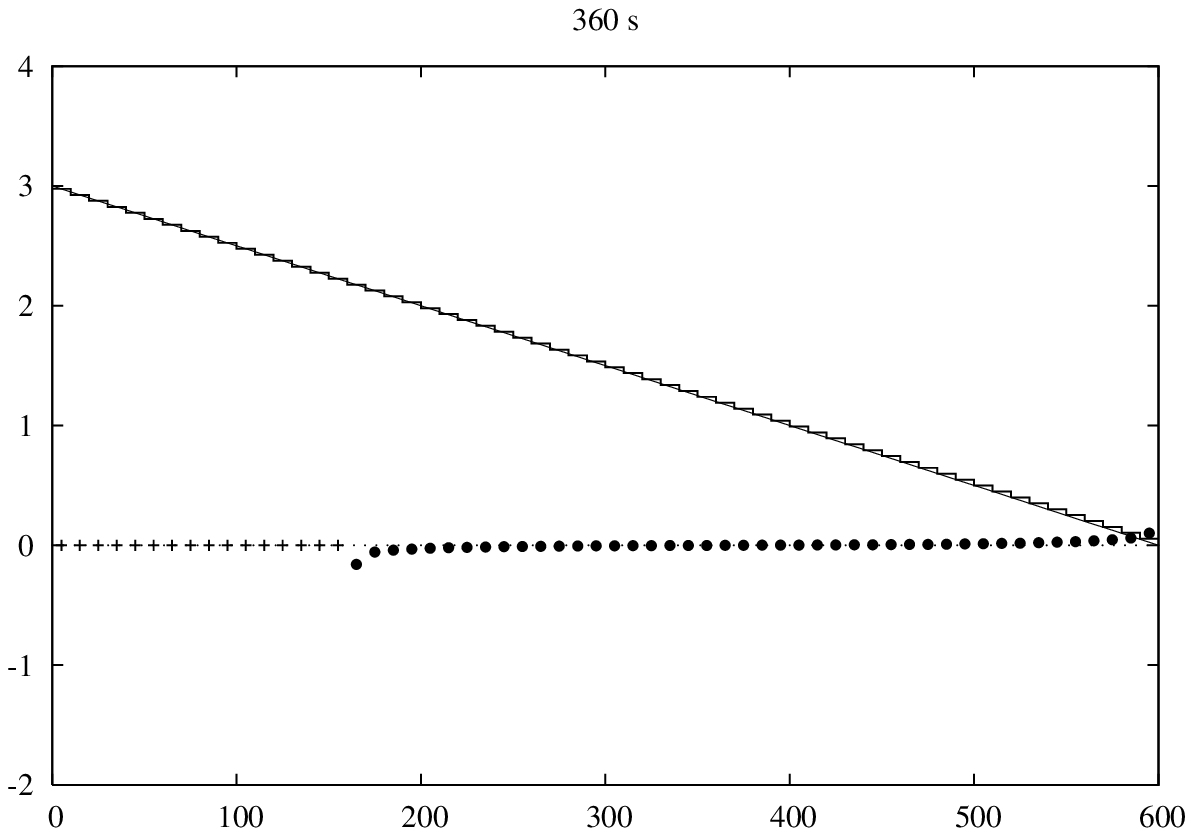}
\caption{TC2 - Free surface (solid line) and interface normal velocity ($cm/15min$) plotted with circles if interface is wet and with crosses if interface is dry.\label{Fig:Couplage10}}
\end{center}
\end{figure}

Figure~\ref{Fig:Couplage10} presents the free surface and
the normal velocity $v^{\star, n}_{\h}$ along the interface
at four characteristic times of the simulation ($10s, 60s,
180s$ and $360s$). The same notation is used as in
Figure~\ref{Fig:Couplage5}. The hydrological response of the
system can be divided into four phases. \\
\textit{1 - Soil saturation $[0,50s]$.} In this phase, which
results from the initial water table position, the $15cm$
top layer is being saturated. The rainfall is totally
absorbed by the soil, the surface water depth is equal to
zero and a Neumann condition is imposed on the all faces of
$\I$ (Figure~\ref{Fig:Couplage10} at $10s$).\\
\textit{2 - Surface runoff occurs on part of $\I$
$[50s,90s]$.} The rainfall is still partially absorbed by
the soil, but a Dirichlet condition is now being imposed on
the part of $\I$ located near the outlet. Interestingly,
infiltration occurs on the most part of the interface since
the normal velocity is negative, but some exfiltration
occurs on the first few faces located near the outlet where
the normal velocity becomes positive
(Figure~\ref{Fig:Couplage10} at $60s$).\\
\textit{3 - Surface runoff occurs on $\I$ $[90s,180s]$.}
Surface runoff occurs on the whole interface and the soil is
totally saturated. A Dirichlet condition is imposed
throughout the interface and the surface water depth is
positive (Figure~\ref{Fig:Couplage10} at $180s$). \\
\textit{4 - Drainage $[180s,360s]$.} When rainfall stops,
the surface water depth returns to zero on the faces located
near the point $\rm{A}$ because of infiltration and surface
runoff. A Neumann boundary condition is imposed on the dry
zone near the point $\rm{A}$ (Figure~\ref{Fig:Couplage10} at
$360s$).  

\begin{figure}[htb]
\begin{center}
\includegraphics[width=6.5cm]{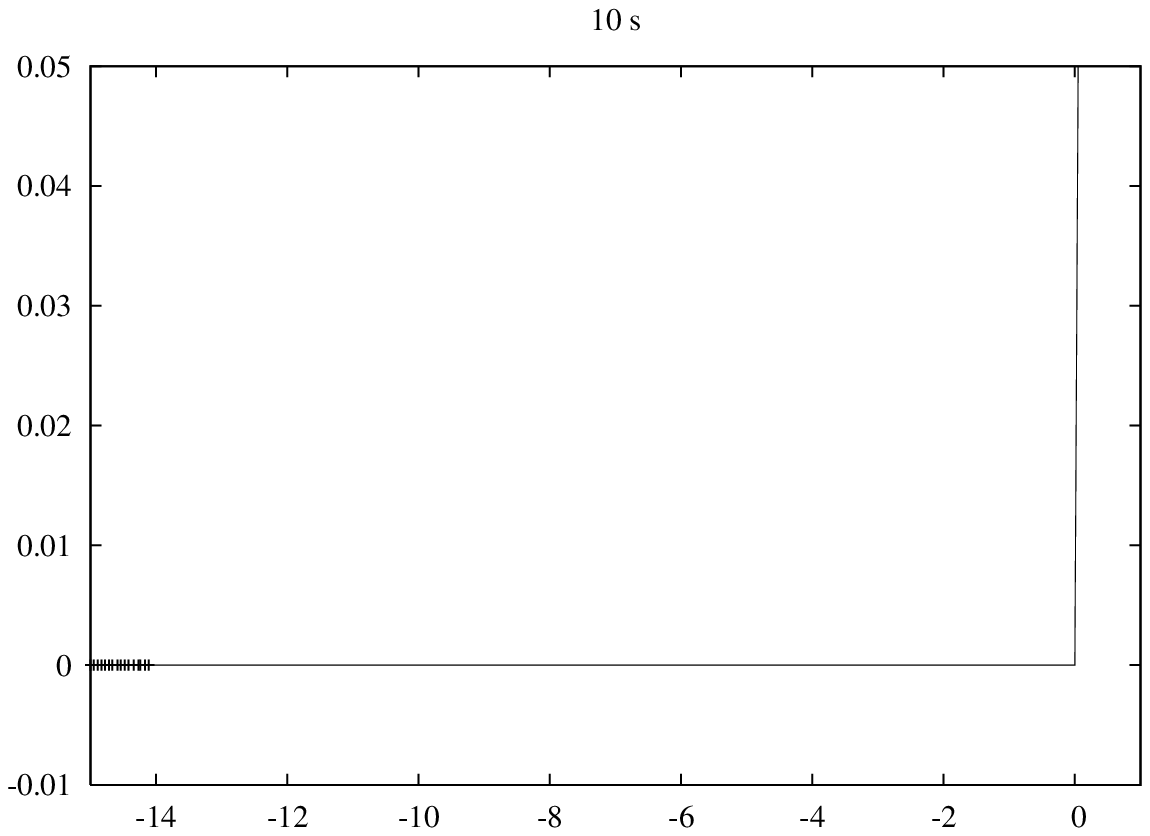}\hspace{3pt}
\includegraphics[width=6.5cm]{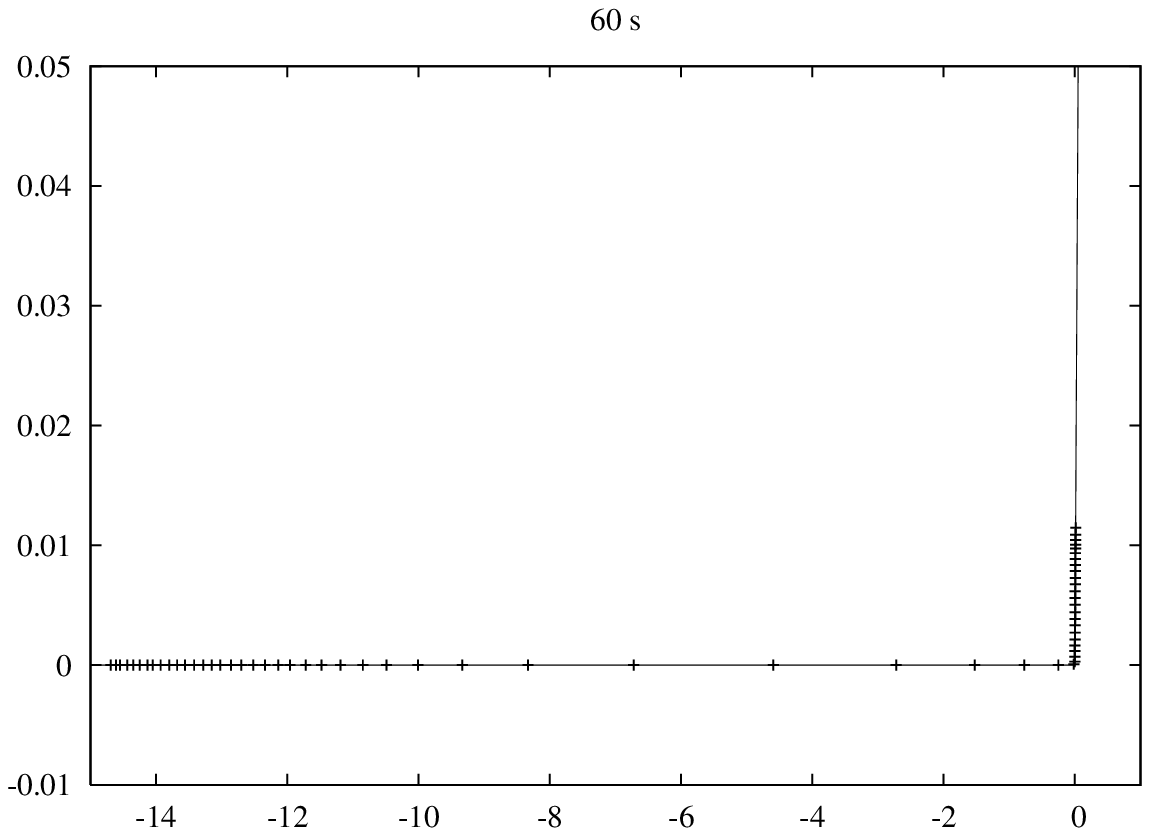}\\[6pt]
\includegraphics[width=6.5cm]{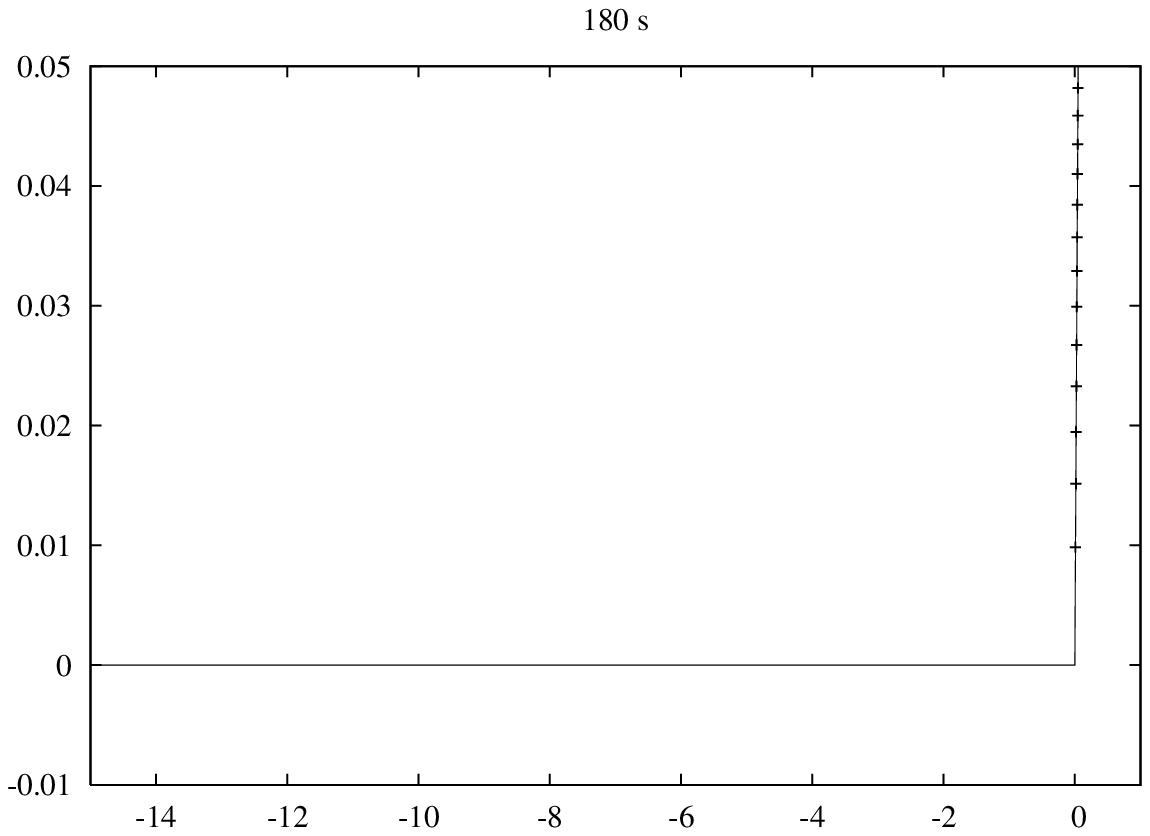}\hspace{3pt}
\includegraphics[width=6.5cm]{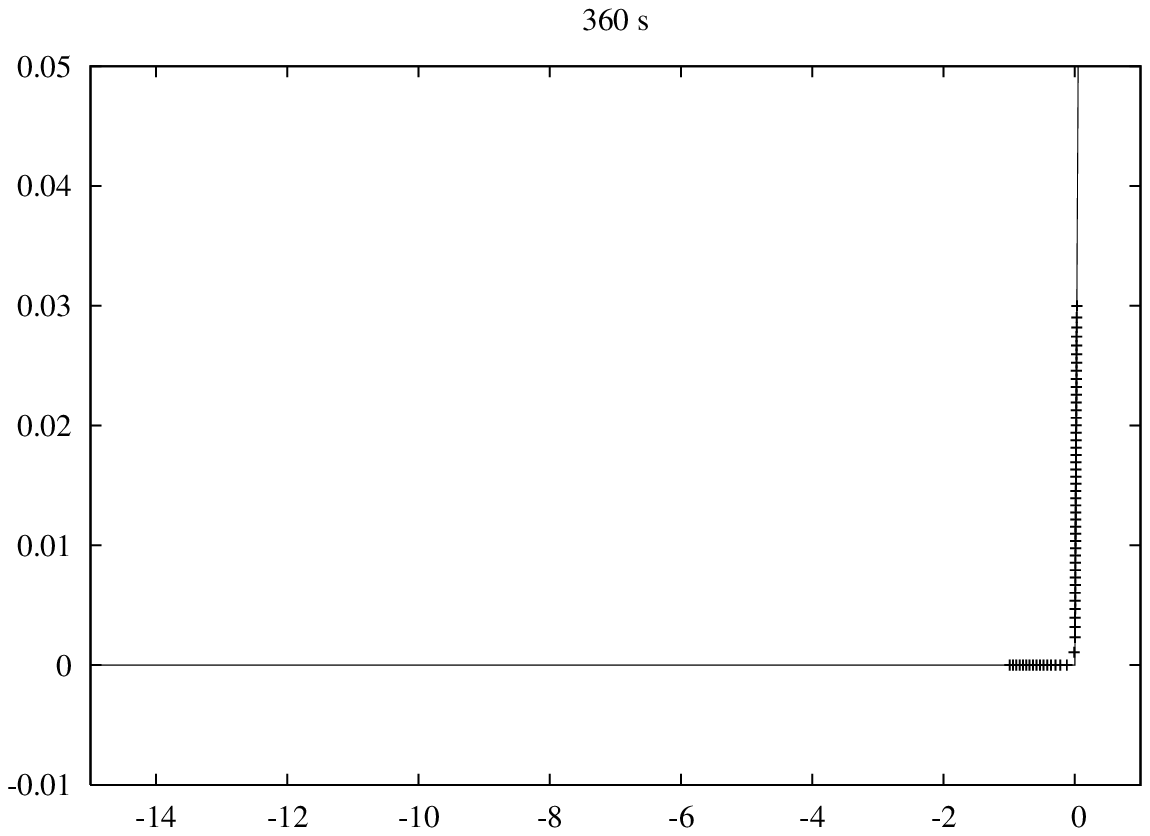}\vspace{-0.3cm}
\caption{TC2 - Cloud of points ($\psi_{\h}^n,h_{\h}^n$) on the admissible set $\A$ at different times.\label{Fig:Couplage11}}
\end{center}
\end{figure}

Figure~\ref{Fig:Couplage11} provides a closer insight at the
issue of staying on the admissible set $\A$. For the same
times as in Figure~\ref{Fig:Couplage10} and for each face of
$\mathcal{F}_{\h}^{\I}$, each couple
($\psi_{\h}^n,h_{\h}^n$) is represented as in
Figure~\ref{Fig:Couplage6}. Note that different scales are
used, so that the branch $\{h=\psi\}$ is almost
vertical. The four phases described above are clearly
illustrated by the position of the cloud of points. At
$10s$, the hydraulic head is negative and the water depth is
equal to zero. The cloud of points is only on the branch
$\{h=0\}$ corresponding to a dry soil. At $60s$, the
hydraulic head is equal to the water depth for some
faces. The cloud of points is located on the two branches
because the soil contains both saturated and unsaturated
zones. At $180s$, the hydraulic head is equal to the water
depth for all the faces. The cloud of points is only on the
branch $\{h=\psi\}$ corresponding to a wet soil. At $360s$,
the hydraulic head becomes again negative where the surface
water depth is equal to zero. The cloud of points is again
located on the two branches. 

\begin{figure}[htb]
\begin{center}
\psfrag{Mass}{\hspace{-0.2cm}\raisebox{0.1cm}{\footnotesize{Mass}}}
\psfrag{Time}{\hspace{-0.25cm}\raisebox{-0.2cm}{\footnotesize{Time}}}
\psfrag{a1}{\hspace{-0.6cm}\raisebox{-0.1cm}{\scriptsize{$\Sigma\M^n_{\rm{in}}$}}}
\psfrag{a2}{\hspace{-0.5cm}\raisebox{0.0cm}{\scriptsize{$\DM^n$}}}
\psfrag{a3}{\hspace{-0.81cm}\raisebox{-0.1cm}{\scriptsize{$\DM^n_{\rm{grnd}}$}}}
\psfrag{a4}{\hspace{-0.8cm}\raisebox{-0.03cm}{\scriptsize{$\DM^n_{\rm{over}}$}}}
\psfrag{a5}{\hspace{-0.7cm}\raisebox{-0.32cm}{\scriptsize{$\Sigma\M^n_{\rm{out}}$}}}
\includegraphics[width=6.5cm]{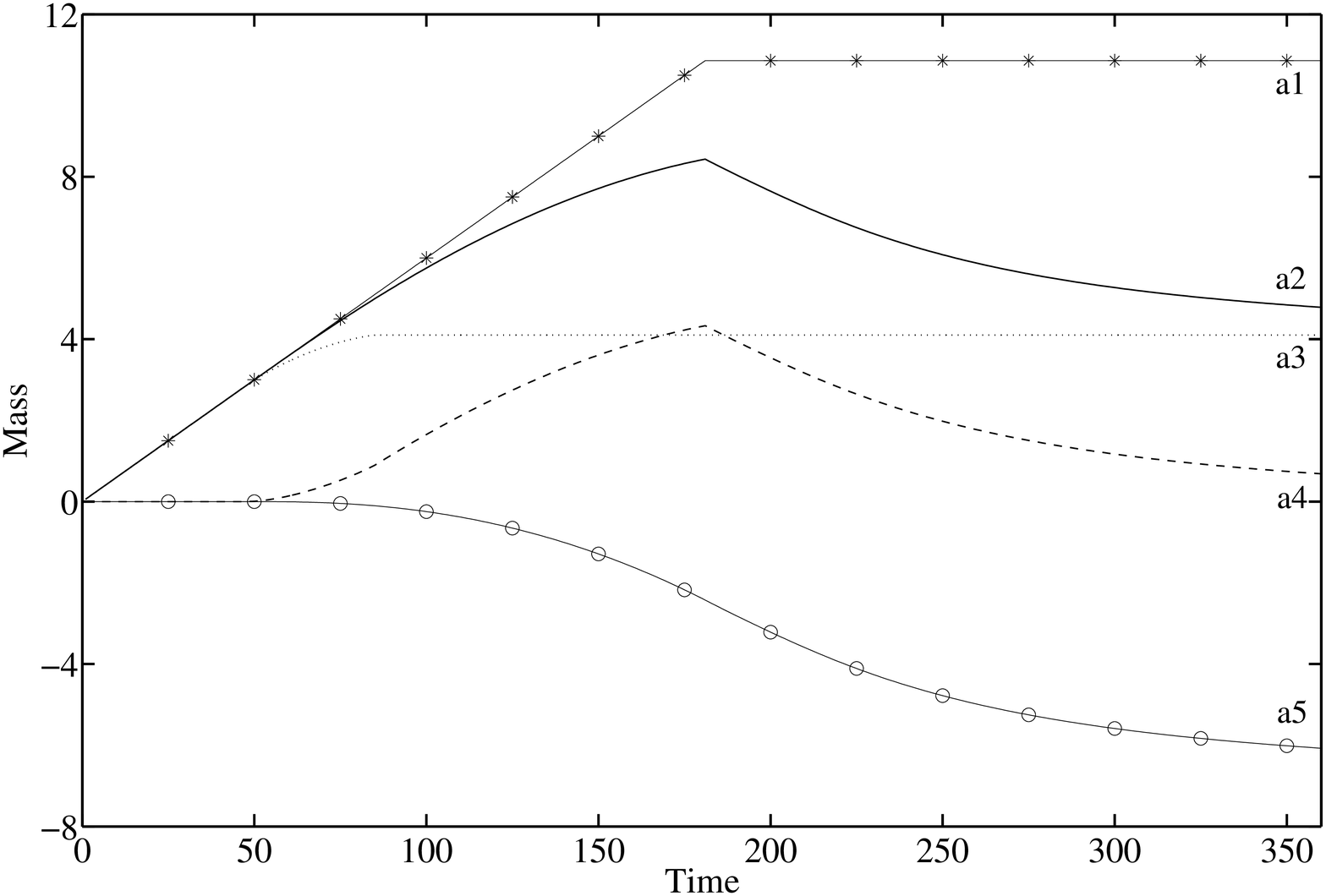}\hspace{0.2cm}
\includegraphics[width=6.5cm]{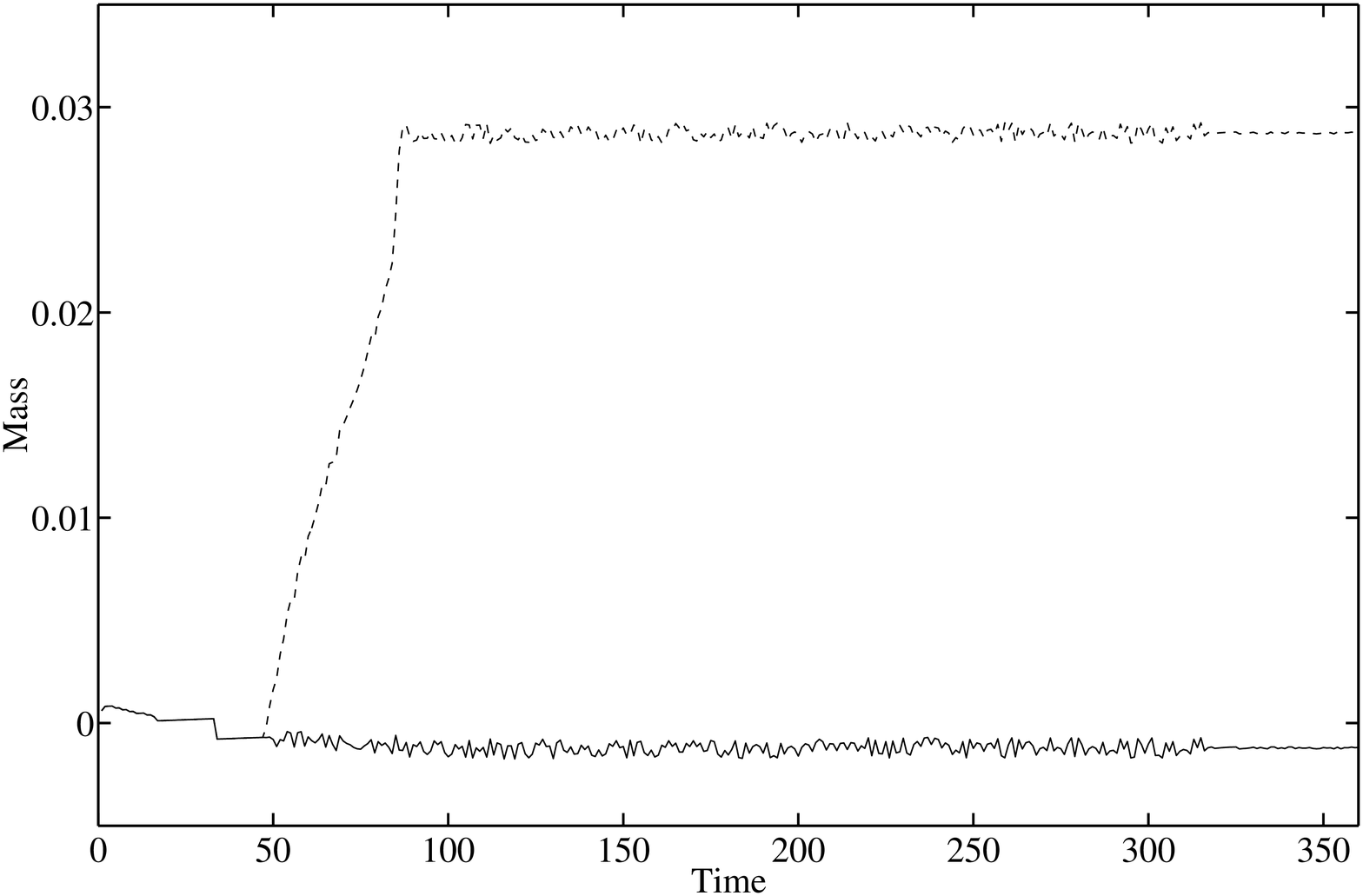}
\caption{TC2 - Left: Mass repartition in the coupled system; Right: Mass balance defect $\mathbb{E}^{n}$ for Algorithm~\ref{alg2} (dashed) and Algorithm~\ref{alg3} (solid).}\label{Fig:Couplage12}
\end{center}
\end{figure}

Figure~\ref{Fig:Couplage12} presents results on mass
conservation issues. Multiplying equation 
(\ref{Theo:MC2}) by the water density $\rho$, summing over
the time intervals in $[0,n\delta t]$, knowing that
$\tilde{F}_{\mathcal{WB}^n}=0$ and the definition of $V^n$,
$F_{\rm{ABr}}^n$ and $\Delta V^n$ yields 
$$
\underbrace{\frac{}{}\rho (V_{\rm{grnd}}^n -V_{\rm{grnd}}^0)}_{\DM^n_{\rm{grnd}}}+ 
\underbrace{\frac{}{}\rho (V_{\rm{over}}^n -V_{\rm{over}}^0)}_{\DM^n_{\rm{over}}}= 
  \sum_{i=1}^n \underbrace{\frac{}{}\rho \delta t  \big( F_{\rm{A}}^{i} + F_{\rm{r}}^{i} \big)}_{\M^i_{\rm{in}}} 
+ \sum_{i=1}^n \underbrace{\frac{}{}\rho\delta t F_{\rm{B}}^{i}}_{\M^i_{\rm{out}}} 
+ \underbrace{\frac{}{}\rho\Delta V^n}_{\mathbb{E}^n}. $$
$\DM^n$ is defined as the total mass variation over the time
interval $[0,n\delta t]$ and is the sum of the total
groundwater mass variation $\DM^n_{\rm{grnd}}$ and the total
overland mass variation $\DM^n_{\rm{over}}$. The quantities
$\sum_{i=1}^n \M^i_{\rm{in}}$, $\sum_{i=1}^n
\M^i_{\rm{out}}$ and $\mathbb{E}^n$ are respectively the
total inflow of water, the total outflow of water and the
total mass balance defect cumulated at time $n\delta t$. The
five quantities $\DM^n, \DM^n_{\rm{grnd}},
\DM^n_{\rm{over}}, \sum_{i=1}^n \M^i_{\rm{in}}$ and
$\sum_{i=1}^n \M^n_{\rm{out}}$ are presented in the left
part of Figure~\ref{Fig:Couplage12}. In particular, this
figure confirms the four phases of the simulation. The
rainfall is totally absorbed by the soil at the 
beginning of the simulation until $50s$ since $\DM^n =
\DM^n_{\rm{grnd}}$. Then, the increase of
$\DM^n_{\rm{grnd}}$ diminishes and $\DM^n_{\rm{over}}$
becomes positive as a result of soil saturation. From 90s to
the end of the simulation, the variations of $\DM^n$ and
$\DM^n_{\rm{over}}$ are the same, corresponding to a
complete saturation of the soil. Moreover, during the last
phase, the total water inflow is constant because the
rainfall stops, so that the total water outflow is the same
as the total mass variation.  
 
Both total mass balance defects obtained with the
single-step and two-step coupling algorithms are compared in
the right part of Figure~\ref{Fig:Couplage12}. While
Algorithm~\ref{alg3} yields a sizable improvement over
Algorithm~\ref{alg2}, it can still be noticed that the mass
balance defect produced by Algorithm~\ref{alg2} is only of
the order of a few percent of the global quantities such as
$\DM^n$.  

\begin{figure}[htb]
\begin{center}
\psfrag{Mass}{\hspace{-0.2cm}\raisebox{0.1cm}{\footnotesize{Mass}}}
\psfrag{Time}{\hspace{-0.25cm}\raisebox{-0.2cm}{\footnotesize{Time}}}
\psfrag{a1}{\hspace{0.2cm}\raisebox{0.1cm}{\scriptsize{$\M^n_{\rm{in}}$}}}
\psfrag{a2}{\hspace{0.cm}\raisebox{0.1cm}{\scriptsize{$\rho\delta t \tilde{F}_{\I}^{+}$}}}
\psfrag{a3}{\hspace{0.cm}\raisebox{-0.1cm}{\scriptsize{$\rho\delta t \tilde{F}_{\I}^{-}$}}}
\psfrag{a4}{\hspace{-0.2cm}\raisebox{-0.1cm}{\scriptsize{$\M^n_{\rm{out}}$}}}
\includegraphics[width=7.5cm,height=6.6cm]{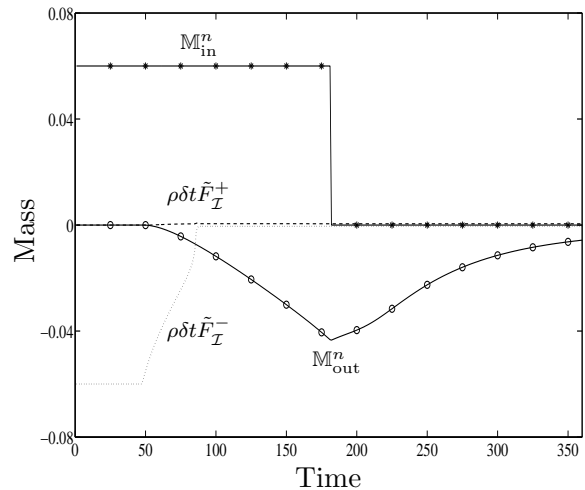}
\caption{TC2 - Mass fluxes in the kinematic wave equation.\label{Fig:Couplage13}}
\end{center}
\end{figure}

Finally, Figure~\ref{Fig:Couplage13} studies in more detail
the mass fluxes in the kinematic wave equation. The mass
flux $\tilde{F}_{\I}^n$ is decomposed into the exfiltration
flux $\tilde{F}_{\I}^{n,+}$ and the infiltration flux
$\tilde{F}_{\I}^{n,-}$ in the form $\tilde{F}_{\I}^n \eqdef
\tilde{F}_{\I}^{n,+} + \tilde{F}_{\I}^{n,+}$, with 
$\tilde{F}_{\I}^{n,+} \eqdef -\int_{\Ih^{n,+}}
v^{\star,n}_{\h}$ and $\tilde{F}_{\I}^{n,-} \eqdef
-\int_{\Ih^{n,-}} v^{\star,n}_{\h}$, where the
time-dependent sets $\Ih^{n,+}$ and $\Ih^{n,-}$ are defined
as follows 
$$\Ih^{n,+} \eqdef \{x\in \I;\ v_{\h}^{\star,n}(x) <0\}\quad\text{and}\quad
  \Ih^{n,-} \eqdef \{x\in \I;\ v_{\h}^{\star,n}(x) >0\}. $$
The four quantities $\rho\delta t \tilde{F}_{\I}^{n,+},
\rho\delta t \tilde{F}_{\I}^{n,-}, \M^n_{\rm{in}}$ and
$\M^n_{\rm{out}}$ are plotted on Figure~\ref{Fig:Couplage13}
as a function of time. 

\subsection{Test case 3 (TC3)}
\label{sec:TC3}
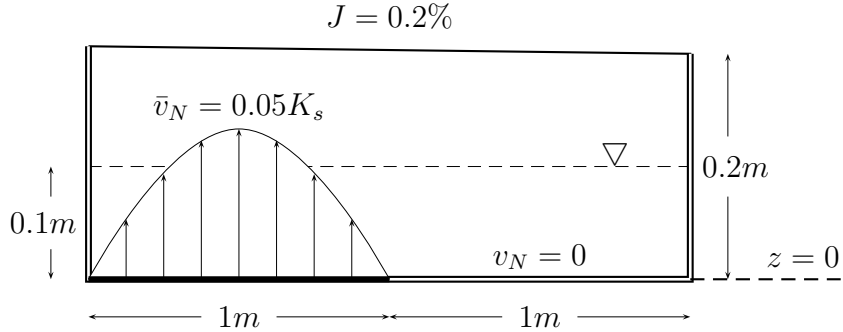
\begin{figure}[htb]
\centering
\psset{unit=1}
\begin{pspicture}(-2,-0.5)(10,3.8)
%--- Boite ---------------------------
\psline(0,3.1)(0,0)(8,0)(8,3)(0,3.1)
\psline[linestyle=dashed,linewidth=0.1mm](0,1.5)(8,1.5)
\psline[doubleline=true](0,3.1)(0,0)(8,0)(8,3)
\psline[linestyle=dashed](8,0)(10,0)
\pscurve[fillstyle=solid,fillcolor=white,linewidth=0.1mm](0,0)(2,2)(4,0)
\psline(0,0)(4,0)
\psline[linewidth=0.1mm]{->}(0.5,0)(0.5,0.8)
\psline[linewidth=0.1mm]{->}(1,0)(1,1.4)
\psline[linewidth=0.1mm]{->}(1.5,0)(1.5,1.85)
\psline[linewidth=0.1mm]{->}(2,0)(2,2)
\psline[linewidth=0.1mm]{->}(2.5,0)(2.5,1.85)
\psline[linewidth=0.1mm]{->}(3,0)(3,1.4)
\psline[linewidth=0.1mm]{->}(3.5,0)(3.5,0.8)
\psline[linewidth=0.8mm](0,0)(4,0)
\psline[linewidth=0.1mm]{<-}(0,-0.5)(1.5,-0.5)
\psline[linewidth=0.1mm]{->}(2.5,-0.5)(4,-0.5)
\psline[linewidth=0.1mm]{<-}(4,-0.5)(5.5,-0.5)
\psline[linewidth=0.1mm]{->}(6.5,-0.5)(8,-0.5)
\psline[linewidth=0.1mm]{<-}(8.5,0.0)(8.5,1.1)
\psline[linewidth=0.1mm]{->}(8.5,1.9)(8.5,3)
\psline[linewidth=0.1mm]{<-}(-0.5,0.0)(-0.5,0.4)
\psline[linewidth=0.1mm]{->}(-0.5,1.1)(-0.5,1.5)
%--- Les indic des points -------------
\rput(4,3.5){$J=0.2\%$}
\rput(2,-0.5){$1m$}
\rput(6,-0.5){$1m$}
\rput(8.6,1.5){$0.2m$}
\rput(-0.6,0.75){$0.1m$}
\rput(2,2.3){$\bar{v}_N=0.05K_s$}
\rput(6,0.3){$v_N=0$}
\rput(7,1.7){$\bigtriangledown$}
\rput(9.5,0.3){$z=0$}
\end{pspicture}
\caption{TC3 - Geometry, initial water table position and flux infiltration in groundwater.\label{Fig:Couplage14}}
\end{figure}

In this third test case, an exfiltration is produced on the
upper part of the domain by injecting water at the
bottom-left part. The geometry is presented in
Figure~\ref{Fig:Couplage14} and the final simulation time is
$T=360s$. The initial condition is an horizontal water table
located at $0.1m$ with an hydrostatic pressure profile and
the boundary condition on the walls corresponds to a zero
flux, 
\begin{align*}
\psi^0 = - z + 0.1m &\quad \text{in}\ \Omega, \\
v_N = 0 &\quad \text{on}\ \W \times[0,T].
\end{align*}
The rainfall intensity is set to zero. An infiltration flux
with a parabolic profile and a mean-value $\bar{v}_N$ equal
to $5\%$ of hydraulic conductivity at saturation is imposed
during $2$ minutes on the left half $\B_l$ of the bottom
($\B_l=\{(x,z) \in \B ,x\in [0,1] \}$ and $\B_r=\{(x,z) \in
\B, x\in [1,2]\}$). This infiltration flux is linear during
the first $10s$, constant on $[10,120]$, and equal to zero
for $t>120s$:
$$v_N(x,t) =
\left\{
\begin{array} {ll} 
x(x-1)\ 0.003K_s\ t, & \text{  if}\ (x,t)\in\B_l\times[0,10],\\
x(x-1)\ 0.03K_s,     & \text{  if}\ (x,t)\in\B_l\times[10,120],\\
0,                   & \text{  otherwise}.
\end{array} 
\right.$$
For the overland flow, the initial condition and the boundary condition are
\begin{align*}
h^0\ = 0 &\quad \text{on}\ \I, \\
h_{\rm{A}}\ = 0 &\quad \text{at}\ \rm{A}\times[0,T].
\end{align*}
A mesh with 1874 triangles (corresponding to a typical
mesh-size of $2.5cm$) and time step $\delta t = \delta t' =
1s$ have been used. We have verified that the interface
normal velocity obtained with a finer mesh (7310 elements)
and a smaller time step ($\delta t = \delta t' = 0.5s$) can
be superimposed to that reported below.  

\begin{figure}[tp]
\begin{center}
\includegraphics[width=6.5cm]{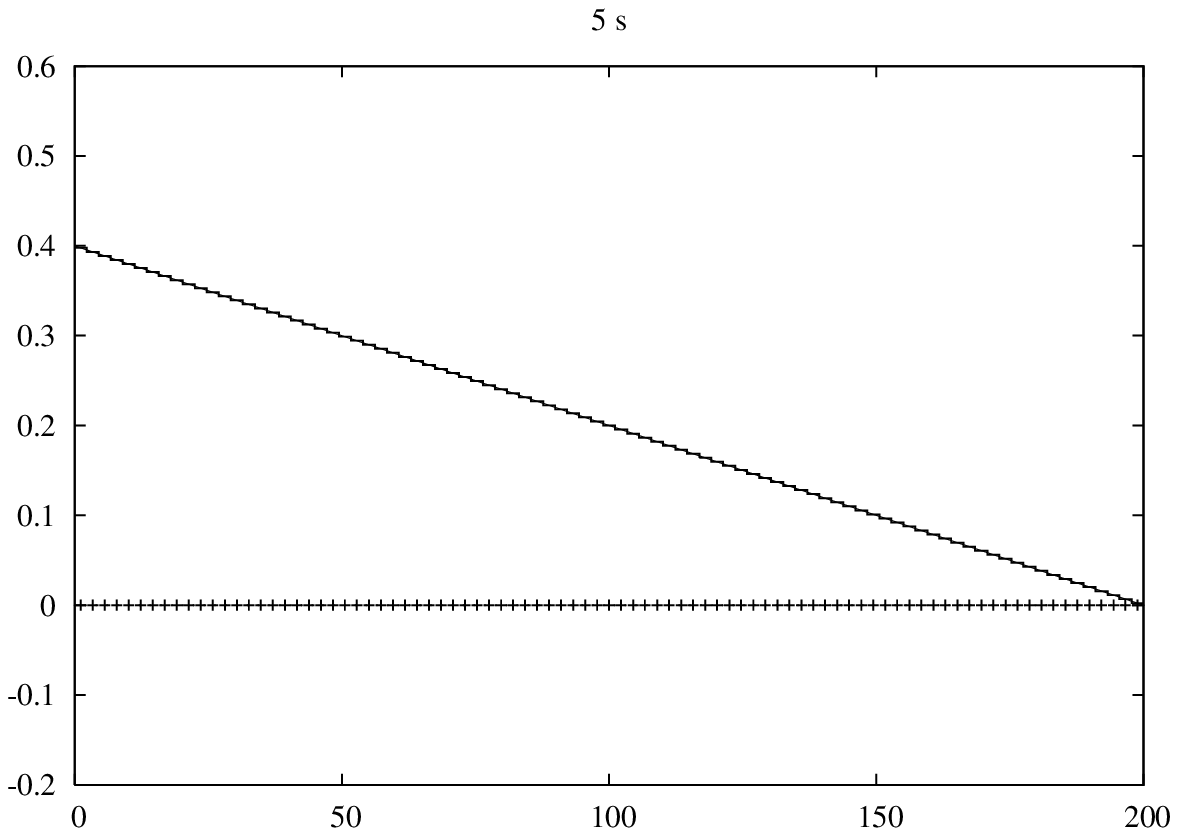}\hspace{3pt}
\includegraphics[width=6.5cm]{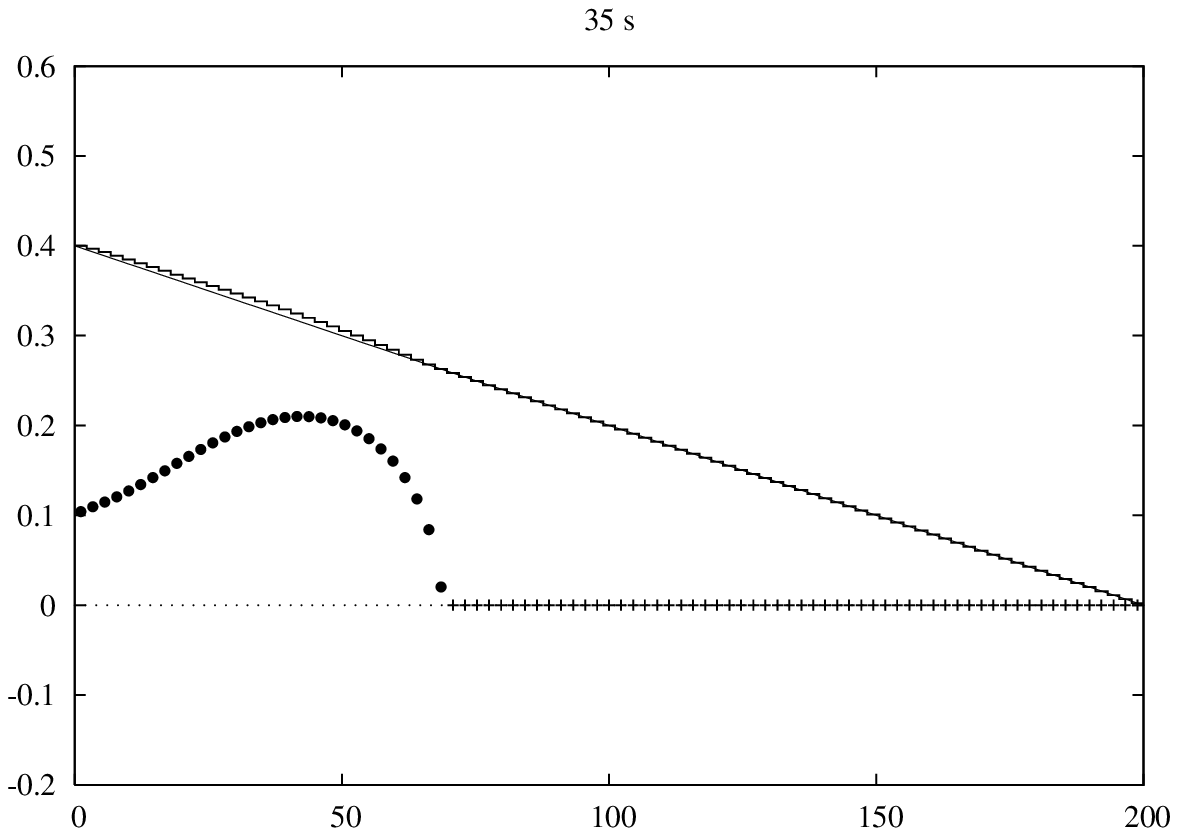}\\[6pt]
\includegraphics[width=6.5cm]{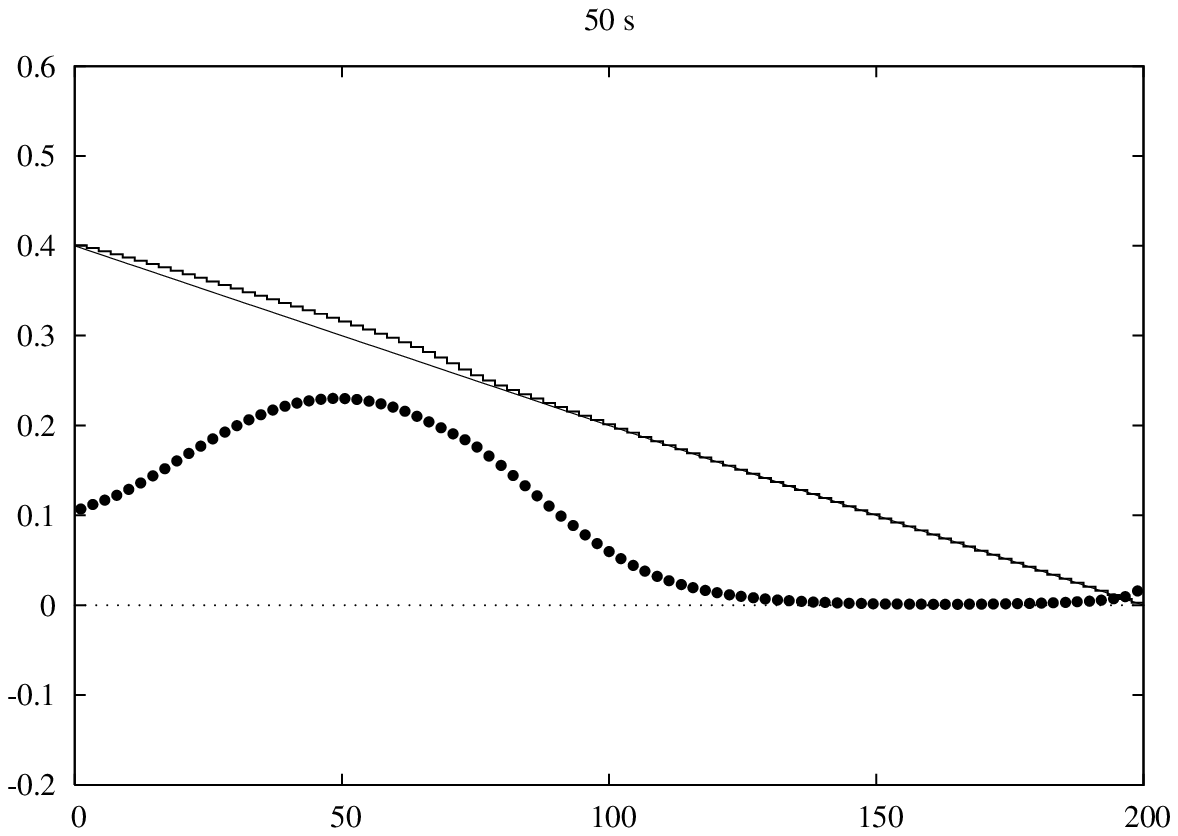}\hspace{3pt}
\includegraphics[width=6.5cm]{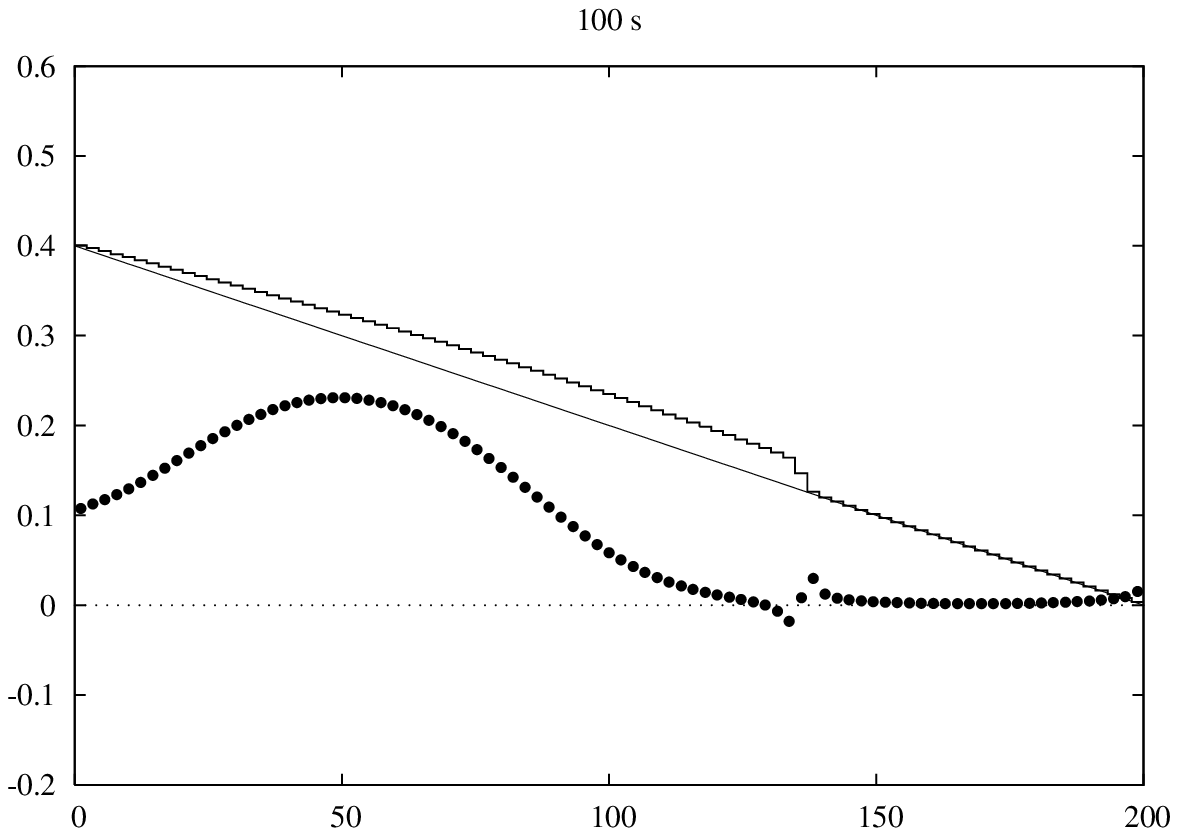}\\[6pt]
\includegraphics[width=6.5cm]{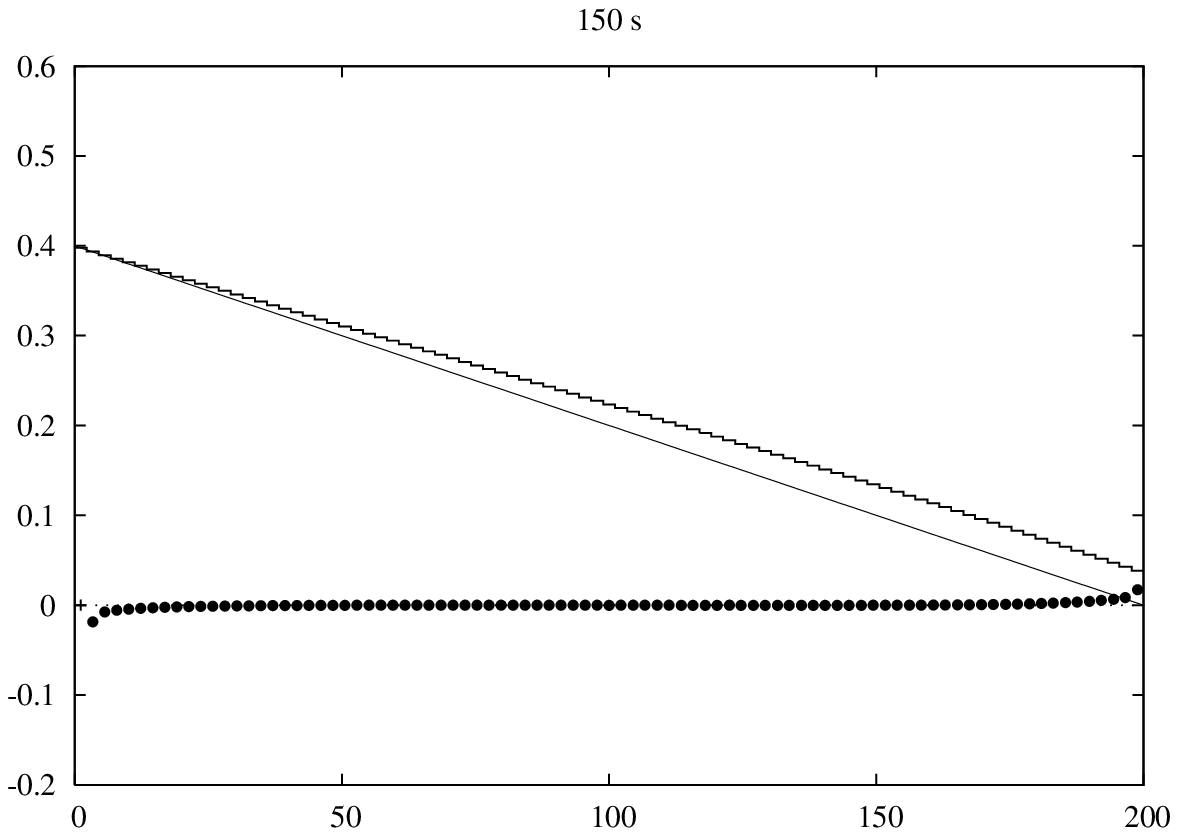}\hspace{3pt}
\includegraphics[width=6.5cm]{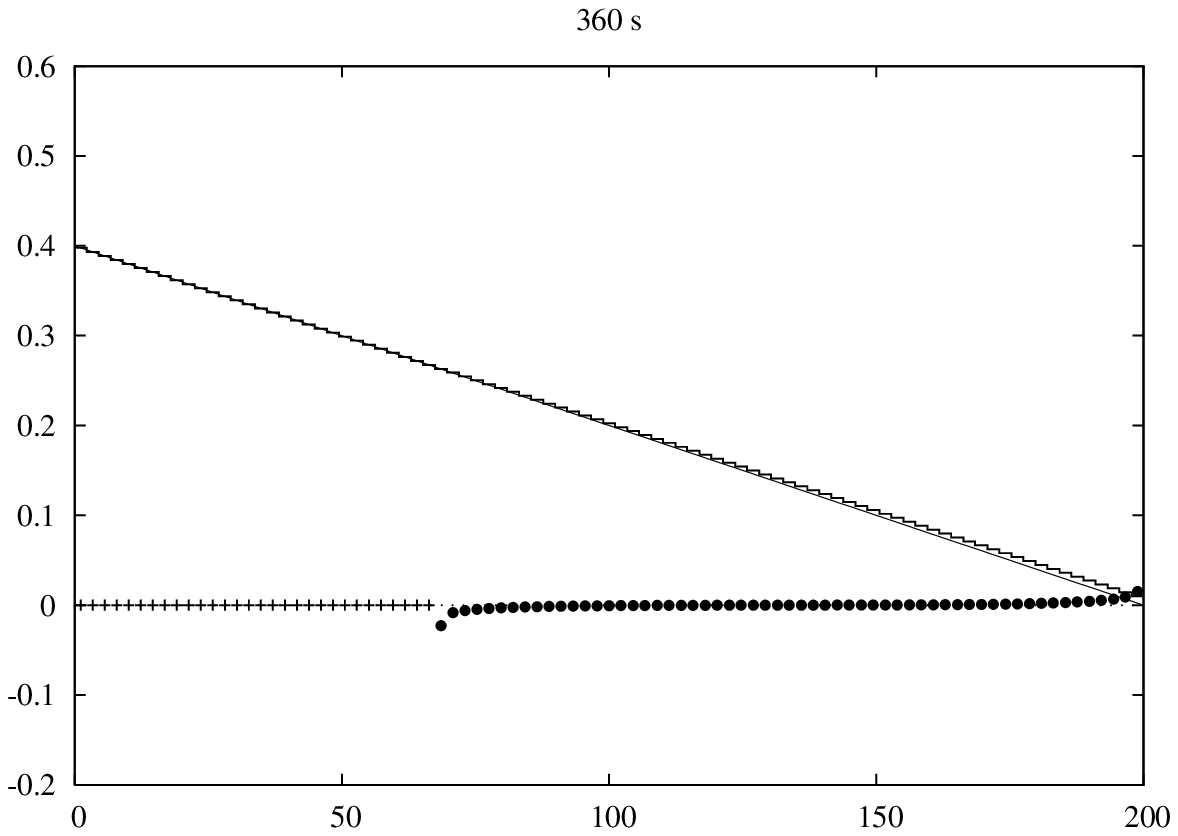}
\caption{TC3 - Free surface (solid line) and interface normal velocity ($cm/6min$) plotted with circles if interface is wet and with crosses if interface is dry.\label{Fig:Couplage15}}
\end{center}
\end{figure}

\begin{figure}[tp]
\begin{center}
\psfrag{Surface water depth}{\hspace{-1.cm}\raisebox{0.1cm}{\footnotesize{Water depth $h_{\h}^n$ (cm)}}}
\psfrag{x}{\hspace{0.cm}\raisebox{0.cm}{\footnotesize{x}}}
\psfrag{t1}{\hspace{-0.3cm}\raisebox{0.3cm}{\scriptsize{35s}}}
\psfrag{t2}{\hspace{-0.1cm}\raisebox{0cm}{\scriptsize{50s}}}
\psfrag{t3}{\hspace{-0.2cm}\raisebox{0cm}{\scriptsize{100s}}}
\psfrag{t4}{\hspace{-0.3cm}\raisebox{0cm}{\scriptsize{150s}}}
\psfrag{t5}{\hspace{-0.cm}\raisebox{0cm}{\scriptsize{360s}}}
\includegraphics[width=6.5cm]{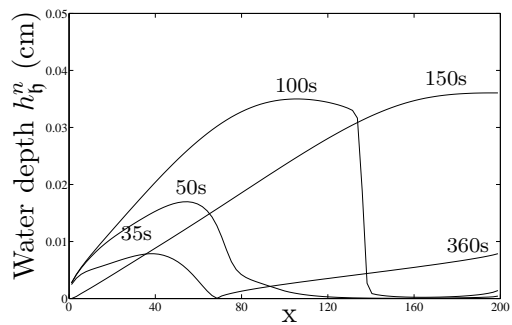}
\caption{TC3 - Surface water depth $h_{\h}^n$ at different times.\label{Fig:Couplage14bis}}
\end{center}
\end{figure}

Figure~\ref{Fig:Couplage15} presents the free surface and
the normal velocity $v^{\star, n}_{\h}$ along the interface
at six characteristic times of the simulation ($5s, 35s,
50s, 100s, 150s$ and $360s$) and
Figure~\ref{Fig:Couplage14bis} presents the surface water 
depth $h_{\h}^n$ at these different times. The same notation
is used as in Figure~\ref{Fig:Couplage5}. The hydrological
response of the system can be divided into six phases. 

\textit{1 - Soil saturation $[0,15s]$.} This phase results
from the initial water table position. At the beginning of
the simulation, the soil is partially saturated and 
the injection at the bottom of the domain increases the
hydraulic head. The interface is totally dry and a zero
Neumann boundary condition is enforced everywhere on $\I$
(Figure~\ref{Fig:Couplage15} at $5s$). \\
\textit{2 - Partial exfiltration $[15s,45s]$.} The soil
becomes saturated in the left part of the domain and the
interface normal velocity positive, so that water begins
to exfiltrate from the faces situated in this saturated
zone. A Dirichlet condition is enforced on those faces
(Figure~\ref{Fig:Couplage15} at $35s$). \\ 
\textit{3 - Full exfiltration $[45,100s]$.} When the soil is
totally saturated, the amount of exfiltrated water is equal
to the amount of injected water. We observe that overland
flow occurs over the whole interface $\I$ and that a
Dirichlet condition is being enforced everywhere. However,
most of the overland flow still remains concentrated near
the upper part of the interface (Figure~\ref{Fig:Couplage15}
at $50s$).\\
\textit{4 - Propagation of the runon wave $[100,120s]$.} In this
phase, the runon wave propagates downstream. It is
worthwhile to notice that a slight part of the surface water
infiltrates back into the soil as indicated by the sign of
the normal velocity near the heading part of the runon wave
(Figure~\ref{Fig:Couplage15} at $100s$).\\
\textit{5 - Outflow $[120,200s]$.} When water injection
ceases at the bottom of the domain, the amount of
exfiltrated water decreases sharply and there is even a
small portion of the interface located near the point
$\rm{A}$ where water infiltrates back into the soil (despite
the boundary condition is of Dirichlet type since the
surface water depth is still positive) while most of the
overland flow reaches the outlet and exits the system
(Figure~\ref{Fig:Couplage15} at $150s$). \\
\textit{6 - Drainage $[200,360s]$.}  The surface water depth
vanishes on the upper part of the interface $\I$ where a
zero Neumann condition is now imposed
(Figure~\ref{Fig:Couplage15} at $360s$). 
\begin{figure}[b]
\begin{center}
\includegraphics[width=6.5cm]{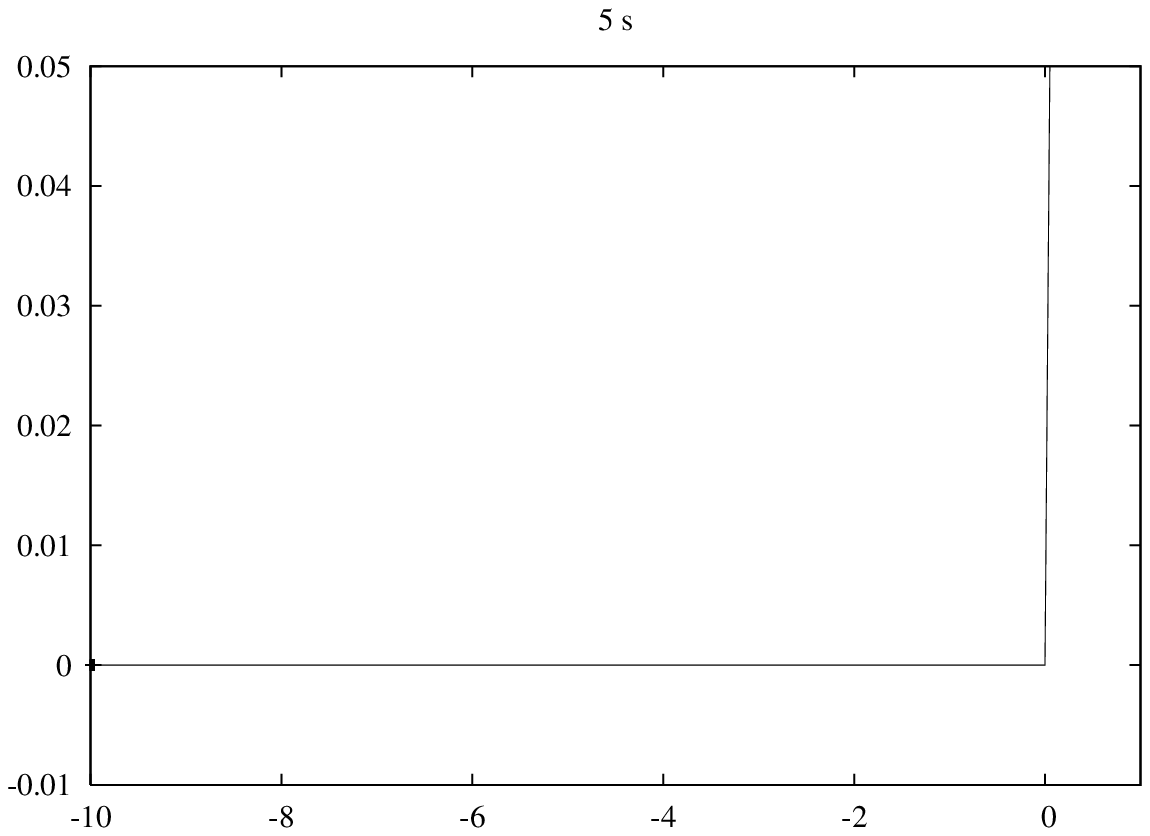}\hspace{3pt}
\includegraphics[width=6.5cm]{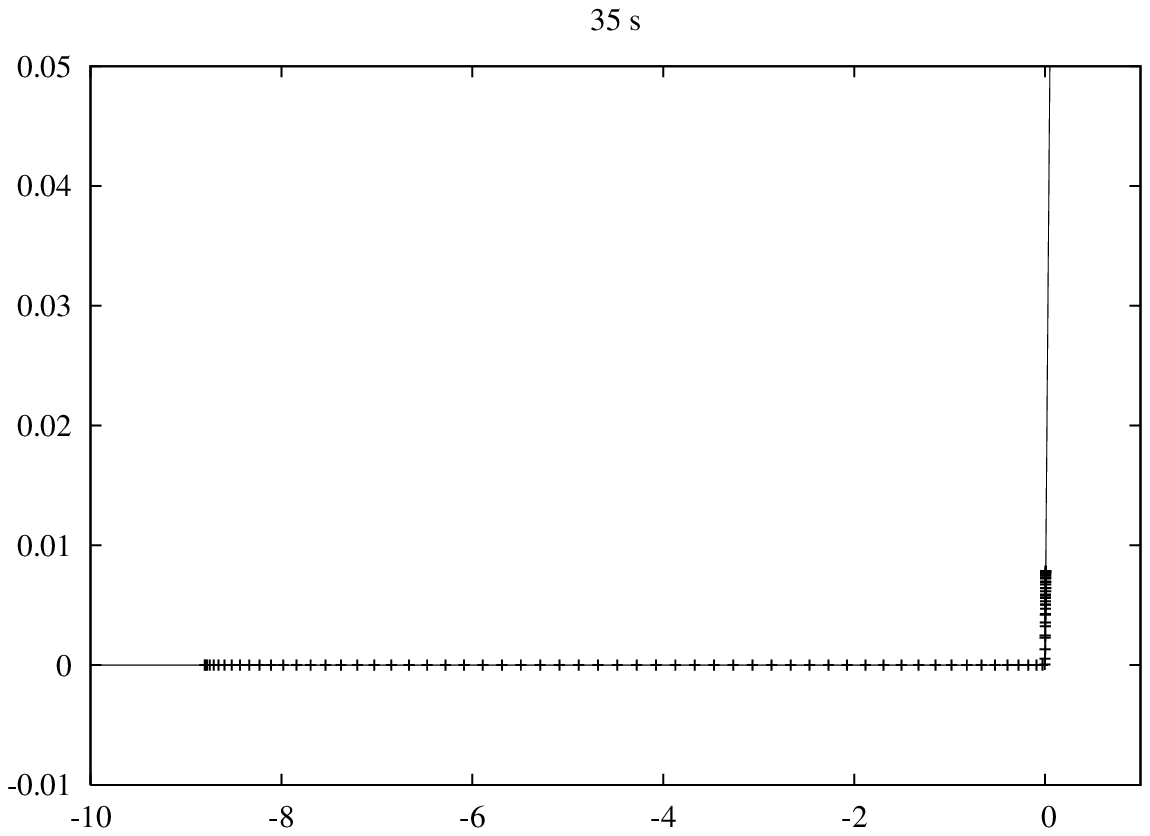}\\[6pt]
\includegraphics[width=6.5cm]{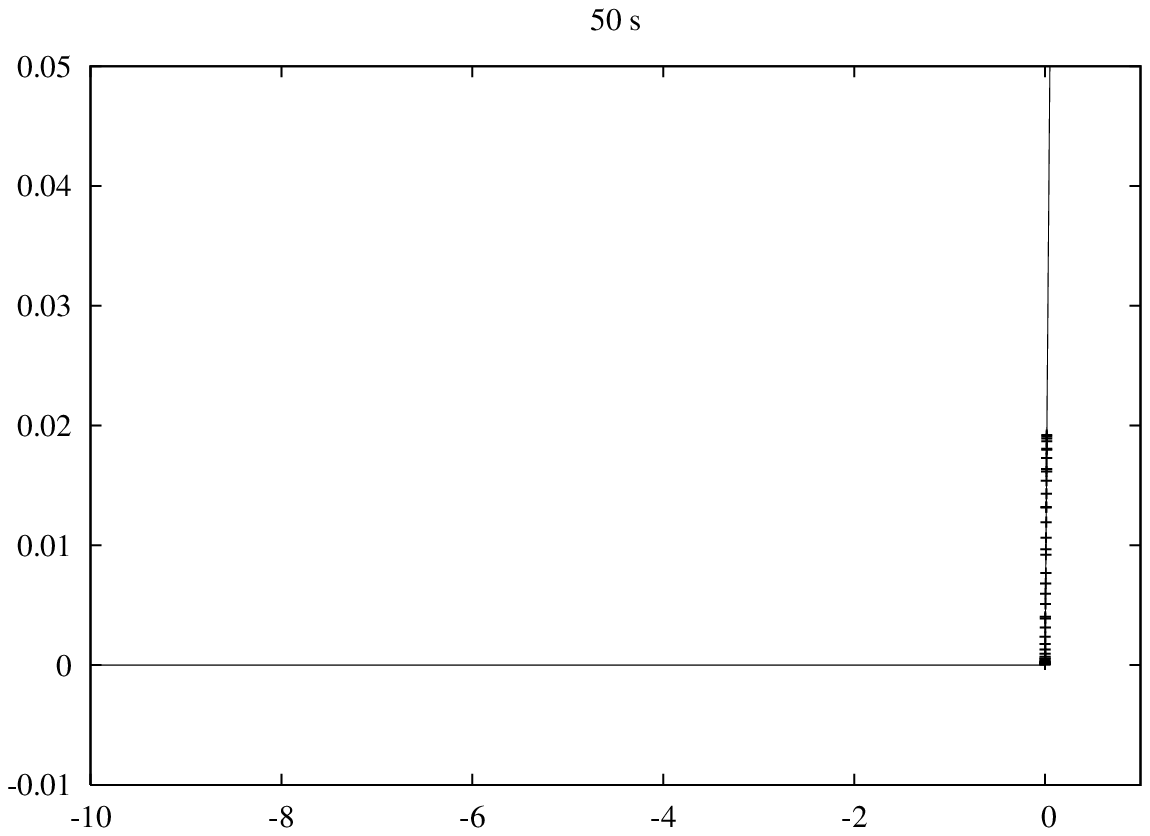}\hspace{3pt}
\includegraphics[width=6.5cm]{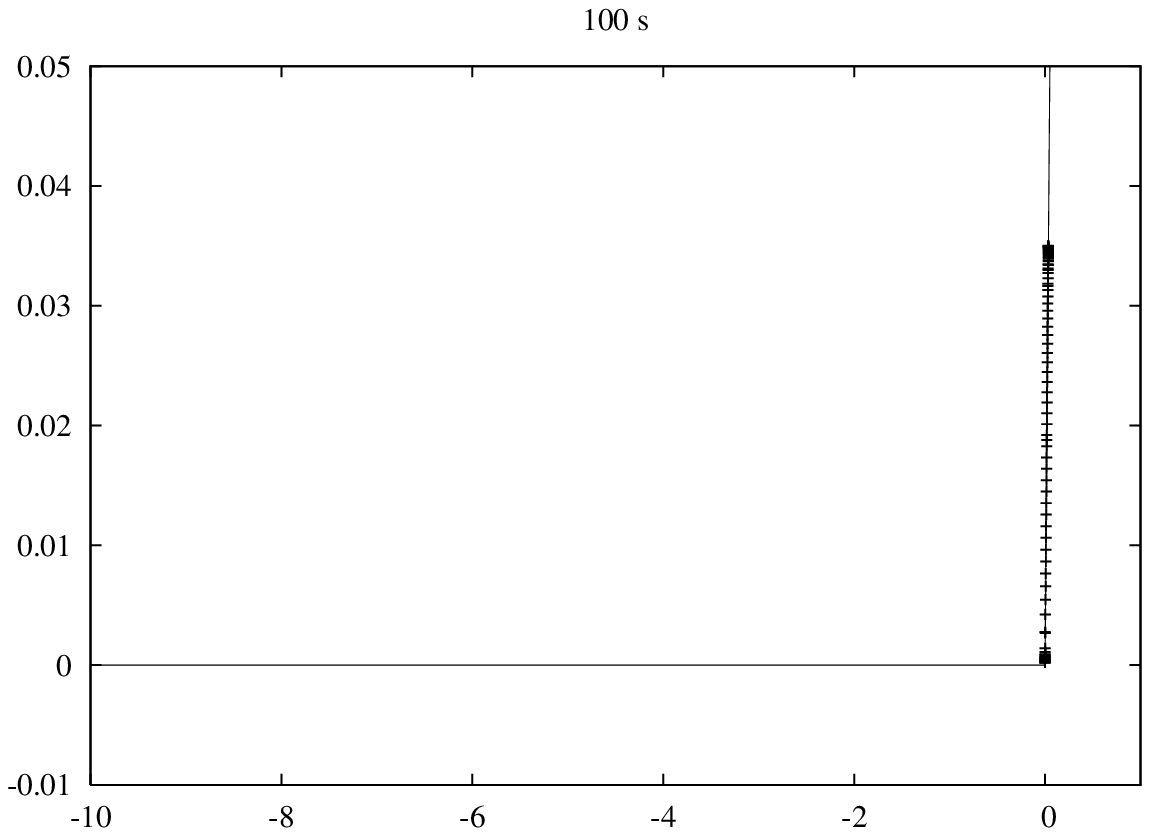}\\[6pt]
\includegraphics[width=6.5cm]{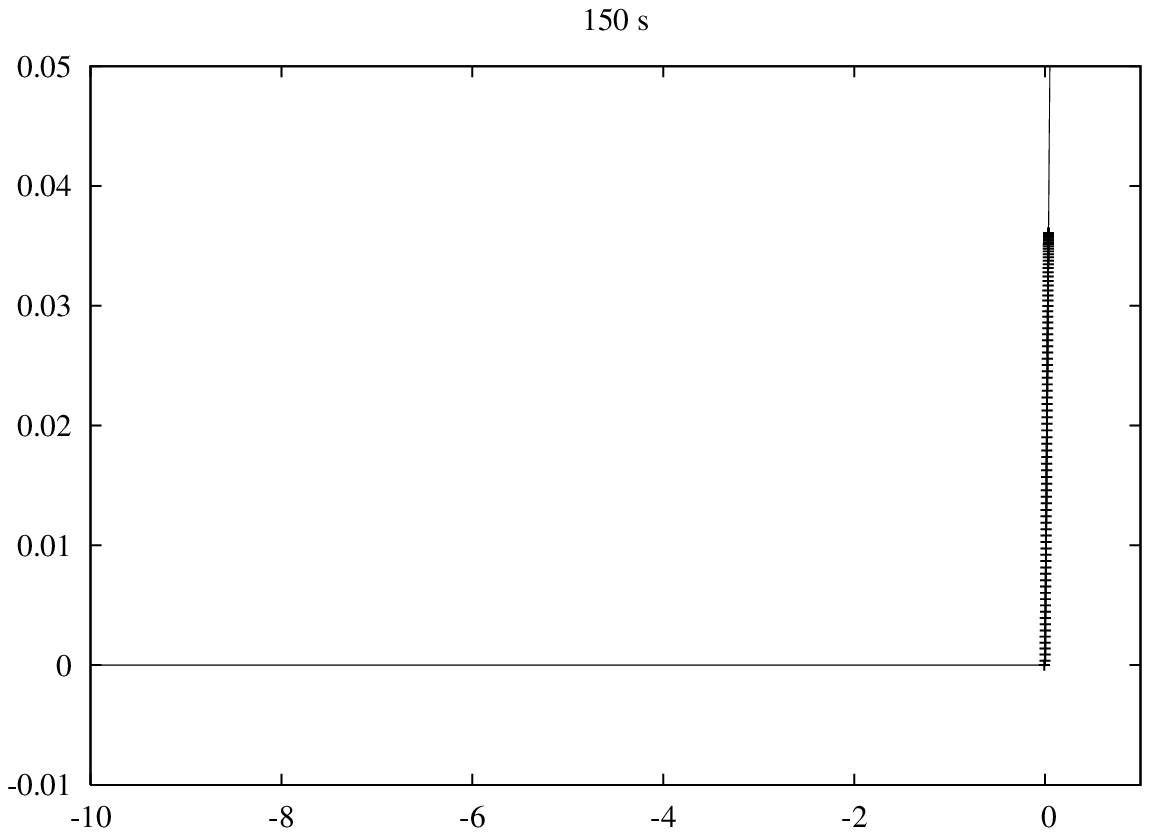}\hspace{3pt}
\includegraphics[width=6.5cm]{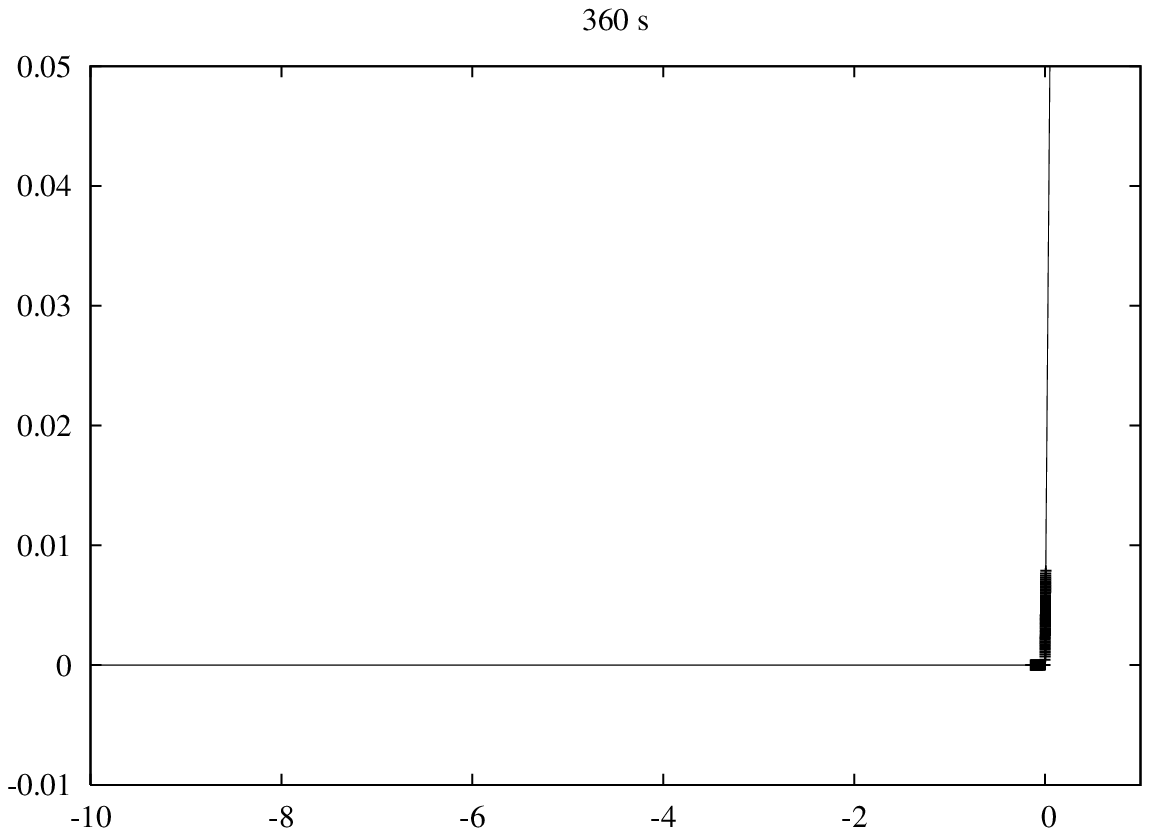}
\caption{TC3 - Cloud of points ($\psi_{\h}^n,h_{\h}^n$) on the admissible set $\A$ at different times.\label{Fig:Couplage16}}
\end{center}
\end{figure}

Figure~\ref{Fig:Couplage16} shows that each couple
($\psi_{\h}^n,h_{\h}^n$) stays on the admissible set $\A$. As
in the previous test case, the phases described above are
clearly illustrated by the position of the cloud of
points. It is located on the branch $\{h=0\}$ when the soil
is unsatured, on the branch $\{h=\psi\}$ when the soil is
saturated and on the two branches when there are both
saturated and unsaturated zones at the interface.

\begin{figure}[htb]
\begin{center}
\psfrag{Mass}{\hspace{-0.2cm}\raisebox{0.1cm}{\footnotesize{Mass}}}
\psfrag{Time}{\hspace{-0.25cm}\raisebox{-0.2cm}{\footnotesize{Time}}}
\psfrag{a1}{\hspace{-0.6cm}\raisebox{-0.1cm}{\scriptsize{$\Sigma\M^n_{\rm{in}}$}}}
\psfrag{a2}{\hspace{-0.5cm}\raisebox{0.0cm}{\scriptsize{$\DM^n$}}}
\psfrag{a3}{\hspace{-3.81cm}\raisebox{0.2cm}{\scriptsize{$\DM^n_{\rm{grnd}}$}}}
\psfrag{a4}{\hspace{-0.8cm}\raisebox{-0.1cm}{\scriptsize{$\DM^n_{\rm{over}}$}}}
\psfrag{a5}{\hspace{-0.7cm}\raisebox{-0.05cm}{\scriptsize{$\Sigma\M^n_{\rm{out}}$}}}
\includegraphics[width=6.5cm]{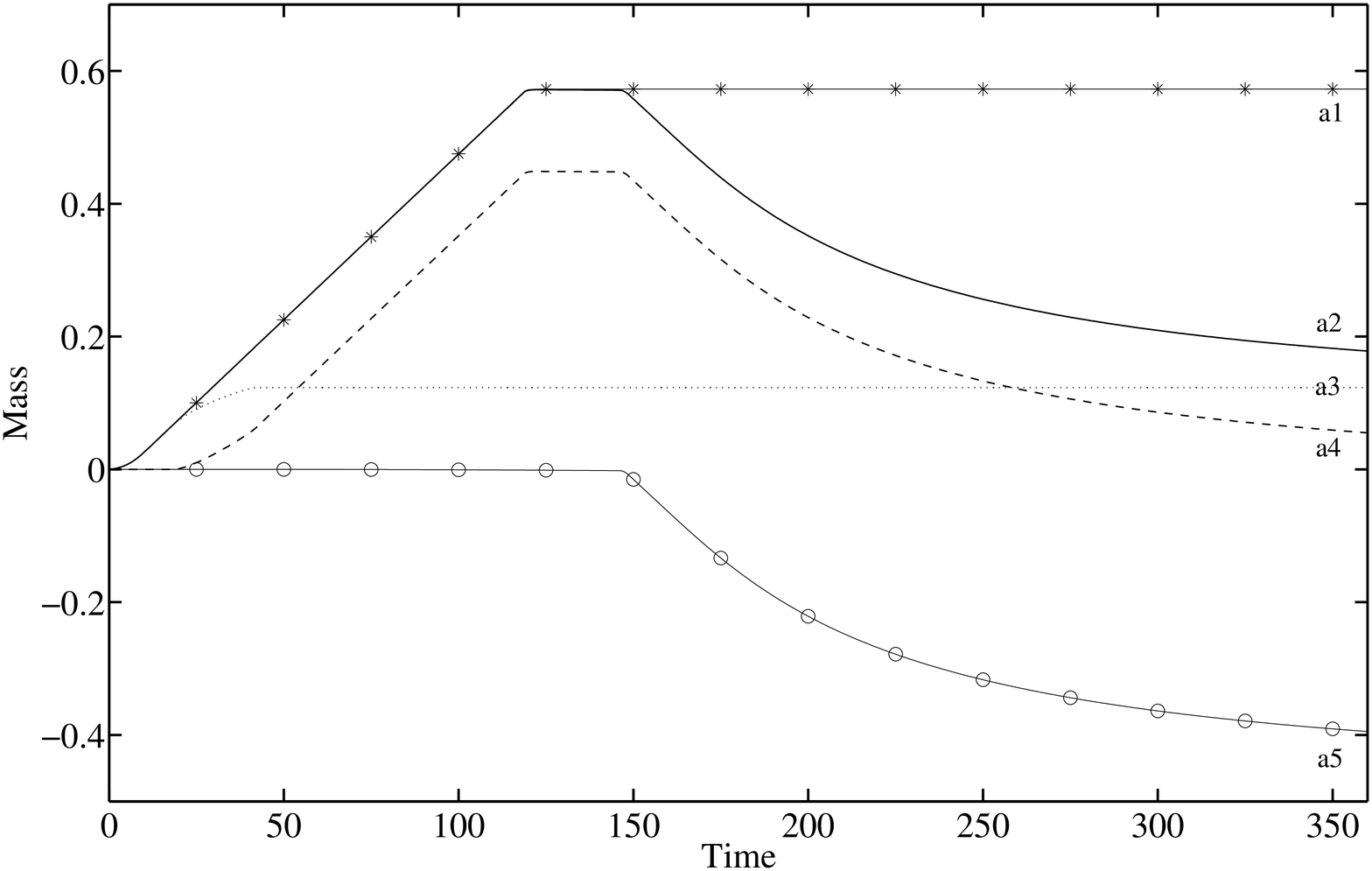}\hspace{0.2cm}
\includegraphics[width=6.5cm]{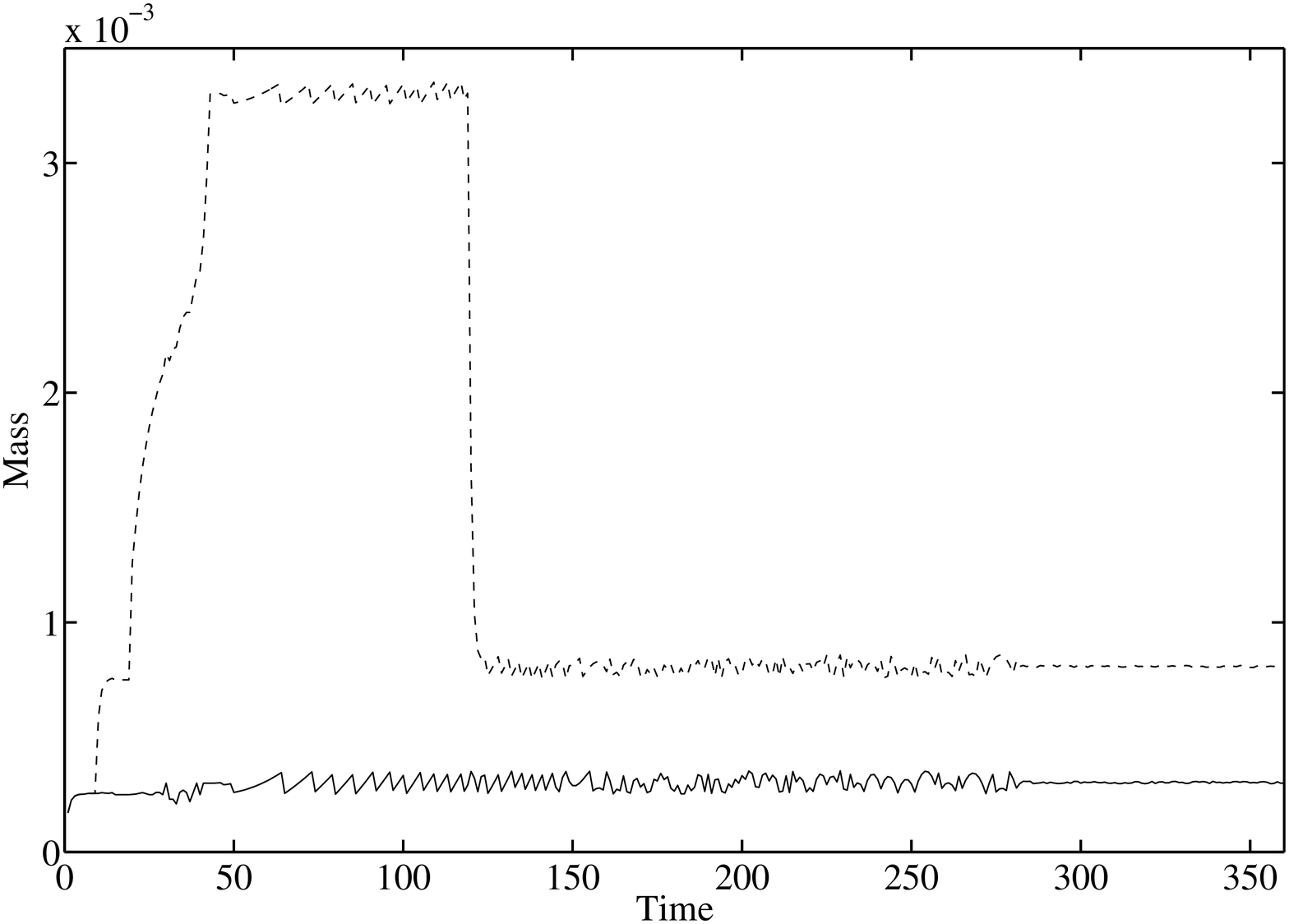}
\caption{TC3 - Left: Mass repartition in the coupled system; Right: Mass balance defect $\mathbb{E}^{n}$ for Algorithm~\ref{alg2} (dashed) and Algorithm~\ref{alg3} (solid).}\label{Fig:Couplage17}
\end{center}
\end{figure}

\begin{figure}[htb]
\begin{center}
\psfrag{Mass}{\hspace{-0.2cm}\raisebox{0.1cm}{\footnotesize{Mass}}}
\psfrag{Time}{\hspace{-0.25cm}\raisebox{-0.2cm}{\footnotesize{Time}}}
\psfrag{b1}{\hspace{0.8cm}\raisebox{-0.05cm}{\scriptsize{$\M^n_{\rm{in}}$}}}
\psfrag{b2}{\hspace{0cm}\raisebox{0cm}{\scriptsize{$\rho\delta t \tilde{F}_{\I}^{+}$}}}
\psfrag{b3}{\hspace{1.2cm}\raisebox{-0.1cm}{\scriptsize{$\rho\delta t \tilde{F}_{\I}^{-}$}}}
\psfrag{b4}{\hspace{-0.2cm}\raisebox{-0.1cm}{\scriptsize{$\M^n_{\rm{out}}$}}}
\includegraphics[width=6.5cm]{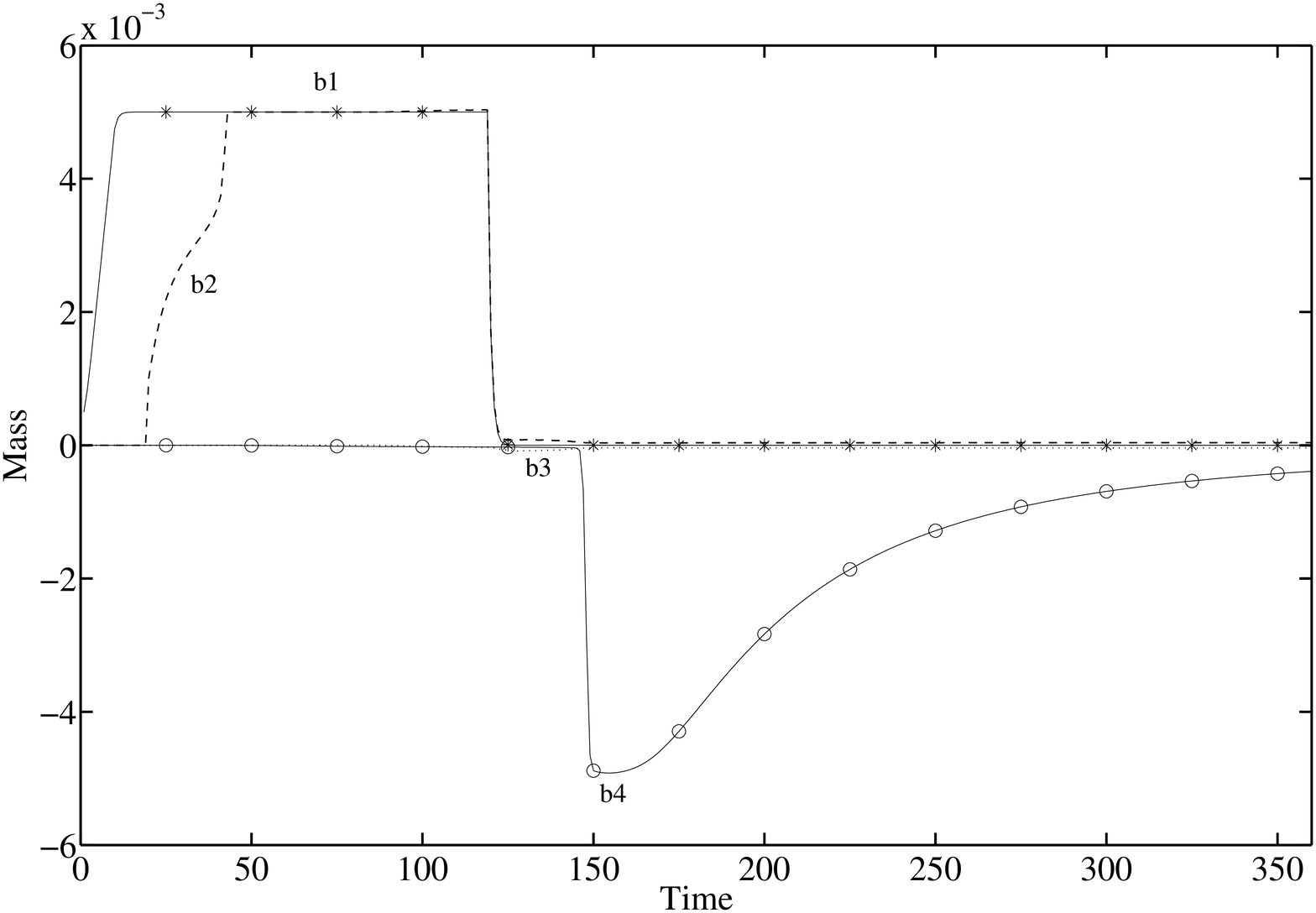}
\caption{TC3 - Mass fluxes in the kinematic wave equation.\label{Fig:Couplage18}}
\end{center}
\end{figure}

Results on Figure~\ref{Fig:Couplage17} and
Figure~\ref{Fig:Couplage18} are similar to the ones of the
previous test case, in particular the comparison of the mass
balance defects for the one-step and the two-step algorithms.  

\section{Conclusion}
In this work we have presented a robust and accurate
numerical method to simulate coupled subsurface and overland
flows governed by Richards' equation and the kinematic wave
equation. Special care was taken to design coupling
algorithms that preserve the overall mass in the system and
that also satisfy the various equality and inequality
constraints imposed at the interface. Extensions of this
work include the use of more complex models, such as the
shallow-water equations, to describe the overland flow and
the possibility to account for drainage pipes in the
soil. Extension to two-dimensional surface flows and
three-dimensional subsurface variably saturated flows can
also be considered. The present algorithms are currently
being tested in more complex and realistic test cases
related to field studies.  

\bibliographystyle{elsart-num-sort}
\bibliography{bibliographie}

\end{document}